\documentclass[10pt]{article}

\usepackage{amsmath}
\usepackage{amsmath, amsthm, amsfonts, amssymb, color}
\usepackage{mathrsfs}
\usepackage{latexsym,bm}
\usepackage{graphicx}
\usepackage[toc,page]{appendix}

\newtheorem{theorem}{Theorem}[section]%edit the theorems
\newtheorem{proposition}[theorem]{Proposition}
\newtheorem{assumption}[theorem]{Assumption}
\newtheorem{lemma}[theorem]{Lemma}
\newtheorem{corollary}[theorem]{Corollary}

\newtheorem{remark}[theorem]{Remark}

\definecolor{wco}{rgb}{0.5,0.2,0.3}

\newcommand{\pp}{\mathbb{P}}

\newcommand{\mf}{\mathcal{M}(E)}
\newcommand{\mfo}{\mathcal{M}^{0}(E)}

\newcommand{\ud}{\mathrm{d}}

\newcommand{\R}{\mathbb{R}}
\newcommand{\e}{\mathrm{e}}

\newcommand{\re}{\mathrm{Re}}
\newcommand{\im}{\mathrm{Im}}
\numberwithin{equation}{section} %the number of the equation will be counted in one section
\begin{document}

\allowdisplaybreaks

\title{\bf Fluctuations of the linear functionals for supercritical non-local branching superprocesses}
\author{
{\bf Ting Yang }\thanks{School of Mathematics and Statistics, Beijing Institute of Technology, Beijing, 100081, P.R.China.
Email: yangt@bit.edu.cn
The research of this author is supported by NSFC (Grant Nos. 12271374 and 12371143).}
}

\date{}
\maketitle

\begin{abstract}
Suppose $\{X_{t}:t\ge 0\}$ is a supercritical superprocess on a Luzin space $E$, with a non-local branching mechanism and probabilities $\mathbb{P}_{\delta_{x}}$, when initiated from a unit mass at $x\in E$. By ``supercritical", we mean that the first moment semigroup of $X_{t}$ exhibits a Perron-Frobenius type behaviour characterized by an eigentriplet $(\lambda_{1},\varphi,\widetilde{\varphi})$, where the principal eigenvalue $\lambda_{1}$ is greater than $0$. Under a second moment condition, we prove that $X_{t}$ satisfies a law of large numbers in the sense that, for any bounded measurable function $f$ on $E$,
$$\lim_{t\to+\infty}\mathrm{e}^{-\lambda_{1}t}\langle f,X_{t}\rangle =\langle f,\widetilde{\varphi}\rangle W^{\varphi}_{\infty}\quad\mbox{ in }L^{2}(\mathbb{P}_{\delta_{x}}),$$
where $W^{\varphi}_{\infty}$ is the limit of the martingale $W^{\varphi}_{t}=\mathrm{e}^{-\lambda_{1}t}\langle \varphi,X_{t}\rangle$.

The main purpose of this paper is to further investigate the fluctuations of the linear functional $\mathrm{e}^{-\lambda_{1}t}\langle f,X_{t}\rangle$ around the limit given by the law of large numbers. To this end, we introduce a parameter $\epsilon(f)$ for a bounded measurable function $f$, which determines the exponent term of the decay rate for the first moment of the fluctuation. Qualitatively, the second-order behaviour of $\langle f,X_{t}\rangle$ depends on the sign of $\epsilon(f)-\lambda_{1}/2$. We prove that, for a suitable test function $f$, the fluctuation of the associated linear functional exhibits distinct asymptotic behaviours depending on the magnitude of $\epsilon(f)$: If $\epsilon(f)\ge \lambda_{1}/2$, the fluctuation converges in distribution to a Gaussian limit under appropriate normalization; If $\epsilon(f)<\lambda_{1}/2$, the fluctuation converges to an $L^{2}$ limit with a larger normalization factor. In particular, when the test function is chosen as the right eigenfunction $\varphi$, we establish a functional central limit theorem. As an application, we consider a multitype superdiffusion in a bounded domain. For this model, we derive limit theorems for the fluctuations of arbitrary linear functionals.
\end{abstract}

\medskip


\medskip


\section{Introduction}\label{sec1}

Superprocesses, which are an important class of measure-valued processes, arise as the short lifetime and high density limits of branching Markov processes (BMPs). There is an extensive body of literature on superprocesses with the so-called local branching mechanisms, such as \cite{Dawson,Dynkin,Eth,Li}. In recent years, research has been extended to the more general framework of non-local branching mechanisms, as seen in \cite{DGL,GHK,Li,PY} and the references therein. This generalization introduces additional complexity but also allows for a richer class of models that can capture the non-local dependencies among individuals.

In this paper, we consider a superprocess $X=\{X_{t}:t\ge 0\}$ with a non-local branching mechanism on a Luzin space $E$. In a general setting where the mean semigroup of $X_{t}$ exhibits a Perron-Frobenius type behaviour with a principal eigenvalue $\lambda_{1}>0$, a positive right eigenfunction $\varphi$ and a left eigenmeasure $\widetilde{\varphi}$, we prove that, under a second moment condition, $X$ satisfies a law of large numbers in the sense that, for any bounded measurable function $f$ on $E$,
$$\lim_{t\to+\infty}\e^{-\lambda_{1}t}\langle f,X_{t}\rangle=\langle f,\widetilde{\varphi}\rangle W^{\varphi}_{\infty}\mbox{ in }L^{2}(\pp_{\delta_{x}}),$$
where $W^{\varphi}_{\infty}$ is the (non-degenerate) limit of the martingale $W^{\varphi}_{t}=\e^{-\lambda_{1}t}\langle \varphi,X_{t}\rangle$.
We are therefore interested in understanding how the linear functional $\e^{-\lambda_{1}t}\langle f,X_{t}\rangle$ fluctuates around this random limit as time increases.

This question has been widely and deeply studied for branching Markov processes (BMPs), see, e.g. \cite{KS66a,KS66b,A69a,A69b,A71,Janson} for multitype branching processes, \cite{AM} for branching Ornstein-Uhlenbeck processes, and \cite{RSZ14a,RSZ17a} for a more general class of
spatial
BMPs, where the mean operator is a Hilbert-Schmidt operator on an $L^{2}$ space. Under a second moment condition, three different regimes were found depending on the spectral properties of the mean semigroup of the process. These three regimes are referred to as the small, critical and large branching cases. In the small and critical branching cases, the fluctuation converges in distribution to a Gaussian limit under appropriate normalization, while in the large branching case, the fluctuation converges to a non-Gaussian limit, either almost surely or in $L^{2}$, with a significantly larger normalization. Very recently, such results have been further generalized in \cite{CZ} and \cite{DH} to spatial BMPs with non-local branching mechanisms, where the offspring can be displaced from the parent particle at branching events.

It is reasonable to conjecture that the aforementioned trichotomy remains valid for superprocesses. Indeed, there are several earlier results in this area, though they are limited to local branching mechanisms. For instance, \cite{M,RSZ14b} obtained the spatial central limit theorems for super-Ornstein-Uhlenbeck processes under a second moment condition, establishing the counterparts to the results in \cite{AM}. Subsequently, \cite{RSZ15,RSZ17b} generalized these results to a larger class of superprocesses, paralleling the results in \cite{RSZ14a}.
More recently, \cite{RSSZ} established the spatial central limit theorems for super-OU process with a special ``stable-like" branching mechanism. This type of branching mechanism allows an infinite second moment, thus resulting in different normalization and convergence to a stable random variable.
Additionally, we note that the recent paper \cite{LRS22} studied the convergence rate of the tail of the martingale $W^{\varphi}_{t}$ under some mild moment conditions.

It is worth highlighting the differences between the approaches for these two classes of processes.
For a BMP $\{Z_{t}:t\ge 0\}$, one can write $Z_{t+s}=\sum_{u\in N_{t}}Z^{u,t}_{s}$, where $N_{t}$ denotes the set of particles alive at time $t$, and $Z^{u,t}_{s}$ denotes the BMP starting from $u\in N_{t}$. Conditionally on $Z_{t}$, $\langle f,Z_{t+s}\rangle$ is the sum of a finite number of independent random variables. Thus, the problem essentially reduces to analyzing the limit behaviour of a sequence of sums of independent random variables. In contrast, since such a clear particle structure does not exist for superprocesses, one cannot use classical techniques, such as the Lindeberg-Feller theorem and the Berry-Esseen theorem, to prove convergence. To overcome this difficulty, \cite{M,RSZ14b} employed the method of skeleton decomposition, which represents the superprocess as an immigration process along a BMP, called the skeleton. An advantage of this method is that, it enables one to transfer results directly from the theory of BMPs to superprocesses. However, for a general superprocess with a non-local branching mechanism, even the existence of a skeleton requires justification, which is a non-trivial task. \cite{RSSZ,RSZ15,RSZ17b} used the excursion measures of superprocesses as an alternative to skeleton decomposition. In order to construct such an excursion measure, one has to impose the ``non-persistence" assumption, i.e., $\pp_{\delta_{x}}\left(X_{t}=0\right)\in (0,1)$ for all $x\in E$ and $t>0$. For spatially independent branching mechanisms, this assumption is implied by the Grey's condition, but for general non-local branching mechanisms, this assumption is difficult to verify.

As we have noted, none of the aforementioned approaches can be directly applied to our current setting of superprocesses with non-local branching mechanisms. In this paper, we take a new approach to investigate the asymptotic behaviour of fluctuations. A key component of our method is a L\'{e}vy-Khintchine type representation of the log-Laplace functional of superprocesses. This representation arises naturally from the fact that $X_{t}$ is an infinitely divisible random measure for every $t>0$. Our proof of the (finite dimensional) weak convergence of the fluctuations relies mainly on an (inductive) argument using the Markov property and this L\'{e}vy-Khintchine type representation. Another tool used in our proofs is a stochastic integral representation of the superprocesses. This representation enables one to decompose superprocesses into martingale measures, making it particularly useful in studying the structural properties of superprocesses.

Our main theorems (Theorem \ref{them1} and Theorem \ref{them2}) show that, for a suitable test function $f$, under a second moment assumption and an additional condition on the decay rate of the first moment of $\e^{-\lambda_{1}t}\langle f,X_{t}\rangle-\langle f,\widetilde{\varphi}\rangle W^{\varphi}_{\infty}$, the fluctuation given by
$$C(t)\big(\e^{-\lambda_{1}t}\langle f,X_{t}\rangle-\langle f,\widetilde{\varphi}\rangle W^{\varphi}_{\infty}\big)$$
 exhibits one of the three possible asymptotic behaviours depending on the magnitude of $\epsilon(f)$ (where $\epsilon(f)$ is a parameter defined in \eqref{eq:epsilon(f)} below):
 \begin{itemize}
 \item{If $\epsilon(f)>\lambda_{1}/2$, then the fluctuation converges in distribution to a Gaussian limit, with the normalization $C(t)=\e^{\frac{\lambda_{1}}{2}t}$.}
 \item{If $\epsilon(f)=\lambda_{1}/2$, then the limit is still Gaussian, but the normalization requires an additional polynomial term.}
 \item{If $\epsilon(f)<\lambda_{1}/2$, then the fluctuation converges to an $L^{2}$ limit, and the normalization depends mainly on the decay rate of the first moment.}
  \end{itemize}
In particular, when the test function is taken to be the right eigenfunction $\varphi$, we prove a functional central limit theorem for the tail of the martingale $W^{\varphi}_{t}$.
We further illustrate our main results by considering the model of multitype superdiffusion in a bounded domain, where the nonlocality plays a crucial role. In this setting, we obtain limit theorems (Proposition \ref{prop:b} and Proposition \ref{prop:d}) for the fluctuations of arbitrary linear functionals.

This paper is organized as follows. In the remainder of Section \ref{sec1}, we formally define the non-local branching superprocesses, present the main assumptions and results, and apply our results to the multitype superdiffusion model. Section \ref{sec2} begins with a review on the stochastic integral representation of superprocesses, followed by an analysis of the asymptotic behaviour of the linear functional
$\langle f,X_{t}\rangle$
and a proof of the law of large numbers. Section \ref{sec3} and Section \ref{sec4} contain the proofs of Theorem \ref{them1} and Theorem \ref{them2}, respectively. Section \ref{sec5} provides the proofs of the results stated in Section \ref{sec1.3}. Some auxiliary results needed in the proofs are proved in the Appendix.

\textbf{Notation:} In this paper we use ``:=" as a way of
definition.
{The letter $c$, with or without subscript, denotes a finite positive constant whose value may vary from place to place.}
$\R$ and $\mathbb{C}$ stand for the sets of real and complex numbers, respectively. All vectors in $\R^{n}$ or $\mathbb{C}^{n}$ will be understood as column vectors.
Let $E$ be a Luzin topological space with Borel $\sigma$-algebra $\mathcal{B}(E)$. Let $E_{\partial}:=E\cup\{\partial\}$ be the one-point compactification of $E$. Any function on $E$ will be automatically extended to $E_{\partial}$ by setting $f(\partial)=0$.
We use $\mathcal{B}_{b}(E)$ (resp. $\mathcal{B}^{+}(E)$) to denote the space of bounded (resp. non-negative) measurable functions on $(E,\mathcal{B}(E))$.
{We say a sequence of functions $\{f_{n}:n\ge 1\}$ on $E$ converges boundedly and pointwise to $f$ if there is a constant $c\ge 0$ such that $\|f_{n}\|_{\infty}\le c$ for $n$ sufficiently large and $f_{n}(x)\to f(x)$ as $n\to+\infty$ for all $x\in E$. We use ``$\stackrel{b}{\to}$" to denote the bounded pointwise convergence.}
Let $\mf$ denote the space of finite Borel measures on $E$ topologized by the weak convergence and $\mfo:=\mf\setminus\{0\}$ where $0$ denotes the null measure. For a measure $\mu$ on $\mathcal{B}(E)$ and a measurable function $f$, let $\langle f,\mu\rangle:=\int_{E}f(x)\mu(\ud x)$ whenever the integral makes sense.
For $p\ge 1$ and a probability space $(\Omega,\mathcal{F},\mathrm{P})$, we use $L^{p}(\mathrm{P})$ to denote the space of $\mathbb{C}$-valued random variables $\eta$ on $\Omega$ satisfying that $\int_{\Omega}|\eta|^{p}d\mathrm{P}<+\infty$.
For a $\mathbb{C}$-valued function $f$, we write $\bar{f}$ for its complex-conjugate, and define $\|f\|_{\infty}:=\sup|f(x)|$.
For any two positive functions $f$ and $g$,
we use $f\stackrel{c}\lesssim g$ (resp. $f\stackrel{c}{\asymp} g$) to denote that there is a positive constant
$c$ such that $f\le c g$ (resp. $c^{-1}g\le f\le c g$) on their common domain of definition. We also
write ``$\lesssim$" (resp. $\asymp$) if $c$ is unimportant or
understood.

\subsection{Non-local branching superprocess}\label{sec1.1}

Suppose $\xi=(\Omega,\mathcal{H},\mathcal{H}_{t},\xi_{t},\Pi_{x})$ is a Borel right process in $E$ with transition semigroup $(P_{t})_{t\ge 0}$. Here $\{\mathcal{H}_{t}:t\ge 0\}$ is the minimal admissible filtration. In this paper, we consider a superprocess $X:=(X_{t})_{t\ge 0}$ with spatial motion $\xi$ and a non-local branching mechanism $\psi$ given by
\begin{align}
\psi(x,f)&=a(x)f(x)+b(x)f(x)^{2}-{\int_{E}f(y)\eta(x,\ud y)}\nonumber\\
&\quad+\int_{\mfo }\left(\e^{-{\langle f,\nu\rangle}}-1+\nu(\{x\})f(x)\right)H(x,{ \ud \nu})\quad\forall x\in E,\ f\in\mathcal{B}^{+}(E),\label{def:branching mechanism}
\end{align}
where $a\in\mathcal{B}_{b}(E)$, $b\in\mathcal{B}^{+}_{b}(E)$, $\eta(x,{ \ud y})$ is a bounded kernel on $E$, $\nu_{x}({ \ud y})$ denotes the restriction of $\nu({ \ud y})$ to $E\setminus \{x\}$, and $\left({\langle 1,\nu\rangle}\wedge {\langle 1,\nu\rangle}^{2}+{\langle 1,\nu_{x}\rangle}\right)H(x,{ \ud \nu})$ is a bounded kernel from $E$ to $\mfo$.
To be specific, $X$ is a $\mf$-valued Markov process satisfying that for every $f\in\mathcal{B}^{+}_{b}(E)$ and $\mu\in\mf$,
\begin{equation}\label{eq1}
\pp_{\mu}\left[\e^{-\langle f,X_{t}\rangle}\right]=\e^{-\langle V_{t}f,\mu\rangle},
\end{equation}
where $V_{t}f(x):=-\log\pp_{\delta_{x}}\left[\e^{-\langle f,X_{t}\rangle}\right]$ is the unique nonnegative locally bounded solution to the integral equation
\begin{equation}\label{eq:integral equation for V_{t}}
V_{t}f(x)=P_{t}f(x)-\Pi_{x}\left[\int_{0}^{t}\psi(\xi_{s},V_{t-s}f){ \ud s}\right].
\end{equation}
By ``locally bounded" we mean that $\sup_{x\in E,t\in [0,T]}|V_{t}f(x)|<+\infty$ for every $T\in (0,+\infty)$.
Such a process is defined in \cite{Li} via its log-Laplace functional and is referred to as the $(P_{t},\psi)$-superprocess.

The branching mechanism given in \eqref{def:branching mechanism} is quite general.
For example, let
\begin{equation}\label{def:purely-local}
\psi^{L}(x,\lambda):=\hat{a}(x)\lambda+b(x)\lambda^{2}+\int_{(0,+\infty)}\left(e^{-\lambda u}-1+\lambda u\right)\pi(x,\ud u),
\end{equation}
for $x\in E$ and $\lambda\ge 0$,  and
\begin{equation*}
\psi^{NL}(x,f):=-{\int_{E}f(y)\eta(x,\ud y)}+\int_{\mfo }\left(\e^{-{\langle f,\nu\rangle}}-1\right)\Gamma(x ,\ud \nu)
\end{equation*}
for $x\in E$ and $f\in\mathcal{B}^{+}(E)$, where $b,\eta$ are as in \eqref{def:branching mechanism}, $\hat{a}\in\mathcal{B}_{b}(E)$,  $(u\wedge u^{2})\pi(x,\ud u)$ is a bounded kernel from $E$ to $(0,+\infty)$, and ${\langle 1,\nu\rangle}\Gamma(x,\ud \nu)$ is a bounded kernel from $E$ to $\mfo $. Then $(x,f)\mapsto \psi^{L}(x,f(x))+\psi^{NL}(x,f)$ is a branching mechanism that can be represented in the form of \eqref{def:branching mechanism} with $a(x)=\hat{a}(x)-\int_{\mfo }\nu(\{x\}) H(x,{ \ud \nu})$ and $H(x,{ \ud \nu})= \Gamma(x,{ \ud \nu})+\int_{(0,+\infty)}\delta_{u\delta_{x}}({ \ud \nu})\pi(x,{ \ud u})$.
A branching mechanism of this type is said to be decomposable with local part $\psi^{L}$ and non-local part $\psi^{NL}$.
In particular, if the non-local part equals $0$, we call such a branching mechanism \textit{purely local}. We emphasize that the branching mechanism considered in this paper is allowed to be non-local and non-decomposable.

By \cite[Theorem 5.12]{Li}, a $(P_{t},\psi)$-superprocess $X$ has a Borel right realization in $\mf$. We denote by $\mathcal{W}^{+}_{0}$ the space of c\`{a}dl\`{a}g paths from $[0,+\infty)$ to $\mf$ having zero as a trap. Here, we assume that $X$ is the coordinate process in $\mathcal{W}^{+}_{0}$ and  $\mathcal{F}:=\{\mathcal{F}_{t}:\ t\in[0,\infty]\}$ is the filtration generated by the coordinate process, which is completed  with the class of $\pp_{\mu}$-negligible measurable sets for every $\mu\in\mf$.

It is known by \cite[Proposition 2.27]{Li} that for every $\mu\in\mf$ and $f\in\mathcal{B}_{b}(E)$,
$$\pp_{\mu}\left[\langle f,X_{t}\rangle\right]=\langle T_{t}f,\mu\rangle,$$
where $T_{t}f(x)$ is the unique locally bounded solution to the integral equation
\begin{equation}\label{eq:many-to-one}
T_{t}f(x)=P_{t}f(x)-\Pi_{x}\left[\int_{0}^{t}a(\xi_{s})T_{t-s}f(\xi_{s}){ \ud s}\right]+\Pi_{x}\left[\int_{0}^{t}m[T_{t-s}f](\xi_{s}){ \ud s}\right],
\end{equation}
where $m[f](y):={\int_{E}f(u)\eta(y,\ud u)}+\int_{\mfo }{\langle f,\nu_{y}\rangle}H(y,{ \ud \nu})$ for $y\in E$.

\subsection{Assumptions and main results}\label{sec1.2}

Throughout this paper, we assume the non-local branching superprocess satisfies the following conditions.

\begin{assumption}\label{AS1}
There exist an eigenvalue $\lambda_{1}>0$, a corresponding positive
bounded function
$\varphi$ and a finite left eigenmeasure $\widetilde{\varphi}$, with normalization $\langle \varphi,\widetilde{\varphi}\rangle=1$, such that for every $f\in\mathcal{B}^{+}_{b}(E)$, $\mu\in\mf$ and $t\ge 0$,
    $$\langle T_{t}\varphi,\mu\rangle=\e^{\lambda_{1} t}\langle \varphi,\mu\rangle\mbox{ and }\langle T_{t}f,\widetilde{\varphi}\rangle=\e^{\lambda_{1} t}\langle f,\widetilde{\varphi}\rangle.$$ Further, we define
$$\triangle_{t}:=\sup_{x\in E,f\in\mathcal{B}^{+}_{1}(E)}\left|\varphi(x)^{-1}\e^{-\lambda_{1} t}T_{t}f(x)-\langle f,\widetilde{\varphi}\rangle\right|\quad\forall t\ge 0.$$
Here $\mathcal{B}^{+}_{1}(E)$ denotes the set of nonnegative measurable functions on $E$ that are bounded by unity.
We assume that $\sup_{t\ge 0}\triangle_{t}<+\infty$ and $\lim_{t\to+\infty}\triangle_{t}=0$.
\end{assumption}

\begin{assumption}\label{AS2}
The operators $f\mapsto\psi(\cdot,f)$ and $f\mapsto -a f+m[f]$ preserve $C^{\xi}_{b}(E)$. Here $C^{\xi}_{b}(E)$ denotes the space of bounded measurable functions on $E$ that are finely continuous with respect to $\xi$,
i.e.,

$C^{\xi}_{b}(E):=\{f\in\mathcal{B}_{b}(E):\ \mbox{the mapping }t\mapsto f(\xi_{t})\mbox{ is a.s. right continuous on }[0,+\infty)\}.$
\end{assumption}

\begin{assumption}\label{AS3}
$\sup_{x\in E}\int_{\mfo }{\langle 1,\nu\rangle}^{2}H(x,{ \ud \nu})<+\infty$.
\end{assumption}

We give some remarks on the assumptions here.
Assumption \ref{AS1} is a Perron-Frobenius type assumption that ensures the existence of the leading eigenvalue and corresponding eigenfunction. The superprocess is called {\it supercritical} when $\lambda_{1}>0$. Without restriction on the sign of $\lambda_{1}$, Assumption \ref{AS1} has recently been referred to as the Asmussen-Hering class of branching processes in \cite{HK}, owing to the foundational results for this class in \cite{AH}.
Assumption \ref{AS1} is a mild assumption on the semigroup convergence, and can be verified for a broad class of superprocesses. For instance, it applies to multitype continuous-state branching processes (\cite[Example 3.7]{PY}), multitype superdiffusions in bounded domains (\cite[Example 3.8]{PY}), and a particular class of super-stable processes in bounded domains (\cite[Example7.3]{RSY}). More examples satisfying Assumption \ref{AS1} can be found in \cite{GHK,LRS22,PY}.

Assumption \ref{AS2} is a regularity condition on the branching mechanism, which is inherited from \cite{PY}. This assumption is only required to ensure the existence of a stochastic integral representation for the non-local branching superprocess (see, Section \ref{sec:stochastic integral} below). Thus, Assumption \ref{AS2} is not indispensable and can be replaced by other conditions that yield the same result, such as Conditions 7.1 and 7.2 of \cite[Chapter 7]{Li}.

Assumption \ref{AS3} is obviously a second moment condition. Under Assumption \ref{AS3}, we define for any $\mathbb{C}$-valued bounded functions $f$ and $g$,
$$\vartheta[f,g](x):=2b(x)f(x)g(x)+\int_{\mfo }{\langle f,\nu\rangle}{\langle g,\nu\rangle}H(x,{ \ud \nu})\mbox{ and }\vartheta[f](x):=\vartheta[f,f](x)\quad\forall x\in E.$$
It follows by \cite[Proposition 2.38]{Li} that for every $t\ge 0$, $\mu\in\mf$ and $f\in\mathcal{B}_{b}(E)$, $\pp_{\mu}\left[\langle f,X_{t}\rangle^{2}\right]$ is finite, and
\begin{equation}\label{eq:second moment}
\pp_{\mu}\left[\langle f,X_{t}\rangle^{2}\right]=\langle T_{t}f,\mu\rangle^{2}+\langle\int_{0}^{t}T_{t-s}(\vartheta[T_{s}f]){ \ud s},\mu\rangle.
\end{equation}
The above identity yields in particular that
\begin{equation}\label{eq:variance}
\mathrm{Var}_{\delta_{x}}\left[\langle f,X_{t}\rangle\right]=\int_{0}^{t}T_{t-s}(\vartheta[T_{s}f])(x){ \ud s}\mbox{ and }\mathrm{Var}_{\mu}\left[\langle f,X_{t}\rangle\right]=\langle \mathrm{Var}_{\delta_{\cdot}}\left[\langle f,X_{t}\rangle\right],\mu\rangle.
\end{equation}

Using the fact that $\varphi$ is the bounded
eigenfunction for the superprocess with corresponding eigenvalue $\lambda_{1}$, one can easily show that
$${ W^{\varphi}_{t}}:=\e^{-\lambda_{1} t}\langle \varphi,X_{t}\rangle \quad\forall t\ge 0$$
is a nonnegative
$\pp_{\mu}$-martingale with respect to the filtration $(\mathcal{F}_{t})_{t\ge 0}$ for every $\mu\in\mf$.
{In particular, $W^{\varphi}_{t}$ is a c\`{a}dl\`{a}g martingale if $\varphi$ is continuous.}
Let ${ W^{\varphi}_{\infty}}$ be the martingale limit.
Define
$$\Theta(x):=\int_{0}^{+\infty}\e^{-2\lambda_{1} s}T_{s}(\vartheta[\varphi])(x){ \ud s}\quad\forall x\in E.$$
By \eqref{eq1.9} below, one can easily verify that $\Theta$ is a nonnegative bounded function on $E$ and $\Theta(x)\le c\varphi(x)$ for some positive constant $c$ independent of $x$.
By \eqref{eq:variance}, we have
$$\mathrm{Var}_{\delta_{x}}[{ W^{\varphi}_{t}}]
=\int_{0}^{t}\e^{-2\lambda_{1} s}T_{s}(\vartheta[\varphi])(x){ \ud s}\mbox{ and }\mathrm{Var}_{\mu}[{ W^{\varphi}_{t}}]=\langle \mathrm{Var}_{\delta_{\cdot}}[{ W^{\varphi}_{t}}],\mu\rangle.$$
Obviously, as $t\to+\infty$, $\mathrm{Var}_{\delta_{x}}[{ W^{\varphi}_{t}}]\uparrow \Theta(x)$ and $\mathrm{Var}_{\mu}[{ W^{\varphi}_{t}}]\uparrow \langle \Theta,\mu\rangle$.
This implies that $({ W^{\varphi}_{t}})_{t\ge 0}$ is an $L^{2}(\pp_{\mu})$-bounded martingale, and thus converges to ${ W^{\varphi}_{\infty}}$ in $L^{2}(\pp_{\mu})$.
Immediately we have
\begin{equation}\label{eq:second moment of Wvarphi}
\pp_{\mu}\left[{ W^{\varphi}_{\infty}}\right]=\langle \varphi,\mu\rangle\mbox{ and }\mathrm{Var}_{\mu}[{ W^{\varphi}_{\infty}}]=\langle\Theta,\mu\rangle.
\end{equation}

We first state a law of large numbers for the non-local branching superprocess.
\begin{theorem}[Law of large numbers]\label{cor:wlln}
Suppose Assumptions \ref{AS1}-\ref{AS3} hold.
For every $f\in\mathcal{B}_{b}(E)$ and $\mu\in\mf$,
$$\lim_{t\to+\infty}\e^{-\lambda_{1} t}\langle f,X_{t}\rangle = \langle f,\widetilde{\varphi}\rangle W^{\varphi}_{\infty}\mbox{ in }L^{2}(\pp_{\mu}).$$
In particular, for every $\sigma>0$, $\lim_{\mathbb{N}\ni n\to+\infty}\e^{-\lambda_{1} n\sigma}\langle f,X_{n\sigma}\rangle=\langle f,\widetilde{\varphi}\rangle W^{\varphi}_{\infty}$ $\pp_{\mu}$-a.s.
\end{theorem}

{
Results of this type have been established under different assumptions for superprocesses with purely local branching mechanisms, see, for example, \cite{CRW08,CRY19,EKW,ET02,EW06,LRS13} and the references therein, and have recently been extended to the more general framework of non-local branching mechanisms in \cite{PY}. Specifically, under an intrinsic ultracontractivity assumption (slightly stronger than Assumption \ref{AS1}) and additional continuity conditions on the mean semigroup, \cite{PY} obtained almost sure convergence as well as $L^{p}$-convergence for $p\in (1,2]$. The approach adopted in this paper is similar to that of \cite{LRS13} and \cite{PY}, relying on a stochastic integral representation of the superprocess.
We believe that, under further assumptions on the test functions or the process, the almost sure convergence along lattice times can be extended to continuous time. We refer the reader to \cite{PY} for such an argument. We do not pursue this in the present paper.
}

 The aim of this paper is to further study the fluctuations of the process around the random limit provided by the law of large numbers.

 We begin by recalling some definitions on the convergence of random variables. Suppose $S$ is a Polish space with Borel $\sigma$-field $\mathcal{S}$. We denote by $\mathcal{P}(S)$ the space of all probability measures on $(S,\mathcal{S})$ endowed with the weak topology. Then $\mathcal{P}(S)$ itself is a Polish space for this topology. We say $\rho_{n}$ converges weakly to $\rho$ in $\mathcal{P}(S)$ if and only if $\int_{S}f d\rho_{n}\to \int_{S}f d\rho$ for all bounded continuous functions $f:S\to \R$.
If $\{\xi^{n}:n\ge 1\}$ is a sequence of $S$-valued random variables,
we say $\xi^{n}$ converges in distribution to $\xi$, and write ``$\xi^{n}\stackrel{d}{\to}\xi$ in $S$", if the law of $\xi^{n}$ converges weakly to the law of $\xi$ in $\mathcal{P}(S)$.

Now we consider $\R$-valued c\`{a}dl\`{a}g processes.
Let $\mathcal{D}[0,+\infty)$ denote the space of c\`{a}dl\`{a}g functions from $[0,+\infty)$ to $\R$ equipped with the Skorokhod topology, and let $\mathcal{P}(\mathcal{D}[0,+\infty))$ be the space of all probability measures on $\mathcal{D}[0,+\infty)$ endowed with the weak topology.
An $\R$-valued c\`{a}dl\`{a}g process $Y$ may be considered as a random variable taking values in $\mathcal{D}[0,+\infty)$, and consequently its law $\mathcal{L}(Y)$ is an element of $\mathcal{P}(\mathcal{D}[0,+\infty))$. Suppose $Y^{n}$, $Y$ are c\`{a}dl\`{a}g processes. We say $Y^{n}$ converges in distribution to $Y$ in $\mathcal{D}[0,+\infty)$, and write ``$Y^{n}\stackrel{d}{\to}Y$ in $\mathcal{D}[0,+\infty)$", if $\mathcal{L}(Y^{n})$ converges weakly to $\mathcal{L}(Y)$ in $\mathcal{P}(\mathcal{D}[0,+\infty))$.

We are now ready to state the
first
main result.
\begin{theorem}[Martingale functional central limit theorem]\label{them1} Suppose Assumptions \ref{AS1}-\ref{AS3} hold.
For $t,s\ge 0$, define
$$Y^{t}_{s}:=\e^{\frac{\lambda_{1}}{2}(t+s)}\left({ W^{\varphi}_{t+s}}-{ W^{\varphi}_{\infty}}\right)
\mbox{ and }Y^{t}:=(Y^{t}_{s})_{s\ge 0}.$$
Then for every $x\in E$, under $\pp_{\delta_{x}}$, as $t\to+\infty$,
$$Y^{t}\stackrel{d}{\to}\sigma_{\varphi}\sqrt{{ W^{\varphi}_{\infty}}}G\mbox{ in the sense of finite dimensional distributions.}$$
Moreover, if, in addition, $\varphi$ is continuous, then
$$Y^{t}\stackrel{d}{\to }\sigma_{\varphi}\sqrt{{ W^{\varphi}_{\infty}}}G\mbox{ in }\mathcal{D}[0,+\infty)\mbox{ as }t\to+\infty.$$
Here $\sigma^{2}_{\varphi}:=\langle \vartheta[\varphi],\widetilde{\varphi}\rangle/\lambda_{1}$, $G:=(G_{s})_{s\ge 0}$ is a continuous Gaussian process with mean $0$ and covariance
$$Cov(G_{s},G_{t})=\e^{-\frac{\lambda_{1}}{2}|t-s|}\quad\forall s,t\ge 0,$$
and $G$ is independent of ${ W^{\varphi}_{\infty}}$.
\end{theorem}

\begin{remark}
Theorem \ref{them1} improves on an earlier result by Liu et al. \cite[Theorem 1.5]{LRS22}. According to \cite[Theorem 1.5]{LRS22}, under a second moment condition and an assumption (\cite[Assumption \ref{AS2}]{LRS22}) that implies Assumption \ref{AS1},
one has $W^{\varphi}_{t}-W^{\varphi}_{\infty}=o(\e^{-\frac{\lambda_{1}}{q}t})$ for $q>2$. Theorem \ref{them1} complements this result by giving the exact convergence rate via a functional central limit theorem. Besides, we allow for non-local branching mechanisms, whereas \cite{LRS22} considers only purely local branching mechanisms.
\end{remark}

Theorem \ref{them1} can be viewed as a statement on how fast the linear functional $\e^{-\lambda_{1}t}\langle \varphi,X_{t}\rangle$ approaches its limit $W^{\varphi}_{\infty}$. In Theorem \ref{them2} below, we provide the convergence rate of $\e^{-\lambda_{1}t}\langle f,X_{t}\rangle$ to its limit for a more general test function $f$.
We now introduce a key parameter that will appear in our main result.
For $f\in\mathcal{B}_{b}(E)$, define
\begin{equation}\label{eq:epsilon(f)}
\epsilon(f):=-\limsup_{t\to+\infty}\frac{\log\sup_{x\in E}\left|\varphi(x)^{-1}\e^{-\lambda_{1} t}T_{t}f(x)-\langle f,\widetilde{\varphi}\rangle\right|}{t}.
\end{equation}
By \eqref{ieq:H3} below and the fact that $\sup_{t\ge 0}\triangle_{t}<+\infty$, we have $\epsilon(f)\in [0,+\infty]$ for every $f\in\mathcal{B}_{b}(E)$. In particular, $\epsilon(\varphi)=+\infty$. In what follows, we use ``$\rightrightarrows$" to denote the uniform convergence on $E$.

\begin{theorem}\label{them2}
Suppose Assumptions \ref{AS1}-\ref{AS3} hold. Suppose $f\in\mathcal{B}_{b}(E)$ and $N$ is a standard normal random variable independent of $W^{\varphi}_{\infty}$. Let $\widehat{f}:=f-\langle f,\widetilde{\varphi}\rangle
\varphi$. Then the following are true.
\begin{description}
\item{(i)} If $\epsilon(f)>\lambda_{1}/2$, then for every $x\in E$, under $\pp_{\delta_{x}}$,
$$\e^{\frac{\lambda_{1}}{2}t}\left(\e^{-\lambda_{1} t}\langle f,X_{t}\rangle-\langle f,\widetilde{\varphi}\rangle W^{\varphi}_{\infty}\right)\stackrel{d}{\to}\rho_{f} \sqrt{W^{\varphi}_{\infty}} N\mbox{ as }t\to+\infty,$$
where $\rho^{2}_{f}:=\int_{0}^{+\infty}\e^{-\lambda_{1} s}\langle \vartheta[T_{s}\widehat{f}],\widetilde{\varphi}\rangle { \ud s}+\langle f,\widetilde{\varphi}\rangle^{2}\sigma^{2}_{\varphi}$.

\item{(ii)} Assume that $\pp_{\delta_{x}}\left(X_{t}=0\mbox{ for some }t\in (0,+\infty)\right)>0$. If there are $f^{*}\in\mathcal{B}_{b}(E)$ and $r\ge 0$ such that
\begin{equation}\label{them2.1}
\frac{\e^{\frac{\lambda_{1}}{2}t}}{t^{r}}\left(\e^{-\lambda_{1} t}T_{t}f-\langle f,\widetilde{\varphi}\rangle \varphi\right)\rightrightarrows f^{*}\mbox{ as }t\to+\infty,
\end{equation}
then under $\pp_{\delta_{x}}$,
$$\frac{1}{\sqrt{t}}\e^{-\frac{\lambda_{1}}{2}t}\langle f^{*},X_{t}\rangle\stackrel{d}{\to}\varrho_{f^{*}}\sqrt{W^{\varphi}_{\infty}}N\mbox{ as }t\to+\infty,$$
where $\varrho^{2}_{f^{*}}:=\langle\vartheta[f^{*}],\widetilde{\varphi}\rangle$.
In particular, if \eqref{them2.1} holds for $r=0$, then
$$\frac{1}{\sqrt{t}}\e^{\frac{\lambda_{1}}{2}t}\left(\e^{-\lambda_{1} t}\langle f,X_{t}\rangle-\langle f,\widetilde{\varphi}\rangle W^{\varphi}_{\infty}\right)\stackrel{d}{\to}\varrho_{f^{*}}\sqrt{W^{\varphi}_{\infty}}N\mbox{ as }t\to+\infty.$$
Assume, in addition, that ${\langle 1,\nu\rangle}^{4}H(y,{ \ud \nu})$ is a bounded kernel from $E$ to $\mfo$.
Then under $\pp_{\delta_{x}}$,
$$\frac{1}{t^{\frac{1}{2}+r}}\e^{\frac{\lambda_{1}}{2}t}\left(\e^{-\lambda_{1} t}\langle f,X_{t}\rangle-\langle f,\widetilde{\varphi}\rangle W^{\varphi}_{\infty}\right)\stackrel{d}{\to}\frac{1}{\sqrt{1+2r}}\varrho_{f^{*}}\sqrt{W^{\varphi}_{\infty}}N\mbox{ as }t\to+\infty.$$

\item{(iii)} If there are $f^{*}\in\mathcal{B}_{b}(E)$, $r\ge 0$ and $\epsilon\in [0,\lambda_{1}/2)$, such that
$$\frac{\e^{\epsilon t}}{t^{r}}\left(\e^{-\lambda_{1} t}T_{t}f-\langle f,\widetilde{\varphi}\rangle \varphi\right)\rightrightarrows f^{*}\mbox{ as }t\to+\infty,$$
then for every $x\in E$, $ W^{f^{*}}_{t}:=\e^{-(\lambda_{1}-\epsilon)t}\langle f^{*},X_{t}\rangle$ is an $L^{2}(\pp_{\delta_{x}})$-bounded $\mathcal{F}_{t}$-martingale, and
$$\frac{\e^{\epsilon t}}{t^{r}}\left(\e^{-\lambda_{1} t}\langle f,X_{t}\rangle-\langle f,\widetilde{\varphi}\rangle W^{\varphi}_{\infty}\right)\to W^{f^{*}}_{\infty}\mbox{ in }L^{2}(\pp_{\delta_{x}})\mbox{ as }t\to+\infty,$$
where $ W^{f^{*}}_{\infty}$ is the limit of $\{ W^{f^{*}}_{t}:t\ge 0\}$. The random variable $ W^{f^{*}}_{\infty}$ is non-degenerate under $\pp_{\delta_{x}}$
if
$\varrho^{2}_{f^{*}}=\langle \vartheta[f^{*}],\widetilde{\varphi}\rangle>0$.
Moreover, under $\pp_{\delta_{x}}$,
$$\e^{(\frac{\lambda_{1}}{2}-\epsilon)t}(
W^{f^{*}}_{t}- W^{f^{*}}_{\infty}
)\stackrel{d}{\to}\frac{1}{\sqrt{\lambda_{1}-2\epsilon}}\varrho_{f^{*}}\sqrt{W^{\varphi}_{\infty}}N
\mbox{ as }t\to+\infty.$$
\end{description}
\end{theorem}

\begin{remark}
(i) We note that the limit function $f^{*}$ appearing in the assumptions of Theorem \ref{them2}(ii) and (iii) is allowed to be $0$. In such cases, $\varrho_{f^{*}}=0$ and the corresponding limit random variable is degenerate. If $f^{*}\not\equiv 0$, then $\epsilon(f)=\lambda_{1}/2$ (resp. $\epsilon(f)=\epsilon\in [0,\lambda_{1}/2)$) in
Theorem \ref{them2}(ii)
(resp. Theorem \ref{them2}(iii)), which corresponds to the critical (resp. large) branching case.

\noindent(ii) Although we have stated the assumptions of Theorem \ref{them2}(ii) and (iii) in the simple form of convergence of $\e^{-\lambda_{1}t}T_{t}f-\langle f,\widetilde{\varphi}\rangle\varphi$ to some limit function $f^{*}$ (with an associated rate of convergence), our approach can be extended to more general situations. For instance, Proposition \ref{prop:a} below presents other convergence conditions that fit within this framework.

{
\noindent (iii) We note that the inductive argument developed in Section \ref{sec:CLT} can be adapted to prove the finite dimensional convergence of the fluctuations in distribution. Therefore, we believe it is possible to prove a functional central limit theorem for the fluctuations, under additional assumptions on the test functions that would ensure tightness. Although we do not pursue this direction here, we refer the reader to \cite{DH} and \cite{RSZ17b} for related approaches.
}
\end{remark}

\subsection{Application to multitype superdiffusion in a bounded domain}\label{sec1.3}

As mentioned in Section \ref{sec1.2}, the examples given in \cite{GHK,LRS22,PY} satisfy Assumption \ref{AS1}. Consequently, under regularity and moment assumptions on the branching mechanisms, our main theorems are applicable.
In this section, we focus on the setting of multitype superdiffusion in a bounded domain. Our primary motivation is to investigate the fluctuations of arbitrary linear functionals for this class of superprocesses.

Suppose that $K\ge 2$, $S=\{1,\ldots,K\}$, $D$ is a bounded $C^{1,1}$ domain in $\R^{d}$ {(that is, the boundary of $D$ can be locally characterized by $C^{1,1}$ functions)}, and $m$ is the counting measure times the Lebesgue measure on $E:=S\times D$.
Suppose $\mathcal{L}^{(i)}$ ($i\in S$) is a second order differential operator of the form
	$$\mathcal{L}^{(i)}=\sum_{n,m=1}^{d}\frac{\partial}{\partial x_{n}}\left(\alpha^{(i)}_{n,m}(x)\frac{\partial}{\partial x_{m}}\right)\quad \mbox{ on }\R^{d},$$
where for every $x\in\R^{d}$, $\left(\alpha^{(i)}_{n,m}(x)\right)_{1\le n,m\leq d}$ \ is a uniformly elliptic symmetric matrix with $\alpha^{(i)}_{n,m}(x)$ being a twice differentiable function, whose second order derivatives are $\gamma$-H\"{o}lder continuous with $\gamma\in (0,1)$. Let $(\xi^{(i)},\Pi^{(i)})$ be a diffusion on $\R^{d}$ with generator $\mathcal{L}^{(i)}$. Let $\xi^{(i),D}$ be the subprocess of $\xi^{(i)}$ killed upon leaving $D$ and $P^{(i),D}_{t}$ be the corresponding semigroup.
For $f\in \mathcal{B}^{+}(E)$, we use the convention that $\boldsymbol{f}(x)=(f_{1}(x),\ldots,f_{K}(x))^{\mathrm{T}}:=(f(1,x),\ldots,f(K,x))^{\mathrm{T}}$.

Let $\xi$ be a Markov process on $E$ with semigroup \  $P_{t}f(i,x):=P^{(i),D}_{t}f_{i}(x)$ for $f\in\mathcal{B}^{+}(E)$ and $(i,x)\in E$.
	Define the branching mechanism
	$$\Psi((i,x),f):=\psi(i, \boldsymbol{f}(x)),\qquad \forall (i,x)\in E,\ f\in\mathcal{B}^{+}(E),$$
	where for $i\in E$ and $\boldsymbol{u}=(u_{1},u_{2},\cdots,u_{K})^{\mathrm{T}}\in [0,+\infty)^{K}$,
		\begin{equation}\nonumber
		\psi(i,\boldsymbol{u}):=a_{i}u_{i}+b_{i}u_{i}^{2}- \boldsymbol{u}\cdot \boldsymbol{\eta}_i +\int_{(0,+\infty)^{K}}\left(\e^{-\boldsymbol{u}\cdot \boldsymbol{y}}-1+\boldsymbol{u}\cdot\boldsymbol{y}\right)\Gamma_{i}(\ud \boldsymbol{y}).
		\end{equation}
Here, $\boldsymbol{u}\cdot \boldsymbol{y}=\sum_{i\in E} u_iy_i$ is the inner product of two vectors, $a_{i}\in (-\infty,+\infty)$, $b_{i}\ge 0$, $\boldsymbol{\eta}_i=(\eta_{i1},\cdots,\eta_{iK})^{\mathrm{T}}\in [0,\infty)^K$, and   $ \Gamma_{i}(\ud \boldsymbol{y})$ is a measure on $(0,+\infty)^{K}$ such that
		$$\int_{(0,+\infty)^{K}}(\boldsymbol{1}\cdot\boldsymbol{y})\wedge (\boldsymbol{1}\cdot\boldsymbol{y})^2 \Gamma_i(\ud  \boldsymbol{y})<+\infty\ \mbox{ and }\ \int_{(0,+\infty)^{K}}y_{j}\Gamma_{i}(\ud \boldsymbol{y})\le \eta_{ij}\quad\mbox{ for }i\not=j\in E.$$
		Without loss of generality
		we can assume that $\eta_{ii}=0$ for all $i\in E$
		(otherwise, we can
		change the value to $a_i$).

As a special case of the model given in Section \ref{sec1.1}, we have a $(P_{t},\Psi)$-superprocess $\{X_{t}:t\ge 0\}$ in $\mf$. For every $i\in E$ and $\mu\in\mf$, we define \ $\mu^{(i)}(A):=\mu(\{i\}\times A)$.\  The map \ $\mu\mapsto (\mu^{(1)},\cdots,\mu^{(K)})^{\mathrm{T}}$ \ is a homeomorphism between $\mf$ and $\mathcal{M}(D)^{K}$. Hence, \ $\{(X^{(1)}_{t},\cdots,X^{(K)}_{t})^{\mathrm{T}}:t\ge 0\}$ is a Markov process in $\mathcal{M}(D)^{K}$, which is called a {\it $K$-type superdiffusion in $D$}.

Let us denote by $T_{t}$ the mean semigroup of $X_{t}$, that is,
$$T_{t}f(u):=\pp_{\delta_{u}}\big[\langle f,X_{t}\rangle\big]=\pp_{\delta_{u}}\big[\sum_{j\in S}\langle f_{j},X^{(j)}_{t}{\rangle}\big]\quad\forall f\in\mathcal{B}^{+}(E),\ u\in E.$$
It is proved in \cite[Example 3.8]{PY} that, $T_{t}$ admits a transition density $\mathfrak{p}(t,u,v)$ with respect to $m$, such that, for every $t>0$, $(u,v)\mapsto \mathfrak{p}(t,u,v)$ is jointly continuous. Moreover, there are positive constants $a_{0}$, $t_{0}$ and $c_{i}$, $i=1,\cdots,4$, such that
	\begin{equation}\label{eg1.1}
	c_{1}\e^{a_{0}t}p_{0}(c_{2}t,x,y)\le \mathfrak{p}(t,(i,x),(j,y))\le c_{3}\e^{a_{0} t}p_{0}(c_{4}t,x,y),\qquad \forall t\in (0,t_{0}],\  (i,x),\ (j,y)\in E.
	\end{equation}
Here $p_{0}(t,x,y)$ denotes the transition density of a killed Brownian motion in $D$.
We recall that there are positive constants $c_{i},i=5,6$ such that
$$p_{0}(t,x,y)
	\le c_{5}\left(\frac{\delta_{D}(x)}{\sqrt{t}}\wedge 1\right)\left(\frac{\delta_{D}(y)}{\sqrt{t}}\wedge 1\right)t^{-\frac{d}{2}}\e^{-\frac{c_{6}|x-y|^{2}}{2}}\quad \forall t\in [0,1),\ x,y\in D.$$
Here $\delta_{D}(x)$ denotes the Euclidean distance between $x$ and the boundary of $D$.
Thus, it follows by \eqref{eg1.1} and the Chapman-Kolmogorov's equation that, for every $t>0$, there exists a constant $C_{1}=C_{1}(t)>0$, such that
\begin{equation}\label{eg1.2}
\mathfrak{p}(t,(i,x),(j,y))\le C_{1}(t)\delta_{D}(x)\delta_{D}(y)\quad\forall (i,x),(j,y)\in E.
\end{equation}
 Since $\int_{E}\int_{E}\mathfrak{p}(t,u,v)^{2}m({ \ud u})m({ \ud v})<+\infty$ for every $t>0$,  $T_{t}$ is a Hilbert-Schmidt operator in $L^{2}(E,m)$ and hence is compact. The same is true for its dual operator $\widehat{T}_{t}$. If we use $\sigma(\mathcal{L)}$ and $\sigma(\widehat{\mathcal{L}})$ to denote the spectrum of the generators of $T_{t}$ and $\widehat{T}_{t}$, respectively, then by Jentzch's theorem, $\lambda_{1}:=\sup \re (\sigma(\mathcal{L}))=\sup \re (\sigma(\widehat{\mathcal{L}}))$ is a simple eigenvalue for both $\mathcal{L}$ and $\widehat{\mathcal{L}}$, and an eigenfunction $\varphi$ of $\mathcal{L}$ associated with $\lambda_{1}$ and an eigenfunction $\widehat{\varphi}$ of $\widehat{\mathcal{L}}$ associated with $\lambda_{1}$ can be chosen strictly positive and continuous on $E$, and
	satisfying
	$\int_{E}\varphi^{2}\ud m=\int_{E}\varphi\widehat{\varphi}\ud m=1$.
We assume $\lambda_{1}>0$.
	It is further proved in \cite{PY} that there are positive constant $c_{i}, i=7,8,9$ and $t_{1}$, such that
	\begin{equation}\label{eg6.3.2}
	c^{-1}_{7}\delta_{D}(x)\le \varphi(i,x),\widehat{\varphi}(i,x)\le c_{7}\delta_{D}(x),\qquad \forall (i,x)\in E,
	\end{equation}
and
$$\Big|\frac{\e^{-\lambda_{1}t}\mathfrak{p}(t,u,v)}{\varphi(u)\widehat{\varphi}(v)}-1\Big|\le c_{8}\e^{-c_{9}t}\quad\forall t\ge t_{1},\ u,v\in E.$$
Based on these facts, Assumptions \ref{AS1} and \ref{AS2} are verified in \cite[Example 3.8]{PY}.
Additionally, we assume the fourth moment condition
\begin{equation}\label{condi:4thmoment}
\int_{(0,+\infty)^{K}}\big(\sum_{j=1}^{K}y_{j}\big)^{4}\Gamma_{i}({ \ud \boldsymbol{y}})<+\infty\quad\forall i\in S.\end{equation}

Let $W^{\varphi}_{\infty}$ be the limit of the nonnegative martingale $W^{\varphi}_{t}:=\e^{-\lambda_{1}t}\langle \varphi,X_{t}\rangle$. Then by the law of large numbers,
for any $f\in\mathcal{B}_{b}(E)$ and $\mu\in\mf$,
$$\lim_{t\to+\infty}\e^{-\lambda_{1}t}\langle f,X_{t}\rangle =W^{\varphi}_{\infty}\int_{E}f\widehat{\varphi}\ud m\mbox{ in }L^{2}(\pp_{\mu})\mbox{ and in probability}.$$
In fact, if $f$ is a compactly supported continuous function, the above convergence also holds almost surely (see, \cite[Example 3.8]{PY}).

We now turn to the following problem: Given $f\in\mathcal{B}_{b}(E)$ such that
\begin{equation}\label{eg1.7}
\int_{E}f\widehat{\varphi}\ud m=0\mbox{ and }f\not\equiv 0,
\end{equation}
we aim to describe the asymptotic behaviour of the linear functional $\langle f,X_{t}\rangle$ in a manner more precise than the estimate $\langle f,X_{t}\rangle =o(\e^{\lambda_{1}t})$.	

To achieve this, we first establish a spectral decomposition for the mean semigroup $T_{t}$.
Let $L^{2}(E,m,\mathbb{C})$ denote the space of $L^{2}(m)$-integrable $\mathbb{C}$-valued functions, with inner product $\langle f,g\rangle_{m}:=\int_{E}f\overline{g}\ud m$ and norm $\|f\|_{2}:=\langle f,f\rangle_{m}$.
By \eqref{eg1.2}, we can easily verify that the conditions (i)-(ii) of \cite[Section 1.2]{RSZ17a} are satisfied, and thus $T_{t}$ and its adjoint semigroup $\widehat{T}_{t}$ are strongly continuous semigroups on $L^{2}(E,m,\mathbb{C})$. Let $\sigma(T_{t})$ be the spectrum of $T_{t}$ in $L^{2}(E,m,\mathbb{C})$. By the compactness of $T_{t}$, the spectrum of the generator $\mathcal{L}$ of $T_{t}$ consists of countably many eigenvalues only. Moreover, if we arrange the eigenvalues $\{\lambda_{k}:k\ge 1\}$ of $\mathcal{L}$ such that $\lambda_{1}>\re \lambda_{2}\ge \re \lambda_{3}\ge \cdots$, then we have $\lim_{k\to+\infty}\re \lambda_{k}=-\infty$ and $\sigma(T_{t})\setminus\{0\}=\{\e^{\lambda_{k}t}:k\ge 1\}$.
We note that if $\e^{\lambda_{k}t}\in\sigma(T_{t})$ then $\e^{\overline{\lambda_{k}}t}\in\sigma(T_{t})$. Hence for each $k$, there exists a unique integer $k'$ such that $\lambda_{k'}=\overline{\lambda_{k}}$.
We will use the spectral theory of (not necessarily symmetric) strongly continuous semigroups, which is developed in \cite{RSZ17a}, as follows.
There exists $t^{*}>0$ such that $\e^{\lambda_{k}t^{*}}\not=\e^{\lambda_{j}t^{*}}$ for any $k\not=j$. For each $k\ge 1$, define $\mathcal{N}_{k,0}:=\{0\}$ and
$$\mathcal{N}_{k,n}:=\{f\in L^{2}(E,m,\mathbb{C}):\ \big(\e^{\lambda_{k}t^{*}}I-T_{t^{*}}\big)^{n}f=0\}\quad\forall n\ge 1.$$
Then $\mathcal{N}_{k,n}$ is a finite dimensional linear subspace of $L^{2}(E,m,\mathbb{C})$, and there is an integer $v_{k}\ge 1$ such that
$$\mathcal{N}_{k,n}\subsetneq \mathcal{N}_{k,n+1}\mbox{ for }n=0,1,\cdots,v_{k}-1,\mbox{ and }\mathcal{N}_{k,n}=\mathcal{N}_{k,n+1}\mbox{ for }n\ge v_{k}.$$
Let $n_{k}:=\mathrm{dim}(\mathcal{N}_{k,v_{k}})$, $r_{k}:=\mathrm{dim}(N_{k,1})$, and
$\mathcal{N}_{k}:=\mathcal{N}_{1,v_{1}}\oplus\mathcal{N}_{2,v_{2}}\oplus\cdots\oplus \mathcal{N}_{k,v_{k}}$.
The same discussion applies to the adjoint semigroup $\widehat{T}_{t}$, and the corresponding objects for $\widehat{T}_{t}$ will be denoted with a $\widehat{\ }$. We note that $\widehat{v}_{k}=v_{k}$, $\widehat{n}_{k}=n_{k}$ and $\widehat{r}_{k}=r_{k}$.

We use the convention that when a operator acts on a vector-valued function, it acts componentwise.
For each $k\ge 1$, there exist a basis $\{\phi^{(k)}_{j}:1\le j\le n_{k}\}$ of $\mathcal{N}_{k,v_{k}}$ and integers $d_{k,j},1\le j\le r_{k}$, such that, each $\phi^{(k)}_{j}$ is a continuous function on $E$, $\sum_{j=1}^{r_{k}}d_{k,j}=n_{k}$, and for $m$-a.e. $u\in E$,
\begin{equation}\label{eg1.3}
(T_{t}\Phi_{k})^{\mathrm{T}}(u)=\e^{\lambda_{k}t}\Phi^{\mathrm{T}}_{k}(u)D_{k}(t),
\end{equation}
where $\Phi_{k}:=(\phi^{(k)}_{1},\cdots,\phi^{(k)}_{n_{k}})^{\mathrm{T}}$, $D_{k}(t)$ is an $n_{k}\times n_{k}$ matrix given by
\begin{equation*}
D_{k}(t):=
\begin{pmatrix}
J_{k,1}(t)&&&\boldsymbol{0}\\
 &J_{k,2}(t)&&\\
 &&\ddots&\\
 \boldsymbol{0}&&&J_{k,r_{k}}(t)
\end{pmatrix}
\end{equation*}
and $J_{k,j}(t)$ is a $d_{k,j}\times d_{k,j}$ matrix given by
\begin{equation*}
J_{k,j}(t):=
\begin{pmatrix}
1&t&\frac{t^{2}}{2!}&\cdots&\frac{t^{d_{k,j}-1}}{(d_{k,j}-1)!}\\
0&1&t&\cdots&\frac{t^{d_{k,j}-2}}{(d_{k,j}-2)!}\\
 &&\ddots&\ddots\\
0&0&\cdots&1&t\\
0&0&\cdots&0&1
\end{pmatrix}.
\end{equation*}
In particular, $D_{1}(t)=1,\ n_{1}=r_{1}=1$ and $\phi^{(1)}_{1}=\varphi$.

We note that by the H\"{o}lder's inequality and \eqref{eg1.2}, for every $t>0$, there is a constant $C_{2}=C_{2}(t)>0$ such that
$$|T_{t}\phi^{(k)}_{j}(u)|\le |\int_{E}\mathfrak{p}(t,u,v)^{2}m({ \ud v})|^{1/2}\|\phi^{(k)}_{j}\|_{2}\le C_{2}\delta_{D}(x)\|\phi^{(k)}_{j}\|_{2}\quad\forall u=(i,x)\in E.$$
Hence by \eqref{eg1.3}, $\phi^{(k)}_{j}$ is a bounded function on $E$. Then by the joint continuity of $\mathfrak{p}(t,u,v)$ in $(u,v)$ and the dominated convergence theorem, we can prove that $u\mapsto T_{t}\phi^{(k)}_{j}(u)$ is continuous on $E$. Thus the equality \eqref{eg1.3} holds indeed for every $u\in E$.

It follows by \cite[Lemma 1.9 and Corollary 1.12]{RSZ17a} that, for each $k\ge 1$, there exists a unique basis
$\{\widehat{\phi}^{(k)}_{1},\cdots,\widehat{\phi}^{(k)}_{n_{k}}\}$ of $\widehat{\mathcal{N}}_{k,v_{k}}$ such that
$$\widehat{T}_{t}\widehat{\Phi}_{k}(u)=\e^{\overline{\lambda_{k}}t}D_{k}(t)\widehat{\Phi}_{k}(u)\mbox{ and }\big(\langle \phi^{(k)}_{j},\widehat{\phi}^{(k)}_{l}\rangle_{m}\big)_{1\le j,l\le n_{k}}=I,$$
where $\widehat{\Phi}_{k}:=(\widehat{\phi}^{(k)}_{1},\cdots,\widehat{\phi}^{(k)}_{n_{k}})^{\mathrm{T}}$. Moreover, for every $g\in L^{2}(E,m,\mathbb{C})$, there exists a unique $\widetilde{g}_{k}\in \widehat{\mathcal{N}_{k}}^{\perp}$, such that,
$$g=\sum_{j=1}^{k}\Phi^{\mathrm{T}}_{j}\langle g,\widehat{\Phi}_{j}\rangle_{m}+\widetilde{g}_{k},$$
where $\langle g,\widehat{\Phi}_{j}\rangle_{m}=(\langle g,\widehat{\phi}^{(j)}_{1}\rangle_{m},\cdots,\langle g,\widehat{\phi}^{(j)}_{n_{j}}\rangle_{m})^{\mathrm{T}}$.
We also note that the bases of $\mathcal{N}_{k,v_{k}}$ and $\widehat{\mathcal{N}}_{k,v_{k}}$ can be chosen to satisfy that $\Phi_{k'}=\overline{\Phi_{k}}$ and $\widehat{\Phi}_{k'}=\overline{\widehat{\Phi}_{k}}$.

For the remainder of this section, we fix a function $f$ that satisfies \eqref{eg1.7}.
Since $\langle f,\widehat{\Phi}_{1}\rangle_{m}=\langle f,\widehat{\varphi}\rangle_{m}=0$, we have
$$f=\sum_{j=2}^{k}\Phi^{\mathrm{T}}_{j}\langle f,\widehat{\Phi}_{j}\rangle_{m}+\widetilde{f}_{k},$$
where $\widetilde{f}_{k}\in\widehat{\mathcal{N}}^{\perp}_{k}$. Then by \eqref{eg1.3},
\begin{equation}\label{eg1.5}
T_{t}f(u)=\sum_{j=2}^{k}\e^{\lambda_{j}t}\Phi^{\mathrm{T}}_{j}(u)D_{j}(t)\langle f,\widehat{\Phi}_{j}\rangle_{m}+T_{t}\widetilde{f}_{k}(u).
\end{equation}
We observe that each component of the vector $D_{j}(t)\langle f,\widehat{\Phi}_{j}\rangle_{m}$ is a polynomial of $t$. Specifically, if we denote by $r^{(j)}_{n}(f)$ the degree of the $n$-th component of $D_{j}(t)\langle f,\widehat{\Phi}_{j}\rangle_{m}$, then $r^{(j)}_{n}(f)=\sup\{l\ge 0:\ \langle f,\widehat{\phi}^{(j)}_{n+l}\rangle_{m}\not=0\}$, and the $n$-th component of $D_{j}(t)\langle f,\widehat{\Phi}_{j}\rangle_{m}$ is given by
$$\langle f,\widehat{\phi}^{(j)}_{n}\rangle_{m}+t\langle f,\widehat{\phi}^{(j)}_{n+1}\rangle_{m}+\frac{t^{2}}{2!}\langle f,\widehat{\phi}^{(j)}_{n+2}\rangle_{m}+\cdots
+\frac{t^{r^{(j)}_{n}(f)}}{r^{(j)}_{n}(f)!}\langle f,\widehat{\phi}^{(j)}_{n+r^{(j)}_{n}(f)}\rangle_{m}.$$
We define three parameters $\alpha(f)$, $s(f)$ and $\gamma(f)$, and an index set $\mathfrak{I}(f)$  for $f$ by the following relations.
\begin{eqnarray*}
\alpha(f)&:=&\sup\{\re\lambda_{j}:\ \langle f,\widehat{\Phi}_{j}\rangle_{m}\not=0,\ j\ge 2\},\\
s(f)&:=&\sup\{j:\ \re\lambda_{j}=\alpha(f)\},\\
\gamma(f)&:=&\sup\{r^{(j)}_{n}(f):\ j\ge 2,\ \re \lambda_{j}=\alpha(f)\mbox{ and }1\le n\le n_{j}\},\\
\mathfrak{I}(f)&:=&\{j:\ \re\lambda_{j}=\alpha(f)\mbox{ and }\max_{1\le n\le n_{j}}r^{(j)}_{n}(f)=\gamma(f)\}.
\end{eqnarray*}
For instance, if we take $f=\phi^{(k)}_{l}$, one of the basis functions of the space $\mathcal{N}_{k,v_{k}}$, for some integer $l\in (\sum_{j=1}^{m-1}d_{k,j},\sum_{j=1}^{m}d_{k,j}]$ (under the convention that $d_{k,0}=0$), then from \eqref{eg1.3}, $(\alpha(f),s(f),\gamma(f),\mathfrak{I}(f))$ takes the value $(\re\lambda_{k},\max\{k,k'\},l-1-\sum_{j=1}^{m-1}d_{k,j},\{k,k'\})$.
We now present a proposition that characterizes the long time behaviour of $T_{t}f$.

\begin{proposition}\label{prop:a}
For $f$ satisfying \eqref{eg1.7},
we have $$\frac{\e^{-\alpha(f)t}T_{t}f}{t^{\gamma(f)}}-\frac{1}{\gamma(f)!}\sum_{j\in \mathfrak{I}(f)}\e^{i t\im \lambda_{j}}F_{j}(f)\rightrightarrows 0\mbox{ on $E$ as }t\to+\infty,$$
where $F_{j}(f)(u):=\sum_{n\in \mathfrak{M}_{j}(f)}\phi^{(j)}_{n}(u)\langle f,\widehat{\phi}^{(j)}_{n+\gamma(f)}\rangle_{m}$ for $j\in \mathfrak{I}(f)$ and $u\in E$.
\end{proposition}
The proof of this proposition is deferred to Section \ref{sec5}.

We recall that
$$\epsilon(f)=-\limsup_{t\to+\infty}\frac{\log\sup_{u\in E}|\varphi(u)^{-1}\e^{-\lambda_{1}t}T_{t}f(u)|}{t}.$$
Using the facts that for each $k\ge 1$, $1\le j\le n_{k}$ and $u=(i,x)\in E$, $\phi^{(k)}_{j}(u)\lesssim \delta_{D}(x)$, $\varphi(u)\asymp\delta_{D}(x)$, and for $t\ge 1$ and $g\in\widehat{\mathcal{N}}_{k}^{\perp}$, $T_{t}g(u)\lesssim \e^{t\re \lambda_{k+1}}\delta_{D}(x)$ (see \eqref{a.2} below), we get from \eqref{eg1.6} below that for every $t\ge 1$,
$$\big|\varphi(u)^{-1}T_{t}f(u)\big|\lesssim \e^{\alpha(f)t}(1+t+t^{2}+\cdots+t^{\gamma(f)}).$$
This implies that $\epsilon(f)\ge \lambda_{1}-\alpha(f)$. As a consequence of Theorem \ref{them2}(i), we obtain the following result.
\begin{proposition}\label{prop:b}
If $\alpha(f)<\lambda_{1}/2$, then for every $u\in E$, under $\pp_{\delta_{u}}$,
$$\e^{-\frac{\lambda_{1}}{2}t}\langle f,X_{t}\rangle\stackrel{d}{\to}\rho_{f}\sqrt{W^{\varphi}_{\infty}}N\mbox{ as }t\to+\infty,$$
where \[\rho^{2}_{f}:=\int_{0}^{+\infty}\e^{-\lambda_{1}s}\langle \vartheta[T_{s}f],\widehat{\varphi}\rangle_{m}{ \ud s},\]
with $\vartheta[f](i,x):=2b_{i}f^{2}_{i}(x)+\int_{(0,+\infty)^{K}}(\sum_{j\in S}f_{j}(x)y_{j})^{2}\,\Gamma_{j}({ \ud \boldsymbol{y}})$, and $N$ is a standard normal random variable independent of $W^{\varphi}_{\infty}$.
\end{proposition}
We thus focus on the case where $\alpha(f)\ge \lambda_{1}/2$.
First, consider the simple situation in which $\mathfrak{I}(f)$ is a singleton set. Suppose $\mathfrak{I}(f)=\{\kappa\}$. Then $\lambda_{\kappa}$ must be a real number.
Otherwise, if there exists a integer $\kappa'\not=\kappa$ such that $\lambda_{\kappa'}=\overline{\lambda_{\kappa}}$, we would have
$\langle f,\widehat{\Phi}_{\kappa'}\rangle_{m}=\overline{\langle f,\widehat{\Phi}_{\kappa}\rangle_{m}}\not=0$, and $\kappa'\in \mathfrak{I}(f)$, which leads to a contradiction.
When $\mathfrak{I}(f)=\{\kappa\}$,  Proposition \ref{prop:a} yields that
$$\frac{\e^{-\alpha(f)t}T_{t}f}{t^{\gamma(f)}}\rightrightarrows\frac{1}{\gamma(f)!}F_{\kappa}(f)\mbox{ on $E$, as }t\to+\infty.$$
Hence, in this case, Theorem \ref{them2}(ii) and (iii) are applicable.

On the other hand, the assumptions of Theorem \ref{them2}(ii) and (iii) do not cover the case where $\mathfrak{I}(f)$ is not a singleton set. In such situations, the function $\sum_{j\in\mathfrak{I}(f)}\e^{it\im \lambda_{j}}F_{j}(f)$ from Proposition \ref{prop:a} depends explicitly on $t$, preventing direct application of those results. Nevertheless, our approach can be extended to handle both cases--whether $\mathfrak{I}(f)$ is a singleton set or not. We thus obtain the following result, whose proof is deferred to Section \ref{sec5}.
\begin{proposition}\label{prop:d}
\item{(i)} If $\alpha(f)=\lambda_{1}/2$, then
$$\frac{\e^{-\lambda_{1}t}}{t^{1+2\gamma(f)}}\varphi(u)^{-1}\mathrm{Var}_{\delta_{u}}[\langle f,X_{t}\rangle]\stackrel{b}{\to}\frac{\varrho^{2}_{f}}{1+2\gamma(f)}\mbox{ as }t\to+\infty,$$
where $\varrho^{2}_{f}=(\gamma(f)!)^{-2}\sum_{j\in \mathfrak{I}(f)}\langle { \vartheta[F_{j}(f),\overline{F_{j}(f)}]},\widehat{\varphi}\rangle_{m}$.

If, in addition, $\pp_{\delta_{u}}\left(X_{t}=0\mbox{ for some }t\in (0,+\infty)\right)>0$, then under $\pp_{\delta_{u}}$,
$$\frac{\e^{-\frac{\lambda_{1}}{2}t}}{t^{\frac{1}{2}+\gamma(f)}}\langle f,X_{t}\rangle\stackrel{d}{\to}\frac{1}{\sqrt{1+2\gamma(f)}}\varrho_{f}\sqrt{W^{\varphi}_{\infty}}N\mbox{ as }t\to+\infty,$$
where $N$ is a standard normal random variable independent of $W^{\varphi}_{\infty}$.
\item{(ii)} If $\alpha(f)>\lambda_{1}/2$, then for every $j\in \mathfrak{I}(f)$ and $u\in E$, $W^{(j)}_{t}(f):=\e^{-\lambda_{j}t}\langle F_{j}(f),X_{t}\rangle$ is an $L^{2}(\pp_{\delta_{u}})$-bounded martingale, thus converges in $L^{2}(\pp_{\delta_{u}})$ and $\pp_{\delta_{u}}$-a.s. to a limit $W^{(j)}_{\infty}(f)$. Moreover,
    $$\frac{\e^{-\alpha(f)t}}{t^{\gamma(f)}}\langle f,X_{t}\rangle-\frac{1}{\gamma(f)!}\sum_{j\in \mathfrak{I}(f)}\e^{it\im\lambda_{j}}W^{(j)}_{\infty}(f)\to 0\mbox{ in }L^{2}(\pp_{\delta_{u}})\mbox{ as }t\to+\infty.$$
\end{proposition}

\section{Preliminaries}\label{sec2}
{From now on, we assume Assumptions \ref{AS1}-\ref{AS3} hold. Other conditions needed will be stated explicitly.}
\subsection{A stochastic integral representation of superprocess}\label{sec:stochastic integral}

Unless stated otherwise, all martingales or local martingales mentioned in this subsection will be relative to the filtration $(\mathcal{F}_{t})_{t\ge 0}$ and probability measure $\pp_{\mu}$ where $\mu\in\mf$.

Let $N(\ud s,\ud \nu)$ be the random measure on $\R^{+}\times \mathcal{M}(E)^{0}$ defined by
\begin{equation*}
N(\ud s,\ud \nu):=\sum_{s\ge 0}1_{\{\triangle X_{s}\not=0\}}\delta_{(s,\triangle X_{s})}(\ud s,\ud \nu).
\end{equation*}
Here, we use the standard notation \ $\triangle X_{s}:=X_{s}-X_{s-}$ \ for the jump of $X$ at time $s$. Let $\widehat{N}(\ud s,\ud \nu)$ be the predictable compensator of $N(\ud s,\ud \nu)$, and $\widetilde{N}(\ud s,\ud \nu):=N(\ud s,\ud \nu)-\widehat{N}(\ud s,\ud \nu)$ be the compensated random measure.
By \cite[Theorem 4.2]{PY} $\widehat{N}$ is given by
$$\widehat{N}({ \ud s},{ \ud \nu})={ \ud s}\int_{E}X_{s-}({ \ud x})H(x,{ \ud \nu}).$$
Let $(\omega,s,\nu)\mapsto F(\omega,s,\nu)$ be a predictable function on $\mathcal{W}^{+}_{0}\times \R^{+}\times \mfo $ such that
\begin{equation}\label{condi:integral N}
\pp_{\mu}\left[\left(\sum_{s\le t}F(s,\triangle X_{s})^{2} 1_{\{\triangle X_{s}\not=0\}}\right)^{1/2}\right]<+\infty,\qquad\forall t\ge 0.
\end{equation}
Then, following \cite[Chapter II.1d]{JS}, one can define the stochastic integral of $F$ with respect to the compensated measure $\widetilde{N}(\ud s,\ud \nu)$, denoted by
$$\int_{0}^{t}\int_{\mfo }F(s,\nu)\widetilde{N}(\ud s,\ud \nu),$$
as the unique purely discontinuous local martingale whose jumps are indistinguishable from the process \ $F(s,\triangle X_{s})1_{\{\triangle X_{s}\not=0\}}$. It is proved in \cite[section 4.2]{PY} that condition \eqref{condi:integral N} is satisfied for $F(\omega,s,\nu)=\int_{E}g(s,x)\nu(\ud x)$,
where $g(s,x)$ is a bounded measurable function on $\R^{+}\times E$. In this case, we also write $\int_{0}^{t}\int_{\mfo }\langle g(s,\cdot),\nu\rangle\widetilde{N}({ \ud s},{ \ud \nu})$ as $\int_{0}^{t}\int_{E}g(s,x)M^{d}({ \ud s},{ \ud x})$.

On the other hand, as has been shown in \cite[section 4.2]{PY}, there exists a martingale measure $M^{c}({ \ud s},{ \ud x})$ on $\R^{+}\times E$ with covariance measure
$${ \ud s}\int_{E}X_{s}({ \ud z})2b(z)\delta_{z}({ \ud x})\delta_{z}({ \ud y}),$$
such that,  if $(\omega,s,x)\mapsto G(\omega,s,x)$ is a predictable function on $\mathcal{W}^{+}_{0}\times \R^{+}\times E$ such that
\begin{equation}\label{condi:integral Mc}
\pp_{\mu}\left[\int_{0}^{t}\ud s\int_{E}2b(x)G^{2}(s,x)X_{s}(\ud x)\right]<+\infty,\qquad \forall t\ge 0,
\end{equation}
then the stochastic integral of $G$ with respect to the martingale measure $M^{c}(\ud s,\ud x)$,
denoted by
$$\int_{0}^{t}\int_{E}G(s,x)M^{c}(\ud s,\ud x),$$
is well-defined, and is the unique square integrable c\`{a}dl\`{a}g martingale with predictable quadratic variation
$$2\int_{0}^{t}\langle b G^{2}(s,\cdot),X_{s}\rangle \ud s.$$

From \eqref{eq:many-to-one}, we have for all $x\in E$, $t\ge 0$ and $f\in\mathcal{B}^{+}_{b}(E)$,
$$T_{t}f(x)\le P_{t}f(x)+\Pi_{x}\left[\int_{0}^{t}a^{-}(\xi_{s})T_{t-s}f(\xi_{s}){ \ud s}\right]+\Pi_{x}\left[\int_{0}^{t}m[T_{t-s}f](\xi_{s}){ \ud s}\right].$$
Hence we get
$$\|T_{t}f\|_{\infty}\le \|f\|_{\infty}+c_{1}\int_{0}^{+\infty}\|T_{s}f\|_{\infty}{ \ud s}\quad\forall t\ge 0,$$
where $c_{1}=\|a^{-}+m[1]\|_{\infty}$. It then follows by the Growall's lemma that
\begin{equation}\label{ieq:many-to-one}
\|T_{t}f\|_{\infty}\le \e^{c_{1}t}\|f\|_{\infty}\quad\forall t\ge 0.
\end{equation}
By the above inequality, we can easily show that \eqref{condi:integral Mc} is satisfied for $G(\omega,s,x)=g(s,x)$ where $g\in  \mathcal{B}_{b}(\R^{+}\times E)$.

Using the above stochastic integrals, \cite[Proposition 4.3]{PY} shows that, for every $f\in \mathcal{B}_{b}(E)$, $t\ge 0$ and $\mu\in\mf$,
\begin{equation}\label{eq:stochastic integral representation}
\langle f,X_{t}\rangle=\langle T_{t}f,X_{0}\rangle+\int_{0}^{t}\int_{E}T_{t-s}f(x)M({ \ud s},{ \ud x})\quad\pp_{\mu}\mbox{-a.s.}
\end{equation}
Here $\int_{0}^{t}\int_{E}T_{t-s}f(x)M({ \ud s},{ \ud x}):=\int_{0}^{t}\int_{E}T_{t-s}f(x)M^{c}({ \ud s},{ \ud x})+\int_{0}^{t}\int_{E}T_{t-s}f(x)M^{d}({ \ud s},{ \ud x})$ is well-defined since $(s,x)\mapsto 1_{\{s\le t\}}T_{t-s}f(x)$ is a bounded measurable function on $\R^{+}\times E$. In particular, taking $f=\varphi$, we have
$${ W^{\varphi}_{t}}={ W^{\varphi}_{0}}+\int_{0}^{t}\int_{E}\e^{-\lambda_{1} s}\varphi(x)M({ \ud s},{ \ud x})\quad \pp_{\mu}\mbox{-a.s.}$$

\subsection{First order asymptotic behavior of linear functionals}
This section is devoted to the proof of Theorem \ref{cor:wlln}. The approach is borrowed from \cite{PY}, which involves a combination of stochastic analysis and the truncation procedure of Asmussen and Hering \cite{AH76}.

Recall that $\stackrel{b}{\to}$ denotes the bounded pointwise convergence.
For example,
{if we let $f^{+}:=f\vee 0$ and $f^{-}:=(-f)\vee 0$, then}
Assumption \ref{AS1} implies that for every $f\in \mathcal{B}_{b}(E)$ and $x\in E$,
\begin{align}\label{ieq:H3}
&\left|\varphi(x)^{-1}\e^{-\lambda_{1} t}T_{t}f(x)-\langle f,\widetilde{\varphi}\rangle\right|\nonumber\\
\le& \left|\varphi(x)^{-1}\e^{-\lambda_{1} t}T_{t}f^{+}(x)-\langle f^{+},\widetilde{\varphi}\rangle\right|+\left|\varphi(x)^{-1}\e^{-\lambda_{1} t}T_{t}f^{-}(x)-\langle f^{-},\widetilde{\varphi}\rangle\right|\nonumber\\
\le&(\|f^{+}\|_{\infty}+\|f^{-}\|_{\infty})\triangle_{t}.
\end{align}
Hence $\e^{-\lambda_{1} t}T_{t}f\stackrel{b}{\to}\varphi\langle f,\widetilde{\varphi}\rangle$ as $t\to+\infty$ for every $f\in\mathcal{B}_{b}(E)$.

\begin{lemma}\label{lem2.0}
Suppose $f_{t},f\in\mathcal{B}_{b}(E)$ and $f_{t}\stackrel{b}{\to}f$ as $t\to+\infty$.
Then for every $\mu\in\mf$,
\begin{equation}\nonumber
\lim_{s\to+\infty}\lim_{t\to +\infty}\e^{-\lambda_{1}(t+s)}\pp_{\mu}\left[\langle f_{t+s},X_{t+s}\rangle|\mathcal{F}_{t}\right]=\langle f,\widetilde{\varphi}\rangle { W^{\varphi}_{\infty}}\qquad \pp_{\mu}\mbox{-a.s. and in }L^{2}(\pp_{\mu}).
\end{equation}
\end{lemma}

\begin{proof}
Without loss of generality we assume $f_{t},f\in\mathcal{B}^{+}_{b}(E)$.
For any $s,t\ge 0$, by the Markov property
\begin{equation*}
\e^{-\lambda_{1}(t+s)}\pp_{\mu}\left[\langle f_{t+s},X_{t+s}\rangle|\mathcal{F}_{t}\right]
=\e^{-\lambda_{1}(t+s)}\pp_{X_{t}}\left[\langle f_{t+s},X_{s}\rangle\right]
=\e^{-\lambda_{1}(t+s)}\langle  T_{s}f_{t+s},X_{t}\rangle.
\end{equation*}
Hence, we have
\begin{eqnarray}
&&\e^{-\lambda_{1}(t+s)}\pp_{\mu}\left[\langle f_{t+s},X_{t+s}\rangle |\mathcal{F}_{t}\right]-\langle f,\widetilde{\varphi}\rangle { W^{\varphi}_{\infty}}\nonumber\\
&=&\e^{-\lambda_{1} t}\langle  \e^{-\lambda_{1} s}T_{s}f_{t+s}-\langle f_{t+s},\widetilde{\varphi}\rangle\varphi,X_{t}\rangle\nonumber\\
&&+\langle f_{t+s},\widetilde{\varphi}\rangle \left({ W^{\varphi}_{t}}-{ W^{\varphi}_{\infty}}\right)+\langle f_{t+s}-f,\widetilde{\varphi}\rangle { W^{\varphi}_{\infty}}.\label{lem2.0.1}
\end{eqnarray}
By \eqref{ieq:H3} we have
\begin{equation*}
\e^{-\lambda_{1} t}\langle  \left|\e^{-\lambda_{1} s}T_{s}f_{t+s}-\langle f_{t+s},\widetilde{\varphi}\rangle\varphi\right|,X_{t}\rangle
\le \e^{-\lambda_{1} t}\langle \varphi\|f_{t+s}\|_{\infty}\triangle_{s},X_{t}\rangle\le (\sup_{r\ge s}\|f_{r}\|_{\infty})\triangle_{s}{ W^{\varphi}_{t}}.
\end{equation*}
Since $\triangle_{s}\to 0$ as $s\to+\infty$ and ${ W^{\varphi}_{t}}\to { W^{\varphi}_{\infty}}$ $\pp_{\mu}$-a.s. and in $L^{2}(\pp_{\mu})$ as $t\to+\infty$, we get from the above inequality that
$$\lim_{s\to+\infty}\lim_{t\to+\infty}\e^{-\lambda_{1} t}\langle  \e^{-\lambda_{1} s}T_{s}f_{t+s}-\langle f_{t+s},\widetilde{\varphi}\rangle\varphi,X_{t}\rangle=0\quad\pp_{\mu}\mbox{-a.s. and in }L^{2}(\pp_{\mu}).$$
Since $\langle f_{t+s}-f,\widetilde{\varphi}\rangle\to 0$ by the dominated convergence theorem, the last two terms in the right hand side of \eqref{lem2.0.1} converge to 0 $\pp_{\mu}$-a.s. and in $L^{2}(\pp_{\mu})$.
Hence we have completed the proof.
\end{proof}

\begin{lemma}\label{lem2.1}
Suppose $\{f_{t}:t\ge 0\}\subset \mathcal{B}_{b}(E)$ and $\sup_{t\ge 0}\|f_{t}\|_{\infty}<+\infty$. Then for every $\mu\in\mf$ and $s\ge 0$,
\begin{equation*}
\lim_{t\to +\infty}\e^{-\lambda_{1} (t+s)}\left(\langle f_{t+s},X_{t+s}\rangle-\pp_{\mu}\left[\langle f_{t+s},X_{t+s}\rangle|\mathcal{F}_{t}\right]\right)=0\quad \mbox{ in }L^{2}(\pp_{\mu}).
\end{equation*}
Moreover, for every $m\in\mathbb{N}$ and $\sigma>0$,
\begin{equation*}
\lim_{\mathbb{N}\ni n\to+\infty}\e^{-\lambda_{1} (m+n)\sigma}\left(\langle f_{(m+n)\sigma},X_{(m+n)\sigma}\rangle-\pp_{\mu}\left[\langle f_{(m+n)\sigma},X_{(m+n)\sigma}\rangle|\mathcal{F}_{n\sigma}\right]\right)=0\quad\pp_{\mu}\mbox{-a.s.}
\end{equation*}
\end{lemma}

\begin{proof}
Without loss of generality we assume $f_{t},f\in\mathcal{B}^{+}_{b}(E)$. For $\tau>0$, we define
\begin{equation*}
I^{\tau}_{s,t} :=\int_{s}^{t}\int_{E} T_{\tau-u}f_{\tau}(x)M(\ud u,\ud x),
\qquad 0\le s\le t\le \tau.
\end{equation*}
and define \ $I^{\tau,c}_{s,t} $ \ and \ $I^{\tau,d}_{s,t} $ analogously with $M(\ud u,\ud x)$ replaced by $M^{c}(\ud u,\ud x)$ and $M^{d}(\ud u,\ud x)$ respectively.
We recall from section \ref{sec:stochastic integral} that $[0,\tau]\ni t\mapsto I^{\tau,d}_{0,t}$ (resp. $[0,\tau]\ni t\mapsto I^{\tau,c}_{0,t}$) is a purely discontinuous local martingale (resp. martingale) with quadratic variation $\sum_{s\le t}\langle T_{\tau-s}f_{\tau},\triangle X_{s}\rangle^{2}1_{\{\triangle X_{s}\not=0\}}$ (resp. with predictable quadratic variation $\int_{0}^{t}\langle 2b(T_{\tau-s}f_{\tau})^{2},X_{s}\rangle { \ud s}$). By the Burkholder-Davis-Gundy inequality, there is some constant $c_{1}>0$ such that
\begin{align}\label{lem2.1.1}
\pp_{\mu}\left[\sup_{t\in [0,\tau]}|I^{\tau,d}_{0,t}|^{2}\right]&\le c_{1}\pp_{\mu}\left[\sum_{s\le \tau}\langle T_{\tau-s}f_{\tau},\triangle X_{s}\rangle^{2}1_{\{\triangle X_{s}\not=0\}}\right]\nonumber\\
&=c_{1}\pp_{\mu}\left[\int_{0}^{\tau}\langle\int_{\mfo }{\langle T_{\tau-s}f_{\tau},\nu\rangle}^{2}H(\cdot,{ \ud \nu}),X_{s}\rangle { \ud s}\right].
\end{align}
Using \eqref{ieq:many-to-one}, Assumption \ref{AS3} and the fact that $\sup_{s\in [0,\tau],\ x\in E}|T_{\tau-s}f_{\tau}(x)|<+\infty$, one can show the right hand side of \eqref{lem2.1.1} is finite.
We note that $\pp_{\mu}\left[|I^{\tau,d}_{0,\zeta}|^{2}\right]\le \pp_{\mu}\left[\sup_{t\in [0,\tau]}|I^{\tau,d}_{0,t}|^{2}\right]$ for any $[0,\tau]$-valued stopping time $\zeta$. It follows that the set of random variables $\{I^{\tau,d}_{0,\zeta}:\ \zeta\mbox{ is $[0,\tau]$-valued stopping time}\}$ is uniformly integrable. Hence by \cite[Proposition I.1.47]{JS} $[0,\tau]\ni t\mapsto I^{\tau,d}_{0,t}$ is a (uniformly integrable) martingale, and consequently $[0,\tau]\ni t\mapsto I^{\tau}_{0,t}$ is a martingale.
Using this and the stochastic integral representation given in \eqref{eq:stochastic integral representation}, we have
\begin{align}
\e^{-\lambda_{1} (t+s)}\left(\langle f_{t+s},X_{t+s}\rangle-\pp_{\mu}\left[\langle f_{t+s},X_{t+s}\rangle|\mathcal{F}_{t}\right]\right)
&=\e^{-\lambda_{1} (t+s)}\left(I^{t+s}_{0,t+s}-\pp_{\mu}\left[I^{t+s}_{0,t+s} |\mathcal{F}_{t}\right]\right)\nonumber\\
&=\e^{-\lambda_{1} (t+s)}\left(I^{t+s}_{0,t+s}-I^{t+s}_{0,t}\right) \nonumber\\
&=\e^{-\lambda_{1} (t+s)}I^{t+s}_{t,t+s} .\label{lem7.3.1}
\end{align}
Using the quadratic variations, we have
\begin{align*}
\pp_{\mu}\left[\left(\e^{-\lambda_{1} (t+s)}I^{t+s,c}_{t,t+s} \right)^{2}\right]
=&\e^{-2\lambda_{1} (t+s)}\pp_{\mu}\left[\int_{t}^{t+s}\langle 2b\left(T_{t+s-r}f_{t+s}\right)^{2},X_{r}\rangle { \ud r}\right]\\
=&\int_{t}^{t+s}\e^{-2\lambda_{1} r}\pp_{\mu}\left[\langle 2b\left(\e^{-\lambda_{1} (t+s-r)}T_{t+s-r}f_{t+s}\right)^{2},X_{r}\rangle\right]{ \ud r},
\end{align*}
and
\begin{align*}
\pp_{\mu}\left[\left(\e^{-\lambda_{1} (t+s)}I^{t+s,d}_{t,t+s} \right)^{2}\right]
=&\e^{-2\lambda_{1}(t+s)}\pp_{\mu}\left[\sum_{t\le r\le t+s}\langle T_{t+s-r}f_{t+s},\triangle X_{r}\rangle^{2}\right]\\
=&\e^{-2\lambda_{1}(t+s)}\pp_{\mu}\left[
\int_{t}^{t+s}
{ \ud r}\int_{E}X_{r-}({ \ud x})\int_{\mfo }\langle T_{t+s-r}
f_{t+s}
,\nu\rangle^{2}H(x,{ \ud \nu})\right]\\
=&\int_{t}^{t+s}\e^{-2\lambda_{1} r}\pp_{\mu}\left[\langle\int_{\mfo }\langle\e^{-\lambda_{1}(t+s-r)}T_{t+s-r}
f_{t+s},\nu\rangle^{2}H(\cdot,{ \ud \nu}),X_{r}\rangle\right]{ \ud r}.
\end{align*}
By \eqref{ieq:H3}, we have
\begin{equation}\label{eq1.9}
\e^{-\lambda_{1} t}T_{t}f(x)\le \varphi(x)\|f\|_{\infty}(\sup_{s\ge 0}\triangle_{s}+\langle 1,\widetilde{\varphi}\rangle)\quad\forall x\in E,\ f\in\mathcal{B}^{+}_{b}(E),\ t\ge 0.
\end{equation}
Thus
$\left|\e^{-\lambda_{1}(t+s-r)}T_{t+s-r}f_{t+s}(x)\right|\le c_{1}\varphi(x)$ for all $x\in E$,
where $c_{1}:=
\sup_{r\ge 0}\|f_{r}\|_{\infty}
(\sup_{r\ge 0}\triangle_{r}+\langle 1,\widetilde{\varphi}\rangle)\in (0,+\infty)$.
It follows that
\begin{align*}
\pp_{\mu}\left[\left(\e^{-\lambda_{1} (t+s)}I^{t+s}_{t,t+s} \right)^{2}\right]&\le2\left(\pp_{\mu}\left[\left(\e^{-\lambda_{1} (t+s)}I^{t+s,c}_{t,t+s} \right)^{2}\right]+\pp_{\mu}\left[\left(\e^{-\lambda_{1} (t+s)}I^{t+s,d}_{t,t+s} \right)^{2}\right]\right)\\
&\le 2c^{2}_{1}\int_{t}^{t+s}\e^{-2\lambda_{1} r}\pp_{\mu}\left[\langle 2b\varphi^{2}+\int_{\mfo }{\langle \varphi,\nu\rangle}^{2}H(\cdot,{ \ud \nu}),X_{r}\rangle\right]{ \ud r}\\
&=2c^{2}_{1}\int_{t}^{t+s}\e^{-2\lambda_{1} r}\pp_{\mu}\left[\langle \vartheta[\varphi],X_{r}\rangle\right]{ \ud r}=2c^{2}_{1}\int_{t}^{t+s}\e^{-2\lambda_{1} r}\langle T_{r}(\vartheta[\varphi]),\mu\rangle { \ud r}\\
&=2c^{2}_{1}\left(\mathrm{Var}_{\mu}[{ W^{\varphi}_{t+s}}]-\mathrm{Var}_{\mu}[{ W^{\varphi}_{t}}]\right),
\end{align*}
for all $t,s\ge 0$. Recall that $\lim_{t\to+\infty}\mathrm{Var}_{\mu}[{ W^{\varphi}_{t}}]=\langle\Theta,\mu\rangle$.
It follows that
$$\lim_{t\to +\infty}\pp_{\mu}\left[\left(\e^{-\lambda_{1} (t+s)}I^{t+s}_{t,t+s} \right)^{2}\right]=0\mbox{ and } \sum_{n=1}^{+\infty}\pp_{\mu}\left[\left(\e^{-\lambda_{1} (m+n)\sigma}I^{(m+n)\sigma}_{n\sigma,(m+n)\sigma} \right)^{2}\right]<+\infty,$$
which in turn implies that
\begin{equation*}
\lim_{t\to +\infty}\e^{-\lambda_{1} (t+s)}I^{t+s}_{t,t+s} =0\mbox{ in }L^{2}(\pp_{\mu}),
\mbox{ and }\lim_{\mathbb{N}\ni n\to+\infty}\e^{-\lambda_{1} (m+n)\sigma}I^{(m+n)\sigma}_{n\sigma,(m+n)\sigma} =0\quad
\pp_{\mu}\mbox{-a.s.}
\end{equation*}
This together with \eqref{lem7.3.1} yields the desired result.
\end{proof}

\begin{proposition}\label{lem2.2}
Suppose $f_{t}$ and $f$ are $\mathbb{C}$-valued bounded measurable functions on $E$ and $f_{t}\stackrel{b}{\to}f$ as $t\to+\infty$, then the following hold.
\begin{description}
\item{(i)} For every $\mu\in\mf$,
$$\e^{-\lambda_{1} t}\langle f_{t},X_{t}\rangle \to \langle f,\widetilde{\varphi}\rangle { W^{\varphi}_{\infty}}\mbox{ in }L^{2}(\pp_{\mu})\mbox{ as }t\to+\infty.$$
\item{(ii)} For every $\sigma>0$ and $\mu\in\mf$,
$$\lim_{\mathbb{N}\ni n\to+\infty}\e^{-\lambda_{1} n\sigma}\langle f_{n\sigma},X_{n\sigma}\rangle=\langle f,\widetilde{\varphi}\rangle { W^{\varphi}_{\infty}}\quad\pp_{\mu}\mbox{-a.s.}$$
\end{description}
\end{proposition}

\begin{proof}
We may write $f_{t}=\re f_{t}+i\im f_{t}$ and $f=\re f+i\im f$. Then $\re f_{t}$, $\re f$, $\im f_{t}$, $\im f$ are $\R$-valued bounded measurable functions with $\re f_{t}\stackrel{b}{\to}\re f$ and $\im f_{t}\stackrel{b}{\to}f$. We have
\begin{eqnarray*}
\e^{-\lambda_{1} t}\langle f_{t},X_{t}\rangle -\langle f,\widetilde{\varphi}\rangle { W^{\varphi}_{\infty}}
&=& \e^{-\lambda_{1} t}\langle \re f_{t},X_{t}\rangle -\langle \re f,\widetilde{\varphi}\rangle { W^{\varphi}_{\infty}}\\
&&\quad+i\left(\e^{-\lambda_{1} t}\langle \im f_{t},X_{t}\rangle -\langle \im f,\widetilde{\varphi}\rangle { W^{\varphi}_{\infty}}\right).
\end{eqnarray*}
Thus it suffices to prove the convergence results in (i) and (ii) hold for all $\R$-valued bounded measurable functions $f_{t},f$ with $f_{t}\stackrel{b}{\to}f$.
This follows immediately from Lemma \ref{lem2.0}, Lemma \ref{lem2.1} and the fact that for all $t,s\ge 0$,
\begin{align*}
&\e^{-\lambda_{1}(t+s)}\langle f_{t+s},X_{t+s}\rangle\\
=&\e^{-\lambda_{1}(t+s)}\left(\langle f_{t+s},X_{t+s}\rangle-\pp_{\mu}\left[\langle f_{t+s},X_{t+s}\rangle\,|\,\mathcal{F}_{t}\right]\right)+\e^{-\lambda_{1}(t+s)}\pp_{\mu}\left[\langle f_{t+s},X_{t+s}\rangle\,|\,\mathcal{F}_{t}\right].
\end{align*}
\end{proof}

As a direct consequence of Proposition \ref{lem2.2}, we obtain Theorem \ref{cor:wlln}.

\section{Martingale functional central limit theorem}\label{sec3}
To begin this section, we review some basic results on the convergence of stochastic processes. For further details, we refer the reader to \cite{JS}.

Recall that $\mathcal{D}[0,+\infty)$ is the space of c\`{a}dl\`{a}g functions from $[0,+\infty)$ to $\R$ equipped with the Skorokhod topology and $\mathcal{P}(\mathcal{D}[0,+\infty))$ is the space of probability measures on $\mathcal{D}[0,+\infty)$. A subset $A$ of $\mathcal{P}(\mathcal{D}[0,+\infty))$ is called tight if for very $\epsilon>0$, there is a compact subset $K$ of $\mathcal{D}[0,+\infty)$ such that $\rho(\mathcal{D}[0,+\infty)\setminus K)\le \epsilon$ for all $\rho\in A$. It is known that a subset $A$ of $\mathcal{P}(\mathcal{D}[0,+\infty))$ is relatively compact if and only if it is tight.
Suppose $Y^{n}$, $Y$ are c\`{a}dl\`{a}g processes. We say that $\{Y^{n}\}$ is tight (resp. C-tight) if the set of laws of $Y^{n}$ is tight (resp. tight and all limit points of the subsequences are laws of continuous processes).
It is known that $Y^{n}\stackrel{d}{\to}Y$ in $\mathcal{D}[0,+\infty)$ if and only if
\begin{description}
\item{(i)} $\{Y^{n}\}$ is tight.
\item{(ii)} There is a dense set $D\subset [0,+\infty)$ such that $(Y^{n}_{t_{1}},Y^{n}_{t_{2}},\cdots,Y^{n}_{t_{k}})^{\mathrm{T}}\stackrel{d}{\to}(Y_{t_{1}},Y_{t_{2}},\cdots,Y_{t_{k}})^{\mathrm{T}}$ in $\R^{k}$ for all $t_{1},t_{2},\cdots,t_{k}\in D$.
\end{description}
Moreover, the dense set $D$ in (ii) can be taken to be $[0,+\infty)\setminus J(Y)$, where the set $J(Y):=\{t\ge 0:\ \mathrm{P}\left(\triangle Y_{t}\not=0\right)>0\}$ is at most countable.

Based on the above argument, the proof of Theorem \ref{them1} is divided into the two lemmas stated below.
\begin{lemma}\label{lem:finite dimensional convergence}
Suppose $k\in\{0,1,\cdots\}$ and $0=s_{0}\le s_{1}\le\cdots\le s_{k}<+\infty$. Then for every $x\in E$,
$$
\left(Y^{t}_{s_{0}},Y^{t}_{s_{1}},\cdots,Y^{t}_{s_{k}}\right)^{\mathrm{T}}
\stackrel{d}{\to}\sigma_{\varphi}\sqrt{W^{\varphi}_{\infty}}\,(G_{s_{0}},G_{s_{1}},\cdots,G_{s_{k}})^{\mathrm{T}}
\mbox{ in }\R^{k+1}\mbox{ as }t\to+\infty,
$$
where $(G_{s_{0}},G_{s_{1}},\cdots,G_{s_{k}})$ is an $\R^{k+1}$-valued Gaussian random variable with mean $0$ and covariance
$$\mathrm{Cov}(G_{s_{i}},G_{s_{j}})=\e^{-\frac{\lambda_{1}}{2}(s_{j}-s_{i})}\quad\forall 0\le i\le j\le k,$$
and is independent of ${ W^{\varphi}_{\infty}}$.
\end{lemma}

\begin{lemma}\label{lem:tightness}
{Suppose $\varphi$ in continuous. Then }$\{Y^{t}:t\ge 0\}$ is tight.
\end{lemma}

\noindent\textbf{Proof of Theorem \ref{them1}:}
The first conclusion is a direct result of Lemma \ref{lem:finite dimensional convergence}.
Moreover, by Lemmas \ref{lem:finite dimensional convergence} and \ref{lem:tightness},
we have
for every $x\in E$, under $\pp_{\delta_{x}},$
$$Y^{t}\stackrel{d}{\to }\sigma_{\varphi}\sqrt{{ W^{\varphi}_{\infty}}}G\mbox{ in }\mathcal{D}[0,+\infty)\mbox{ as }t\to+\infty,$$
where $G:=(G_{s})_{s\ge 0}$ is centered Gaussian process with covariance
$$Cov(G_{s},G_{t})=\e^{-\frac{\lambda_{1}}{2}|t-s|}\quad\forall s,t\ge 0,$$
and $G$ is independent of ${ W^{\varphi}_{\infty}}$.
We note that for every $t,s\ge 0$, $G_{t}-G_{s}\sim N(0,\mathrm{Var}[G_{t}-G_{s}])$, where
\begin{equation*}
\mathrm{Var}[G_{t}-G_{s}]=\mathrm{Var}[G_{t}]+\mathrm{Var}[G_{s}]-2\mathrm{Cov}(G_{t},G_{s})
=2 \left(1-\e^{-\frac{\lambda_{1}}{2}|t-s|}\right).
\end{equation*}
Thus
$$\mathrm{E}\left[|G_{t}-G_{s}|^{4}\right]=3\mathrm{E}\left[|G_{t}-G_{s}|^{2}\right]^{2}=3\mathrm{Var}[G_{t}-G_{s}]^{2}=12 \left(1-\e^{-\frac{\lambda_{1}}{2}|t-s|}\right)^{2}\le 3\lambda_{1}^{2}|t-s|^{2}.$$
By the Kolmogorov-Chentsov criterion, the process $(G_{s})_{s\ge 0}$ has a continuous version.
{Hence we complete the proof for the second conclusion.}\qed

\subsection{Proof of Lemma \ref{lem:finite dimensional convergence}}\label{sec:CLT}
{Our proof of Lemma \ref{lem:finite dimensional convergence} follows two main steps: we first establish the convergence in one dimensional distribution, and then extend it to all finite dimensions via an inductive argument.}
{\subsubsection{One dimensional convergence}}
Let $Q_{t}(\mu,{ \ud \nu})$ denote the transition semigroup of the $(P_{t},\psi)$-superprocess. Obviously by \eqref{eq1}, $Q_{t}$ satisfies that
$$\int_{\mf}\e^{-{\langle f,\nu\rangle}}Q_{t}(\mu,{ \ud \nu})=\e^{-\langle V_{t}f,\mu\rangle}\quad\forall \mu\in\mf,\ f\in\mathcal{B}^{+}_{b}(E).$$
where $V_{t}f(x)=-\ln \int_{\mf}\e^{-{\langle f,\nu\rangle}}Q_{t}(\delta_{x},{ \ud \nu})$. The semigroup $(Q_{t})_{t\ge 0}$ is said to have the \textit{regular branching property}.
This property implies that $Q_{t}(\mu,{ \ud \nu})$ is an infinitely divisible distribution on $\mf$, and hence $X_{t}$, under $\pp_{\mu}$, is an infinitely divisible random measure.
Consequently, ${ W^{\varphi}_{t}}=\e^{-\lambda_{1} t}\langle \varphi,X_{t}\rangle$ is an infinitely divisible random variable. Since ${ W^{\varphi}_{\infty}}$ is the $L^{2}(\pp_{\mu})$-limit of ${ W^{\varphi}_{t}}$, it follows by \cite[Lemma 7.8]{Sato} that ${ W^{\varphi}_{\infty}}$ is also an infinitely divisible random variable.

For every $x\in E$, let $\phi(x,\theta)$ be the characteristic exponent of the distribution of ${ W^{\varphi}_{\infty}}$ under $\pp_{\delta_{x}}$, that is,
$$\e^{-\phi(x,\theta)}=\pp_{\delta_{x}}\left[\e^{i\theta { W^{\varphi}_{\infty}}}\right]\quad\forall \theta\in\R.$$
By the regular branching property, it holds that for every $\mu\in\mf$,
\begin{equation}\label{eq:charac}
\pp_{\mu}\left[\e^{i\theta { W^{\varphi}_{\infty}}}\right]=\e^{-\langle \phi(\cdot,\theta),\mu\rangle}.
\end{equation}
Since ${ W^{\varphi}_{\infty}}$ is nonnegative and infinitely divisible, $\phi(x,\theta)$ can be represented by
\begin{equation}\label{eq:Phi(x,theta)}
\phi(x,\theta)=-i d(x)\theta+\int_{(0,+\infty)}\left(1-\e^{i\theta r}\right)\Upsilon(x,{ \ud r}),
\end{equation}
where $d(x)\ge 0$ and $(1\wedge r)\Upsilon(x,{ \ud r})$ is a bounded measure on $(0,+\infty)$.
Moreover, by \eqref{eq:charac} we have
\begin{align*}
\pp_{\delta_{x}}\left[{ W^{\varphi}_{\infty}}\right]&=-i\left.\frac{\partial}{\partial \theta}\e^{-\phi(x,\theta)}\right|_{\theta=0}=d(x)+\int_{(0,+\infty)}r\Upsilon(x,{ \ud r}),\\
\pp_{\delta_{x}}\left[ (W^{\varphi}_{\infty})^{2}\right]&=-\left.\frac{\partial^{2}}{\partial \theta^{2}}\e^{-\phi(x,\theta)}\right|_{\theta=0}=\int_{(0,+\infty)}r^{2}\Upsilon(x,{ \ud r})+\pp_{\delta_{x}}\left[{ W^{\varphi}_{\infty}}\right]^{2}.
\end{align*}
Recall that $\pp_{\delta_{x}}\left[{ W^{\varphi}_{\infty}}\right]=\varphi(x)$ and $\mathrm{Var}_{\delta_{x}}[{ W^{\varphi}_{\infty}}]=\Theta(x)$. We have
$$d(x)=\varphi(x)-\int_{(0,+\infty)}r\Upsilon(x,{ \ud r})\mbox{ and }\int_{(0,+\infty)}r^{2}\Upsilon(x,{ \ud r})=\Theta(x).$$
Inserting this to \eqref{eq:Phi(x,theta)} yields that
\begin{equation}
\phi(x,\theta)
=-i\varphi(x)\theta+\frac{1}{2}\Theta(x)\theta^{2}+\int_{(0,+\infty)}\left(1-\e^{i\theta r}+i\theta r-\frac{1}{2}\theta^{2}r^{2}\right)\Upsilon(x,{ \ud r}).\label{eq:Phi2}
\end{equation}

\begin{lemma}\label{prop1}
Suppose Assumptions \ref{AS1}-\ref{AS3} hold.
For every $x\in E$, under $\pp_{\delta_{x}}$,
$$\e^{\frac{\lambda_{1}}{2} t}\left({ W^{\varphi}_{t}}-{ W^{\varphi}_{\infty}}\right)\stackrel{d}{\to}\sigma_{\varphi}\sqrt{W^{\varphi}_{\infty}}N \mbox{ as }t\to +\infty,$$
where ${ \sigma^{2}_{\varphi}}=\langle \vartheta[\varphi],\widetilde{\varphi}\rangle/\lambda_{1}$, and $N$ is a standard normal random variable independent of ${ W^{\varphi}_{\infty}}$.
\end{lemma}

\begin{proof}
It suffices to show that for every $\theta\in\R$ and $x\in E$,
\begin{equation}\label{prop1.1}
\lim_{t\to+\infty}\pp_{\delta_{x}}\left[\exp\{i\theta \e^{\frac{\lambda_{1}}{2}t}\left({ W^{\varphi}_{t}}-{ W^{\varphi}_{\infty}}\right)\}\right]=\pp_{\delta_{x}}\left[\exp\{-\frac{1}{2}\sigma^{2}_{\varphi}\theta^{2}{ W^{\varphi}_{\infty}}\}\right].
\end{equation}

By the Markov property and \eqref{eq:charac}, we have
\begin{align}
\pp_{\delta_{x}}\left[\exp\{i\theta \e^{\frac{\lambda_{1}}{2}t}\left({ W^{\varphi}_{t}}-{ W^{\varphi}_{\infty}}\right)\}\right]
&=\pp_{\delta_{x}}\left[\exp\{i\theta \e^{\frac{\lambda_{1}}{2}t}{ W^{\varphi}_{t}}\}\pp_{X_{t}}\left[\exp\{-i\theta\e^{-\frac{\lambda_{1}}{2}t}{ W^{\varphi}_{\infty}}\}\right]\right]\nonumber\\
&=\pp_{\delta_{x}}\left[\exp\{\langle i\theta\varphi\e^{-\frac{\lambda_{1}}{2}t}-\phi(\cdot,-\theta\e^{-\frac{\lambda_{1}}{2}t}),X_{t}\rangle\}\right].\label{prop1.2}
\end{align}
By \eqref{eq:Phi2}, we have for $x\in E$,
$$i\theta\varphi(x)\e^{-\frac{\lambda_{1}}{2}t}-\phi(x,-\theta\e^{-\frac{\lambda_{1}}{2}t})=-\frac{1}{2}\theta^{2}\e^{-\lambda_{1} t}\Theta(x)+I_{\theta}(t,x),$$
where
$$I_{\theta}(t,x):=\int_{(0,+\infty)}\left(\e^{-i\theta\e^{-\frac{\lambda_{1}}{2}t}r}-1+i\theta\e^{-\frac{\lambda_{1}}{2}t}r+\frac{1}{2}\theta^{2}\e^{-\lambda_{1} t}r^{2}\right)\Upsilon(x,{ \ud r}).$$
It follows by Proposition \ref{lem2.2}(i) that
\begin{equation}\label{prop1.3}
-\frac{1}{2}\theta^{2}\e^{-\lambda_{1} t}\langle \Theta,X_{t}\rangle\to -\frac{1}{2}\theta^{2}\langle \Theta,\widetilde{\varphi}\rangle { W^{\varphi}_{\infty}}=-\frac{1}{2}\theta^{2}{ \sigma^{2}_{\varphi}}{ W^{\varphi}_{\infty}}\mbox{ in }L^{2}(\pp_{\delta_{x}})\mbox{ as }t\to+\infty.
\end{equation}
Here we use the fact that
\begin{equation}\label{eq:<Theta,widetild{Varphi}>}
\langle \Theta,\widetilde{\varphi}\rangle=\int_{0}^{+\infty}\e^{-2\lambda_{1} s}\langle T_{s}(\vartheta[\varphi]),\widetilde{\varphi}\rangle { \ud s}=\int_{0}^{+\infty}\e^{-\lambda_{1} s}\langle \vartheta[\varphi],\widetilde{\phi}\rangle { \ud s}={ \sigma^{2}_{\varphi}}.
\end{equation}
On the other hand, by the inequality (cf. \cite[Lemma A.1]{RSSZ})
\begin{equation}\label{eq:inequality2}
\left|\e^{-z}-1+z-\frac{z^{2}}{2}\right|\le|z|^{2}\left(\frac{|z|}{6}\wedge 1\right)\quad\forall z\in\mathbb{C}^{+}{:=\{z\in\mathbb{C}:\ \re z\ge 0\}},
\end{equation}
we have
\begin{align*}
|I_{\theta}(t,x)|&\le\int_{(0,+\infty)}\left|\e^{-i\theta\e^{-\frac{\lambda_{1}}{2}t}r}-1+i\theta\e^{-\frac{\lambda_{1}}{2}t}r+\frac{1}{2}\theta^{2}\e^{-\lambda_{1} t}r^{2}\right|\Upsilon(x,{ \ud r})\\
&\le\theta^{2}\e^{-\lambda_{1} t}\int_{(0,+\infty)}r^{2}\left(\frac{|\theta|\e^{-\frac{\lambda_{1}}{2}t}r}{6}\wedge 1\right)\Upsilon(x,{ \ud r})\\
&=:\theta^{2}\e^{-\lambda_{1} t}g_{\theta}(t,x).
\end{align*}
Since $\int_{(0,+\infty)}r^{2}\Upsilon(x,{ \ud r})=\Theta(x)$ is a nonnegative bounded function on $E$, one can easily show by the dominated convergence theorem that $g_{\theta}(t,\cdot)\stackrel{b}{\to} 0$ as $t\to+\infty$. Hence by Proposition \ref{lem2.2}(i), $\e^{-\lambda_{1} t}\langle g_{\theta}(t,\cdot),X_{t}\rangle \to 0$ in $L^{2}(\pp_{\delta_{x}})$, and consequently, $\langle I_{\theta}(t,\cdot),X_{t}\rangle\to 0$ in $L^{2}(\pp_{\delta_{x}})$.
This together with \eqref{prop1.3} yields that
$$-\frac{1}{2}\theta^{2}\e^{-\lambda_{1} t}\langle \Theta,X_{t}\rangle+\langle I_{\theta}(t,\cdot),X_{t}\rangle\to -\frac{1}{2}\theta^{2}{ \sigma^{2}_{\varphi}}{ W^{\varphi}_{\infty}}\mbox{ in } L^{2}(\pp_{\delta_{x}}).$$ Hence \eqref{prop1.1} follows by letting $t\to+\infty$ in \eqref{prop1.2}.
\end{proof}

{\subsubsection{Finite dimensional convergence}}
Recall $Q_{t}$ and $V_{t}$ defined in the beginning of this subsection.
Since the semigroup $(Q_{t})_{t\ge 0}$ has the regular branching property,
it follows by \cite[Theorem 2.4 and Theorem 1.36]{Li} that $(V_{t})_{t\ge 0}$ is a cumulant semigroup on $\mathcal{B}^{+}_{b}(E)$, and has the canonical representation
\begin{equation}\label{eq:canonical representation for V_{t}}
V_{t}f(x)={\int_{E}f(y)\Lambda_{t}(x,\ud y)}+\int_{\mfo }\left(1-\e^{-{\langle f,\nu\rangle}}\right)L_{t}(x,{ \ud \nu})\quad\forall f\in\mathcal{B}^{+}_{b}(E),\ t\ge 0,\ x\in E.
\end{equation}
Here $\Lambda_{t}(x,{ \ud y})$ is a bounded kernel on $E$ and $(1\wedge {\langle 1,\nu\rangle})L_{t}(x,{ \ud \nu})$ is a bounded kernel from $E$ to $\mfo $.
It follows from \eqref{eq1} and \eqref{eq:canonical representation for V_{t}} that for every $\mu\in\mf$, $t\ge 0$ and $f\in\mathcal{B}^{+}_{b}(E)$,
\begin{align}\label{eq2}
&\pp_{\mu}\left[\e^{-\langle f,X_{t}\rangle}\right]\nonumber\\
=&\exp\big\{-\big(\int_{E}{\int_{E}f(y)\Lambda_{t}(x,\ud y)}\mu({ \ud x})+\int_{E}\mu({ \ud x})\int_{\mfo }\left(1-\e^{-{\langle f,\nu\rangle}}\right)L_{t}(x,{ \ud \nu})\big)\big\}.
\end{align}
This implies that, for every $t\ge 0$ and $\mu\in\mf$, the distribution of the random measure $X_{t}$ under $\pp_{\mu}$ is equal to that of
$$\int_{E}\Lambda_{t}(x,\cdot)\mu({ \ud x})+\int_{\mf}\nu(\cdot)N^{\mu}_{t}({ \ud \nu}),$$
where $N^{\mu}_{t}({ \ud \nu})$ denotes a Poisson random measure on $\mf$ with intensity $\int_{E}L_{t}(x,{ \ud \nu})\mu({ \ud x})$.
Lemma \ref{lemA.1} in the appendix implies that, for every $\mathbb{C}^{+}$-valued bounded measurable function $f$ on E, the expectation of $\exp\{-\int_{\mf}{\langle f,\nu\rangle}N^{\mu}_{t}({ \ud \nu})\}$ exists, and is equal to $\exp\{-\int_{E}\mu({ \ud x})\int_{\mfo }\left(1-\e^{-{\langle f,\nu\rangle}}\right)L_{t}(x,{ \ud \nu})\}$. From this, one can see that, the restriction that $f$ is real-valued in \eqref{eq2} is unnecessary, and \eqref{eq2} is indeed valid for all $\mathbb{C}^{+}$-valued bounded measurable functions $f$ on $E$.

We recall that under Assumption \ref{AS3}, the second moment of $\langle f,X_{t}\rangle$ is finite for every $f\in\mathcal{B}_{b}(E)$. Using the facts that $\pp_{\delta_{x}}\left[\langle f,X_{t}\rangle\right]=\left.-\frac{\partial}{\partial \theta}\e^{-V_{t}(\theta f)(x)}\right|_{\theta=0}$ and $\pp_{\delta_{x}}\left[\langle f,X_{t}\rangle^{2}\right]=\left.\frac{\partial^{2}}{\partial \theta^{2}}\e^{-V_{t}(\theta f)(x)}\right|_{\theta=0}$, we get from \eqref{eq:canonical representation for V_{t}} and \cite[Proposition 1.38]{Li} that
$\int_{\mfo }{\langle 1,\nu\rangle}^{2}L_{t}(x,{ \ud \nu})<+\infty$, and
\begin{align}\label{eq:second moment2}
T_{t}f(x)&=\pp_{\delta_{x}}\left[\langle f,X_{t}\rangle\right]={\int_{E}f(y)\Lambda_{t}(x,\ud y)}+\int_{\mfo }{\langle f,\nu\rangle}L_{t}(x,{ \ud \nu}),\nonumber\\
\mathrm{Var}_{\delta_{x}}[\langle f,X_{t}\rangle]&=\int_{\mfo }{\langle f,\nu\rangle}^{2}L_{t}(x,{ \ud \nu}),
\end{align}
for every $x\in E$, $t\ge 0$ and $f\in\mathcal{B}^{+}_{b}(E)$.
The above identities are also valid
for $\mathbb{C}^{+}$-valued bounded measurable function $f$ on E. Inserting \eqref{eq:second moment2} to \eqref{eq2}, we get that, for every $\mathbb{C}^{+}$-valued bounded measurable function $f$,
\begin{eqnarray}\label{eq3}
\pp_{\mu}\left[\e^{-\langle f,X_{t}\rangle}\right]
&=&\exp\{-\langle T_{t}f-\frac{1}{2}\mathrm{Var}_{\delta_{\cdot}}[\langle f,X_{t}\rangle]+{ \mathfrak{L}_{t}}[f],\mu\rangle\}\nonumber\\
&=&\exp\{-\langle T_{t}f+\mathfrak{H}_{t}[f],\mu\rangle\},
\end{eqnarray}
where
$$\mathfrak{L}_{t}[f](x):=\int_{\mfo }\left(1-\e^{-{\langle f,\nu\rangle}}-{\langle f,\nu\rangle}+\frac{1}{2}{\langle f,\nu\rangle}^{2}\right)L_{t}(x,{ \ud \nu}),$$
and
$${ \mathfrak{H}_{t}}[f](x):=-\frac{1}{2}\mathrm{Var}_{\delta_{x}}[\langle f,X_{t}\rangle]+{ \mathfrak{L}_{t}}[f](x)=\int_{\mfo }\left(1-\e^{-{\langle f,\nu\rangle}}-{\langle f,\nu\rangle}\right)L_{t}(x,{ \ud \nu}).$$

\begin{lemma}\label{lem5.1}
Suppose $f_{t},f$ are $\mathbb{C}^{+}$-valued bounded measurable functions on $E$ such that $f_{t}\stackrel{b}{\to}f$ as $t\to+\infty$. Then the following hold for every $s>0$.
\begin{description}
\item{(i)} $\int_{0}^{s}T_{s-r}(\vartheta[T_{r}f_{t}]){ \ud r}\stackrel{b}{\to}\int_{0}^{s}T_{s-r}(\vartheta[T_{r}f]){ \ud r}\mbox{ as }t\to+\infty.$

\item{(ii)} $\e^{\lambda_{1} t}\mathfrak{L}_{s}[\e^{-\frac{\lambda_{1}}{2}t}f_{t}]\stackrel{b}{\to} 0 \mbox{ as } t\to +\infty.$

\item{(iii)} $\e^{\lambda_{1} t}\mathfrak{H}_{s}[\e^{-\frac{\lambda_{1}}{2}t}f_{t}]\stackrel{b}{\to}-\frac{1}{2}\int_{0}^{s}T_{s-r}(\vartheta[T_{r}f]){ \ud r}\mbox{ as }t\to+\infty$.
\end{description}
\end{lemma}

\begin{proof}
(i) Let $g_{t}$ and $g$ denote, respectively, the integrals in the left and right hand sides of (i).
By \eqref{ieq:many-to-one} we have
\begin{equation*}
|g_{t}(x)|\le\int_{0}^{s}T_{s-r}(\vartheta[T_{r}|f_{t}|])(x){ \ud r}
\le\|f_{t}\|^{2}_{\infty}\|\vartheta[1]\|_{\infty}\int_{0}^{s}\e^{c_{1}(s+r)}{ \ud r}.
\end{equation*}
Since $\{\|f_{t}\|_{\infty},t\ge t_{0}\}$ is bounded from above for some $t_{0}\ge 0$, so is $\{\|g_{t}\|_{\infty},t\ge t_{0}\}$. We only need to show that $g_{t}$ converges pointwise to $g$.
We have for all $x\in E$,
\begin{equation}\label{lem5.1.1}
|g_{t}(x)-g(x)|\le\int_{0}^{s}T_{s-r}\left(\left|\vartheta[T_{r}f_{t}]-\vartheta[T_{r}f]\right|\right)(x){ \ud r}.
\end{equation}
If
\begin{equation}\label{lem5.1.2}
\left|\vartheta[T_{r}f_{t}]-\vartheta[T_{r}f]\right|\stackrel{b}{\to}0\mbox{ as }t\to+\infty,
\end{equation}
then by \eqref{ieq:many-to-one} and the dominated convergence theorem, $T_{s-r}\left(\left|\vartheta[T_{r}f_{t}]-\vartheta[T_{r}f]\right|\right)\stackrel{b}{\to}0$, and
the integral in the right hand side of \eqref{lem5.1.1} converges to $0$, and thus $g_{t}(x)\to g(x)$ as $t\to+\infty$.
So it suffices to show \eqref{lem5.1.2}. By definition we have
\begin{align}
&\vartheta[T_{r}f_{t}](x)-\vartheta[T_{r}f](x)\nonumber\\
=&2b(x)(T_{r}f_{t}(x)^{2}-T_{r}f(x)^{2})+\int_{\mfo}\left({\langle T_{r}f_{t},\nu\rangle}^{2}-{\langle T_{r}f,\nu\rangle}^{2}\right) H(x,{ \ud \nu})\nonumber\\
=&2b(x)(T_{r}f_{t}(x)+T_{r}f(x))(T_{r}f_{t}(x)-T_{r}f(x))+\int_{\mfo}\langle T_{r}f_{t}-T_{r}f,\nu\rangle\langle T_{r}f_{t}+T_{r}f,\nu\rangle H(x,{ \ud \nu}).\label{lem5.1.3}
\end{align}
Since $f_{t}\stackrel{b}{\to}f$, it follows by \eqref{ieq:many-to-one} and the dominated convergence that $T_{r}f_{t}\stackrel{b}{\to}T_{r}f$ for every $r\ge 0$.
Recall that under Assumption \ref{AS3}, $\int_{\mfo}{\langle 1,\nu\rangle}^{2} H(x,{ \ud \nu})$ is a bounded function on $E$. Thus one can show by the dominated convergence theorem that $\int_{\mfo}{\langle 1,\nu\rangle}\langle|T_{r}f_{t}-T_{r}f|,\nu\rangle H(x,{ \ud \nu}){\to}0$ as $t\to+\infty$ for every $x\in E$.
It follows by \eqref{lem5.1.3} that
\begin{align*}
\left|\vartheta[T_{r}f_{t}](x)-\vartheta[T_{r}f](x)\right|&\le \|2b\|_{\infty}\|T_{r}f_{t}+T_{r}f\|_{\infty}\left|T_{r}f_{t}(x)-T_{r}f(x)\right|\\
&\quad+\|T_{r}f_{t}+T_{r}f\|_{\infty}\int_{\mfo}{\langle 1,\nu\rangle}\langle |T_{r}f_{t}-T_{r}f|,\nu\rangle H(x,{ \ud \nu})\to 0
\end{align*}
as $t\to+\infty$. Hence we prove \eqref{lem5.1.2}.

\smallskip

\noindent(ii) By definition,we have
$$\mathfrak{L}_{s}[\e^{-\frac{\lambda_{1}}{2}t}f_{t}](x)=\int_{\mfo }\left(1-\exp\{-\e^{-\frac{\lambda_{1}}{2}t}{\langle f_{t},\nu\rangle}\}-\e^{-\frac{\lambda_{1}}{2}t}{\langle f_{t},\nu\rangle}+\frac{1}{2}\e^{-\lambda_{1} t}{\langle f_{t},\nu\rangle}^{2}\right)L_{s}(x,{ \ud \nu}),$$
Using \eqref{eq:inequality2}, we have
$$\left|1-\exp\{-\e^{-\frac{\lambda_{1}}{2}t}{\langle f_{t},\nu\rangle}\}-\e^{-\frac{\lambda_{1}}{2}t}{\langle f_{t},\nu\rangle}+\frac{1}{2}\e^{-\lambda_{1} t}{\langle f_{t},\nu\rangle}^{2}\right|\le \e^{-\lambda_{1} t}|{\langle f_{t},\nu\rangle}|^{2}\left(\frac{\e^{-\frac{\lambda_{1}}{2}t}|{\langle f_{t},\nu\rangle}|}{6}\wedge 1\right).$$
Thus
$$|\e^{\lambda_{1} t}\mathfrak{L}_{s}[\e^{-\frac{\lambda_{1}}{2}t}f_{t}](x)|\le\int_{\mfo }|{\langle f_{t},\nu\rangle}|^{2}\left(\frac{\e^{-\frac{\lambda_{1}}{2}t}|{\langle f_{t},\nu\rangle}|}{6}\wedge 1\right)L_{s}(x,{ \ud \nu}).$$
Recall that $\int_{\mfo }{\langle 1,\nu\rangle}^{2}L_{s}(x,{ \ud \nu})<+\infty$ and $\{\|f_{t}\|_{\infty}:\ t\ge t_{0}\}$ is bounded from above. One can show by the dominated convergence theorem that the integral in the right hand side converges boundedly and pointwise to $0$ as $t\to+\infty$, and hence one gets (ii).

\smallskip

\noindent(iii) By definition we have
\begin{eqnarray*}
\mathfrak{H}_{s}[\e^{-\frac{\lambda_{1}}{2}t}f_{t}](x)
&=&\mathfrak{L}_{s}[\e^{-\frac{\lambda_{1}}{2}t}f_{t}](x)-\frac{1}{2}\e^{-\lambda_{1} t}\mathrm{Var}_{\delta_{x}}[\langle f_{t},X_{s}\rangle]\\
&=&\mathfrak{L}_{s}[\e^{-\frac{\lambda_{1}}{2}t}f_{t}](x)-\frac{1}{2}\e^{-\lambda_{1} t}\int_{0}^{s}T_{s-r}(\vartheta[T_{r}f_{t}])(x){ \ud r}.
\end{eqnarray*}
Hence (iii) is a direct result of (i) and (ii).
\end{proof}

\begin{lemma}\label{lem5.2}
Suppose $k\in \{0,1,2,\cdots\}$, $0=s_{0}\le s_{1}\le s_{2}\le\cdots\le s_{k}<+\infty$ and $\{\theta_{0},\theta_{1},\cdots,\theta_{k-1}\}\subset \R$. If $F_{t},f$ are $\mathbb{C}^{+}$-valued bounded measurable functions such that $\e^{\lambda_{1} t}F_{t}-f\stackrel{b}{\to}0$ as $t\to+\infty$, then for every $x\in E$,
\begin{multline}\label{lem5.2.0}
\pp_{\delta_{x}}\big[\exp\big\{i\theta_{0}\e^{\frac{\lambda_{1}}{2}t}{ W^{\varphi}_{t}}+i\theta_{1}\e^{\frac{\lambda_{1}}{2}(t+s_{1})}{ W^{\varphi}_{t+s_{1}}}+\cdots
+i\theta_{k-1}\e^{\frac{\lambda_{1}}{2}(t+s_{k-1})}W^{\varphi}_{t+s_{k-1}}\\
-i\left(\theta_{0}+\theta_{1}\e^{\frac{\lambda_{1}}{2}s_{1}}+\cdots+\theta_{k-1}\e^{\frac{\lambda_{1}}{2}s_{k-1}}\right)\e^{\frac{\lambda_{1}}{2}t}{ W^{\varphi}_{t+s_{k}}}-\langle F_{t},X_{t+s_{k}}\rangle\big\}\big]\\
\longrightarrow \pp_{\delta_{x}}\big[\exp\big\{-{ W^{\varphi}_{\infty}}\big[\frac{1}{2}\sigma^{2}_{\varphi}\big(\sum_{j=0}^{k-1}(1-\e^{-\lambda_{1}(s_{k}-s_{j})})\theta_{j}^{2}\\
+2\sum_{0\le i<j\le k-1}(1-\e^{-\lambda_{1}(s_{k}-s_{j})})\e^{-\frac{\lambda_{1}}{2}(s_{j}-s_{i})}\theta_{i}\theta_{j}\big)+
\e^{\lambda_{1} s_{k}}\langle f,\widetilde{\varphi}\rangle\big]\big\}\big]
\end{multline}
as $t\to+\infty$.
\end{lemma}

\begin{proof}
We prove this lemma by induction. When $k=0$, one has
\begin{eqnarray*}
\mbox{LHS of }\eqref{lem5.2.0}&=&\pp_{\delta_{x}}\left[\exp\{-\langle F_{t},X_{t}\rangle\}\right]\\
&=&\pp_{\delta_{x}}\left[\exp\{-\e^{-\lambda_{1} t}\langle \e^{\lambda_{1} t}F_{t}-f,X_{t}\rangle-\e^{-\lambda_{1} t}\langle f,X_{t}\rangle\}\right].
\end{eqnarray*}
 Given that $\e^{\lambda_{1} t}F_{t}-f\stackrel{b}{\to}0$, by Proposition \ref{lem2.2}(i) we have $\e^{-\lambda_{1} t}\langle \e^{\lambda_{1} t}F_{t}-f,X_{t}\rangle\to 0$ in $L^{2}(\pp_{\delta_{x}})$ and
$\e^{-\lambda_{1} t}\langle f,X_{t}\rangle\to \langle f,\widetilde{\varphi}\rangle { W^{\varphi}_{\infty}}$ in $L^{2}(\pp_{\delta_{x}})$. Hence by the dominated convergence theorem,
$$\mbox{LHS of }\eqref{lem5.2.0}\to \pp_{\delta_{x}}\left[\exp\{-{ W^{\varphi}_{\infty}}\langle f,\widetilde{\varphi}\rangle\}\right]\mbox{ as }t\to+\infty.$$
This shows \eqref{lem5.2.0} is valid for $k=0$.

Now we suppose \eqref{lem5.2.0} holds for $k=n\in \{0,1,\cdots\}$. We shall show that it also holds for $k=n+1$.
When $k=n+1$, we have by the Markov property
\begin{align}
&\mbox{LHS of }\eqref{lem5.2.0}\\
=&\pp_{\delta_{x}}\big[\exp\big\{i\theta_{0}\e^{\frac{\lambda_{1}}{2}t}{ W^{\varphi}_{t}}+i\theta_{1}\e^{\frac{\lambda_{1}}{2}(t+s_{1})}{ W^{\varphi}_{t+s_{1}}}+\cdots
+i\theta_{n}\e^{\frac{\lambda_{1}}{2}(t+s_{n})}{ W^{\varphi}_{t+s_{n}}}\nonumber\\
&\quad-i\big(\theta_{0}+\theta_{1}\e^{\frac{\lambda_{1}}{2}s_{1}}+\cdots+\theta_{n}\e^{\frac{\lambda_{1}}{2}s_{n}}\big)\e^{\frac{\lambda_{1}}{2}t}{ W^{\varphi}_{t+s_{n+1}}}-\langle F_{t},X_{t+s_{n+1}}\rangle\big\}\big]\nonumber\\
=&\pp_{\delta_{x}}\big[\exp\big\{i\theta_{0}\e^{\frac{\lambda_{1}}{2}t}{ W^{\varphi}_{t}}+i\theta_{1}\e^{\frac{\lambda_{1}}{2}(t+s_{1})}{ W^{\varphi}_{t+s_{1}}}+\cdots
+i\theta_{n}\e^{\frac{\lambda_{1}}{2}(t+s_{n})}{ W^{\varphi}_{t+s_{n}}}\big\}\nonumber\\
&\quad\pp_{X_{t+s_{n}}}\big[\exp\big\{-\langle i\big(\theta_{0}+\theta_{1}\e^{\frac{\lambda_{1}}{2}s_{1}}+\cdots+\theta_{n}\e^{\frac{\lambda_{1}}{2}s_{n}}\big)\e^{-\lambda_{1} s_{n+1}-\frac{\lambda_{1}}{2}t}\varphi+F_{t},X_{s_{n+1}-s_{n}}\rangle\big\}\big].\label{lem5.2.1}
\end{align}
By \eqref{eq3} we have
\begin{align*}
&\pp_{X_{t+s_{n}}}\big[\exp\big\{-\langle i\big(\theta_{0}+\theta_{1}\e^{\frac{\lambda_{1}}{2}s_{1}}+\cdots+\theta_{n}\e^{\frac{\lambda_{1}}{2}s_{n}}\big)\e^{-\lambda_{1} s_{n+1}-\frac{\lambda_{1}}{2}t}\varphi+F_{t},X_{s_{n+1}-s_{n}}\rangle\big\}\big]\\
=&\exp\big\{-\langle T_{s_{n+1}-s_{n}}\big(i\big(\theta_{0}+\theta_{1}\e^{\frac{\lambda_{1}}{2}s_{1}}+\cdots+\theta_{n}\e^{\frac{\lambda_{1}}{2}s_{n}}\big)\e^{-\lambda_{1} s_{n+1}-\frac{\lambda_{1}}{2}t}\varphi+F_{t}\big)\\
&\quad\quad+\mathfrak{H}_{s_{n+1}-s_{n}}[i\big(\theta_{0}+\theta_{1}\e^{\frac{\lambda_{1}}{2}s_{1}}+\cdots+\theta_{n}\e^{\frac{\lambda_{1}}{2}s_{n}}\big)\e^{-\lambda_{1} s_{n+1}-\frac{\lambda_{1}}{2}t}\varphi+F_{t}],X_{t+s_{n}}\rangle\big\}\\
=&\exp\big\{-i\big(\theta_{0}+\theta_{1}\e^{\frac{\lambda_{1}}{2}s_{1}}+\cdots+\theta_{n}\e^{\frac{\lambda_{1}}{2}s_{n}}\big)\e^{\frac{\lambda_{1}}{2}t}{ W^{\varphi}_{t+s_{n}}}\\
&\quad-\langle T_{s_{n+1}-s_{n}}F_{t}+\mathfrak{H}_{s_{n+1}-s_{n}}[i\big(\theta_{0}+\theta_{1}\e^{\frac{\lambda_{1}}{2}s_{1}}+\cdots+\theta_{n}\e^{\frac{\lambda_{1}}{2}s_{n}}\big)\e^{-\lambda_{1} s_{n+1}-\frac{\lambda_{1}}{2}t}\varphi+F_{t}],X_{t+s_{n}}\rangle\big\}.
\end{align*}
In the second equality, we use the fact that $T_{s_{n+1}-s_{n}}\varphi=\e^{\lambda_{1}(s_{n+1}-s_{n})}\varphi$. Inserting this to \eqref{lem5.2.1} yields that
\begin{align}
&\mbox{LHS of }\eqref{lem5.2.0}\\
=&\pp_{\delta_{x}}\big[\exp\big\{i\theta_{0}\e^{\frac{\lambda_{1}}{2}t}{ W^{\varphi}_{t}}+i\theta_{1}\e^{\frac{\lambda_{1}}{2}(t+s_{1})}{ W^{\varphi}_{t+s_{1}}}+\cdots
+i\theta_{n-1}\e^{\frac{\lambda_{1}}{2}(t+s_{n-1})}{ W^{\varphi}_{t+s_{n-1}}}\nonumber\\
&-i\big(\theta_{0}+\theta_{1}\e^{\frac{\lambda_{1}}{2}s_{1}}+\cdots+\theta_{n-1}\e^{\frac{\lambda_{1}}{2}s_{n-1}}\big)\e^{\frac{\lambda_{1}}{2}t}{ W^{\varphi}_{t+s_{n}}}\nonumber\\
&-\langle T_{s_{n+1}-s_{n}}F_{t}+\mathfrak{H}_{s_{n+1}-s_{n}}\big[i\big(\theta_{0}+\theta_{1}\e^{\frac{\lambda_{1}}{2}s_{1}}+\cdots+\theta_{n}\e^{\frac{\lambda_{1}}{2}s_{n}}\big)\e^{-\lambda_{1} s_{n+1}-\frac{\lambda_{1}}{2}t}\varphi+F_{t}\big],X_{t+s_{n}}\rangle\big\}\big]\nonumber\\
=:&\pp_{\delta_{x}}\big[\exp\big\{i\theta_{0}\e^{\frac{\lambda_{1}}{2}t}{ W^{\varphi}_{t}}+i\theta_{1}\e^{\frac{\lambda_{1}}{2}(t+s_{1})}{ W^{\varphi}_{t+s_{1}}}+\cdots
+i\theta_{n-1}\e^{\frac{\lambda_{1}}{2}(t+s_{n-1})}{ W^{\varphi}_{t+s_{n-1}}}\nonumber\\
&-i\big(\theta_{0}+\theta_{1}\e^{\frac{\lambda_{1}}{2}s_{1}}+\cdots+\theta_{n-1}\e^{\frac{\lambda_{1}}{2}s_{n-1}}\big)\e^{\frac{\lambda_{1}}{2}t}{ W^{\varphi}_{t+s_{n}}}\nonumber\\
&-\langle \widetilde{F}_{t},X_{t+s_{n}}\rangle\big\}\big].\label{lem5.2.2}
\end{align}
We note that
$$\e^{\lambda_{1} t}T_{s_{n+1}-s_{n}}F_{t}=T_{s_{n+1}-s_{n}}\left(\e^{\lambda_{1} t}F_{t}-f\right)+T_{s_{n+1}-s_{n}}f.$$
Since $\e^{\lambda_{1} t}F_{t}-f\stackrel{b}{\to}0$, it follows by \eqref{ieq:many-to-one} and the dominated convergence that \[T_{s_{n+1}-s_{n}}\left(\e^{\lambda_{1} t}F_{t}-f\right)\stackrel{b}{\to}0\mbox{ as }t\to+\infty.\]
Hence
\begin{equation}\label{lem5.2.4}
\e^{\lambda_{1} t}T_{s_{n+1}-s_{n}}F_{t}-T_{s_{n+1}-s_{n}}f\stackrel{b}{\to}0\mbox{ as }t\to+\infty.
\end{equation}
On the other hand, one may write
$$i\big(\theta_{0}+\theta_{1}\e^{\frac{\lambda_{1}}{2}s_{1}}+\cdots+\theta_{n}\e^{\frac{\lambda_{1}}{2}s_{n}}\big)\e^{-\lambda_{1} s_{n+1}-\frac{\lambda_{1}}{2}t}\varphi+F_{t}=\e^{-\frac{\lambda_{1}}{2}t}g_{t},$$
where
$$g_{t}:=i\big(\theta_{0}+\theta_{1}\e^{\frac{\lambda_{1}}{2}s_{1}}+\cdots+\theta_{n}\e^{\frac{\lambda_{1}}{2}s_{n}}\big)\e^{-\lambda_{1} s_{n+1}}\varphi+\e^{-\frac{\lambda_{1}}{2}t}\big(\e^{\lambda_{1} t}F_{t}-f\big)+\e^{-\frac{\lambda_{1}}{2}t}f.$$
It is easy to see that $g_{t}$ is a $\mathbb{C}^{+}$-valued bounded function satisfying that $g_{t}\stackrel{b}{\to}g:=i\big(\theta_{0}+\theta_{1}\e^{\frac{\lambda_{1}}{2}s_{1}}+\cdots+\theta_{n}\e^{\frac{\lambda_{1}}{2}s_{n}}\big)\e^{-\lambda_{1} s_{n+1}}\varphi$. Thus by Lemma \ref{lem5.1}(iii),
\begin{align}
&\e^{\lambda_{1} t}\mathfrak{H}_{s_{n+1}-s_{n}}[\e^{-\frac{\lambda_{1}}{2}t}g_{t}]\nonumber\\
\stackrel{b}{\to}&-\frac{1}{2}\int_{0}^{s_{n+1}-s_{n}}T_{s_{n+1}-s_{n}-r}\big(\vartheta[T_{r}g]\big){ \ud r}\nonumber\\
=&\frac{1}{2}\e^{-2\lambda_{1} s_{n+1}}\big(\theta_{0}+\theta_{1}\e^{\frac{\lambda_{1}}{2}s_{1}}+\cdots+\theta_{n}\e^{\frac{\lambda_{1}}{2}s_{n}}\big)^{2}\int_{0}^{s_{n+1}-s_{n}}\e^{2\lambda_{1} r}
T_{s_{n+1}-s_{n}-r}(\vartheta[\varphi]){ \ud r}.\label{lem5.2.5}
\end{align}
We note that $\widetilde{F}_{t}=T_{s_{n+1}-s_{n}}F_{t}+\mathfrak{H}_{s_{n+1}-s_{n}}[\e^{-\frac{\lambda_{1}}{2}t}g_{t}]$. So \eqref{lem5.2.4} and \eqref{lem5.2.5} yield that
$\e^{\lambda_{1} t}\widetilde{F}_{t}-\widetilde{f}\stackrel{b}{\to}0$, where
$$\widetilde{f}:=T_{s_{n+1}-s_{n}}f+\frac{1}{2}\e^{-2\lambda_{1} s_{n+1}}\big(\theta_{0}+\theta_{1}\e^{\frac{\lambda_{1}}{2}s_{1}}+\cdots+\theta_{n}\e^{\frac{\lambda_{1}}{2}s_{n}}\big)^{2}\int_{0}^{s_{n+1}-s_{n}}\e^{2\lambda_{1} r}
T_{s_{n+1}-s_{n}-r}(\vartheta[\varphi]){ \ud r}.$$
Hence we can apply the formula \eqref{lem5.2.0} for $k=n$ to the right hand side of \eqref{lem5.2.2} to get
\begin{multline}
\mbox{RHS of }\eqref{lem5.2.2}
\to\pp_{\delta_{x}}\big[\exp\big\{-{ W^{\varphi}_{\infty}}\big[\frac{1}{2}\sigma^{2}_{\varphi}\big(\sum_{j=0}^{n-1}(1-\e^{-\lambda_{1}(s_{n}-s_{j})})\theta^{2}_{j}
\\\quad+2\sum_{0\le i<j\le n-1}(1-\e^{-\lambda_{1}(s_{n}-s_{j})})\e^{-\frac{\lambda_{1}}{2}(s_{j}-s_{i})}\theta_{i}\theta_{j}\big)
+\e^{\lambda_{1} s_{n}}\langle \widetilde{f},\widetilde{\varphi}\rangle\big]\big\}\big].\label{lem5.2.3}
\end{multline}
Note that
\begin{align*}
&\langle \widetilde{f},\widetilde{\varphi}\rangle\\
=&\langle T_{s_{n+1}-s_{n}}f,\widetilde{\varphi}\rangle\\
&+\frac{1}{2}\e^{-2\lambda_{1} s_{n+1}}\big(\theta_{0}+\theta_{1}\e^{\frac{\lambda_{1}}{2}s_{1}}+\cdots+\theta_{n}\e^{\frac{\lambda_{1}}{2}s_{n}}\big)^{2}\int_{0}^{s_{n+1}-s_{n}}\e^{2\lambda_{1} r}
\langle T_{s_{n+1}-s_{n}-r}(\vartheta[\varphi]),\widetilde{\varphi}\rangle { \ud r}\\
=&\e^{\lambda_{1}(s_{n+1}-s_{n})}\langle f,\widetilde{\varphi}\rangle\\
&+\frac{1}{2}\e^{-2\lambda_{1} s_{n+1}}\big(\theta_{0}+\theta_{1}\e^{\frac{\lambda_{1}}{2}s_{1}}+\cdots+\theta_{n}\e^{\frac{\lambda_{1}}{2}s_{n}}\big)^{2}\int_{0}^{s_{n+1}-s_{n}}\e^{2\lambda_{1} r}\e^{\lambda_{1}(s_{n+1}-s_{n}-r)}{ \langle \vartheta[\varphi],\widetilde{\varphi}\rangle} { \ud r}\\
=&\e^{\lambda_{1}(s_{n+1}-s_{n})}\langle f,\widetilde{\varphi}\rangle+\frac{1}{2}\sigma^{2}_{\varphi}\big(\theta_{0}+\theta_{1}\e^{\frac{\lambda_{1}}{2}s_{1}}+\cdots+\theta_{n}\e^{\frac{\lambda_{1}}{2}s_{n}}\big)^{2}\e^{-2\lambda_{1} s_{n}}(1-\e^{-\lambda_{1}(s_{n+1}-s_{n})}).
\end{align*}
Inserting this to the right hand side of \eqref{lem5.2.3}, we can show that the right hand side of \eqref{lem5.2.3} is equal to that of \eqref{lem5.2.0} for $k=n+1$. Hence we have completed the proof.
\end{proof}

\noindent\textbf{Proof of Lemma \ref{lem:finite dimensional convergence}:} We need to prove that the characteristic function of the random vector $(Y^{t}_{s_{0}},\cdots,Y^{t}_{s_{k}})^{\mathrm{T}}$ converges pointwise to that of the claimed limit distribution.
Suppose $\{\theta_{0},\theta_{1},\cdots,\theta_{k}\}\subset\R$.  By the Markov property and \eqref{eq:charac}, we have
\begin{align}
&\pp_{\delta_{x}}\left[\exp\left\{i\theta_{0}Y^{t}_{s_{0}}+\cdots+i\theta_{k}
Y^{t}_{s_{k}}
\right\}\right]\nonumber\\
=&\pp_{\delta_{x}}\big[\exp\big\{i\theta_{0}\e^{\frac{\lambda_{1}}{2}t}{ W^{\varphi}_{t}}+\cdots
+i\theta_{k}\e^{\frac{\lambda_{1}}{2}(t+s_{k})}{ W^{\varphi}_{t+s_{k}}}\big\}\pp_{X_{t+s_{k}}}\big[\exp\big\{-iC(\mathbf{\theta},\mathbf{s})
\e^{-\frac{\lambda_{1}}{2}t}{ W^{\varphi}_{\infty}}\big\}\big]\big]\nonumber\\
=&\pp_{\delta_{x}}\big[\exp\big\{i\theta_{0}\e^{\frac{\lambda_{1}}{2}t}{ W^{\varphi}_{t}}+\cdots
+i\theta_{k}\e^{\frac{\lambda_{1}}{2}(t+s_{k})}{ W^{\varphi}_{t+s_{k}}}-\langle \phi(\cdot,-C(\mathbf{\theta},\mathbf{s})\e^{-\frac{\lambda_{1}}{2}t}),X_{t+s_{k}}\rangle\big\}\big],\label{prop2.1}
\end{align}
where $\mathbf{\theta}:=(\theta_{0},\cdots,\theta_{k})^{\mathrm{T}}$, $\mathbf{s}:=(s_{0},\cdots,s_{k})^{\mathrm{T}}$ and $C(\mathbf{\theta},\mathbf{s}):=(\theta_{0}+\theta_{1}\e^{\frac{\lambda_{1}}{2}s_{1}}+\cdots+\theta_{k}\e^{\frac{\lambda_{1}}{2}s_{k}})\e^{-\lambda_{1} s_{k}}$.
By \eqref{eq:Phi2} one may write $\phi(x,-C(\mathbf{\theta},\mathbf{s})\e^{-\frac{\lambda_{1}}{2}t})$ as
$$i C(\mathbf{\theta},\mathbf{s})\e^{-\frac{\lambda_{1}}{2}t}\varphi(x)+F_{\mathbf{\theta},\mathbf{s}}(t,x),$$
where
$$F_{\mathbf{\theta},\mathbf{s}}(t,x):=\frac{1}{2}C(\mathbf{\theta},\mathbf{s})^{2}\e^{-\lambda_{1} t}\Theta(x)+I_{\mathbf{\theta},\mathbf{s}}(t,x),$$
and
$$I_{\mathbf{\theta},\mathbf{s}}(t,x):=\int_{(0,+\infty)}\left(1-\e^{-i C(\mathbf{\theta},\mathbf{s})\e^{-\frac{\lambda_{1}}{2}t}r}-i C(\mathbf{\theta},\mathbf{s})\e^{-\frac{\lambda_{1}}{2}t}r-\frac{1}{2}C(\mathbf{\theta},\mathbf{s})^{2}\e^{-\lambda_{1} t}r^{2}\right)\Upsilon(x,{ \ud r}).$$
Then one has
$$\langle \phi(\cdot,-C(\mathbf{\theta},\mathbf{s})\e^{-\frac{\lambda_{1}}{2}t}),X_{t+s_{k}}\rangle=i C(\mathbf{\theta},\mathbf{s})\e^{\frac{\lambda_{1}}{2}t+\lambda_{1} s_{k}}{ W^{\varphi}_{t+s_{k}}}+\langle F_{\mathbf{\theta},\mathbf{s}}(t,\cdot),X_{t+s_{k}}\rangle.$$
Inserting this to the right hand side of \eqref{prop2.1} yields that
\begin{align}
\mbox{RHS of }\eqref{prop2.1}
&=\pp_{\delta_{x}}\big[\exp\big\{i\theta_{0}\e^{\frac{\lambda_{1}}{2}t}{ W^{\varphi}_{t}}+\cdots
+i\theta_{k-1}\e^{\frac{\lambda_{1}}{2}(t+s_{k-1})}W^{\varphi}_{t+s_{k-1}}\nonumber\\
&\quad-i(\theta_{0}+\theta_{1}\e^{\frac{\lambda_{1}}{2}s_{1}}+\cdots+\theta_{k-1}\e^{\frac{\lambda_{1}}{2}s_{k-1}})\e^{\frac{\lambda_{1}}{2}t}{ W^{\varphi}_{t+s_{k}}}-\langle  F_{\mathbf{\theta},\mathbf{s}}(t,\cdot),X_{t+s_{k}}\rangle\big\}\big].\label{prop2.2}
\end{align}
It follows by \eqref{eq:inequality2} that for $t\ge 0$ and $x\in E$,
\begin{equation*}
|I_{\mathbf{\theta},\mathbf{s}}(t,x)|\le C(\mathbf{\theta},\mathbf{s})^{2}\e^{-\lambda_{1} t}
\int_{(0,+\infty)}r^{2}\left(\frac{|C(\mathbf{\theta},\mathbf{s})|\e^{-\frac{\lambda_{1}}{2}t}r}{6}\wedge 1\right)\Upsilon(x,{ \ud r}).
\end{equation*}
Since $\int_{(0,+\infty)}r^{2}\Upsilon(x,{ \ud r})=\Theta(x)$ is a nonnegative bounded function on $E$,
one can verify by the above inequality that $\e^{\lambda_{1} t}I_{\mathbf{\theta},\mathbf{s}}(t,x)\stackrel{b}{\to}0$ as $t\to+\infty$. Hence we have
$$\e^{\lambda_{1} t}F_{\mathbf{\theta},\mathbf{s}}(t,x)-\frac{1}{2}C(\mathbf{\theta},\mathbf{s})^{2}\Theta(x)=\e^{\lambda_{1} t}I_{\mathbf{\theta},\mathbf{s}}(t,x)\stackrel{b}{\to}0\mbox{ as }t\to+\infty.$$
Applying Lemma \ref{lem5.2} to the right hand side of \eqref{prop2.2} we get that as $t\to+\infty$,
\begin{align}
&\mbox{RHS of }\eqref{prop2.1}\longrightarrow \pp_{\delta_{x}}\big[\exp\big\{-{ W^{\varphi}_{\infty}}\big[\frac{1}{2}\sigma^{2}_{\varphi}\big(\sum_{j=0}^{k-1}(1-\e^{-\lambda_{1}(s_{k}-s_{j})})\theta_{j}^{2}\nonumber\\
&+2\sum_{0\le i<j\le k-1}(1-\e^{-\lambda_{1}(s_{k}-s_{j})})\e^{-\frac{\lambda_{1}}{2}(s_{j}-s_{i})}\theta_{i}\theta_{j}\big)
+\frac{1}{2}C(\mathbf{\theta},\mathbf{s})^{2}\e^{\lambda_{1} s_{k}}\langle \Theta,\widetilde{\varphi}\rangle\big]\big\}\big]\nonumber\\
=&\pp_{\delta_{x}}\left[\exp\big\{-\frac{1}{2}\sigma^{2}_{\varphi}{ W^{\varphi}_{\infty}}\big[\sum_{j=0}^{k}\theta^{2}_{j}+2\sum_{0\le i<j\le k}\e^{-\frac{\lambda_{1}}{2}(s_{j}-s_{i})}\theta_{i}\theta_{j}\big]\big\}\right].\label{prop2.1.3}
\end{align}
In the final equality we use the fact that $\langle \Theta,\widetilde{\varphi}\rangle={ \sigma^{2}_{\varphi}}$ to simplify the computation.
Since the right hand side of \eqref{prop2.1.3} is the limit of a sequence of characteristic functions, it must be the characteristic function of some $\R^{k+1}$-valued random variable. We denote such a random variable by $\eta=(\eta_{0},\cdots,\eta_{k})^{\mathrm{T}}$.  It follows that for $0\le j\le l\le k$,
$\mathrm{E}[\eta_{j}]=-i\left.\frac{\partial}{\partial \theta_{j}}\mathrm{E}\left[\e^{i\theta\cdot\eta}\right]\right|_{\theta=0}=0,$
and
$$\mathrm{Cov}(\eta_{j},\eta_{l})=\mathrm{E}[\eta_{j}\eta_{l}]=-\left.\frac{\partial^{2}}{\partial \theta_{j}\partial \theta_{l}}\mathrm{E}\left[\e^{i\theta\cdot\eta}\right]\right|_{\theta=0}={ \sigma^{2}_{\varphi}}\e^{-\frac{\lambda_{1}}{2}(s_{l}-s_{j})}\varphi(x).$$
The matrix $\Xi:=\left(\e^{-\frac{\lambda_{1}}{2}|s_{l}-s_{j}|}\right)_{0\le l,j\le k}$ is positive definite, since ${ \sigma^{2}_{\varphi}}\varphi(x)\Xi$ is the covariance matrix of $\eta$.
Let $(G_{s_{0}},\cdots,G_{s_{k}})^{\mathrm{T}}$ be an $\R^{k+1}$-valued Gaussian random variable, independent of ${ W^{\varphi}_{\infty}}$, with mean $0$ and covariance matrix $\Xi$.
It is easy to show that the characteristic function of ${ \sigma_{\varphi}}\sqrt{W^{\varphi}_{\infty}}(G_{s_{0}},\cdots,G_{s_{k}})^{\mathrm{T}}$ is equal to the right hand side of \eqref{prop2.1.3}.
Hence we have proved this result.\qed

\subsection{Proof of Lemma \ref{lem:tightness}}

{  We recall the definition of C-tightness given at the beginning of Section \ref{sec3}. Our proof of tightness follows the approach of \cite[Section 3.2]{RSZ17b}, with the key tool being \cite[Corollary VI.3.33]{JS}, which states that the sum of a tight sequence and a C-tight sequence remains tight. We decompose the process as
$$Y^{t}_{s}=Y^{t}_{1}(s)+Y^{t}_{2}(s)\quad\forall t,s\ge 0,$$
where
$$Y^{t}_{1}(s):=\e^{\frac{\lambda_{1}}{2}(t+s)}\left({ W^{\varphi}_{t+s}}-{ W^{\varphi}_{t}}\right)\mbox{ and }Y^{t}_{2}(s):=\e^{\frac{\lambda_{1}}{2}(t+s)}\left({ W^{\varphi}_{t}}-{ W^{\varphi}_{\infty}}\right)=\e^{\frac{\lambda_{1}}{2}s}Y^{t}_{0}.$$
It then suffices to show that $\{Y^{t}_{1}:\ t\ge 0\}$ is tight and that $\{Y^{t}_{2}:\ t\ge 0\}$ is C-tight.}

\begin{lemma}\label{lem5.3}
Suppose $\tau_{1},\tau_{2}$ are $\mathcal{F}$-stopping times with $0\le \tau_{1}\le \tau_{2}$. Then for every $\mu\in\mf$,
$$\pp_{\mu}\left[\left|{ W^{\varphi}_{\tau_{2}}}-{ W^{\varphi}_{\tau_{1}}}\right|^{2}\right]\le\pp_{\mu}\left[\int_{\tau_{1}}^{\tau_{2}}\e^{-2\lambda_{1} s}\langle \vartheta[\varphi],X_{s}\rangle { \ud s}\right].$$
\end{lemma}

\begin{proof}
By \eqref{eq:stochastic integral representation}, the martingale ${ W^{\varphi}_{t}}$ has the following representation. For every $t\ge 0$ and $\mu\in\mf$,
\begin{equation}
{ W^{\varphi}_{t}}=\langle \varphi,X_{0}\rangle+M^{c}_{t}+M^{d}_{t}\quad\pp_{\mu}\mbox{-a.s.},
\end{equation}
where $M^{c}_{t}$ is a continuous square-integrable martingale with predictable quadratic variation
$$\langle M^{c}\rangle_{t}=\int_{0}^{t}\e^{-2\lambda_{1} s}\langle 2b\varphi^{2},X_{s}\rangle { \ud s},$$
and
$M^{d}_{t}=\int_{0}^{t}\int_{\mfo }\e^{-\lambda_{1} s}{\langle \varphi,\nu\rangle}\widetilde{N}({ \ud s},{ \ud \nu})$
is a purely-discontinuous square integrable martingale, whose jumps are indistinguishable from the process $\e^{-\lambda_{1} s}{\langle \varphi,\triangle X_{s}\rangle}1_{\{\triangle X_{s}\not=0\}}$.
Let $[W^{\varphi}]_{t}$ denote the quadratic variation process of ${ W^{\varphi}_{t}}$, that is, $[W^{\varphi}]_{t}=({ W^{\varphi}_{0}})^{2}+\langle M^{c}\rangle_{t}+\sum_{s\le t}(\triangle M^{d}_{s})^{2}$.

Define $M^{1,2}_{t}:=W^{\varphi}_{\tau_{2}\wedge t}-W^{\varphi}_{\tau_{1}\wedge t}$ for $t\ge 0$. By \cite[Theorem 9.2 and Theorem 9.3]{HWY}, $M^{1,2}_{t}$ can be represented as a stochastic integral with respect to the martingale ${ W^{\varphi}_{t}}$ by
$$M^{1,2}_{t}=\int_{0}^{t}1_{\{\tau_{1}<r\le \tau_{2}\}}\ud W^{\varphi}_{r},$$
and thus $M^{1,2}_{t}$ is a local martingale with quadratic variation
$$[M^{1,2}]_{t}=\int_{\tau_{1}\wedge t}^{\tau_{2}\wedge t}\ud [W^{\varphi}]_{r}=\int_{\tau_{1}\wedge t}^{\tau_{2}\wedge t}\ud \langle M^{c}\rangle _{r}+\sum_{\tau_{1}\wedge t<r\le \tau_{2}\wedge t}(\triangle M^{d}_{r})^{2}.$$
Since $(M^{1,2}_{t})^{2}-[M^{1,2}]_{t}$ is a local martingale, one can select a sequence of stopping times $s_{n}$ such that $s_{n}\uparrow+\infty$ a.s., and $((M^{1,2}_{t\wedge s_{n}})^{2}-[M^{1,2}]_{t\wedge s_{n}})_{t\ge 0}$ is a (uniformly integrable) martingale. Hence one has $\pp_{\mu}\left[(M^{1,2}_{t\wedge s_{n}})^{2}\right]=\pp_{\mu}\left[[M^{1,2}]_{t\wedge s_{n}}\right]$. Letting $n\to+\infty$ we get $\pp_{\mu}\left[(M^{1,2}_{t})^{2}\right]\le\pp_{\mu}\left[[M^{1,2}]_{t}\right]$ for every $t\ge 0$. Thus by Fatou's lemma, we have
\begin{align}
\pp_{\mu}\left[\left({ W^{\varphi}_{\tau_{2}}}-{ W^{\varphi}_{\tau_{1}}}\right)^{2}\right]&=\pp_{\mu}\left[\lim_{t\to+\infty}\left(M^{1,2}_{t}\right)^{2}\right]\le\liminf_{t\to+\infty}\pp_{\mu}\left[\left(M^{1,2}_{t}\right)^{2}\right]\nonumber\\
&\le\liminf_{t\to+\infty}\pp_{\mu}\left[[M^{1,2}]_{t}\right]=\pp_{\mu}\left[\int_{\tau_{1}}^{\tau_{2}}\ud \langle M^{c}\rangle_{r}\right]+\pp_{\mu}\left[\sum_{\tau_{1}<r\le \tau_{2}}(\triangle M^{d}_{r})^{2}\right]\nonumber\\
&=\pp_{\mu}\left[\int_{\tau_{1}}^{\tau_{2}}\e^{-2\lambda_{1} s}\langle 2b\varphi^{2},X_{s}\rangle { \ud s}\right]
+\pp_{\mu}\left[\sum_{\tau_{1}<s\le \tau_{2}}\e^{-2\lambda_{1} s}{\langle \varphi,\triangle X_{s}\rangle}^{2}\right].\label{lem5.3.1}
\end{align}
We note that
\begin{align*}
\pp_{\mu}\left[\sum_{\tau_{1}<s\le \tau_{2}}\e^{-2\lambda_{1} s}{\langle \varphi,\triangle X_{s}\rangle}^{2}\right]&=\pp_{\mu}\left[\int_{\tau_{1}}^{\tau_{2}}\e^{-2\lambda_{1} s}{\langle \varphi,\nu\rangle}^{2}N({ \ud s},{ \ud \nu})\right]\\
&=\pp_{\mu}\left[\int_{\tau_{1}}^{\tau_{2}}\e^{-2\lambda_{1} s}\langle\int_{\mfo }{\langle \varphi,\nu\rangle}^{2}H(\cdot,{ \ud \nu}),X_{s}\rangle { \ud s}\right].
\end{align*}
Putting this back into the right hand side of \eqref{lem5.3.1}, we get that
\begin{equation*}
\pp_{\mu}\left[\left({ W^{\varphi}_{\tau_{2}}}-{ W^{\varphi}_{\tau_{1}}}\right)^{2}\right]\le\pp_{\mu}\left[\int_{\tau_{1}}^{\tau_{2}}\e^{-2\lambda_{1} s}\langle \vartheta[\varphi],X_{s}\rangle { \ud s}\right],
\end{equation*}
Where $\vartheta[\varphi](x)=2b(x)\varphi(x)^{2}+\int_{\mfo }{\langle \varphi,\nu\rangle}^{2}H(x,{ \ud \nu})$.
\end{proof}

\noindent\textbf{Proof of Lemma \ref{lem:tightness}:}
\noindent{\it Step 1.} We shall prove $\{Y^{t}_{2}:=(Y^{t}_{2}(s))_{s\ge 0}:\ t\ge 0\}$ is C-tight.

First we show that $\{Y^{t}_{2}: t\ge 0\}$ is tight. By \cite[Theorem VI.3.21]{JS}, we have to show the following hold for every $x\in E$.
\begin{description}
\item{(i)} For any $N>0$,
$$\lim_{K\to+\infty}\limsup_{t\to+\infty}\pp_{\delta_{x}}\left(\sup_{s\in [0,N]}|Y^{t}_{2}(s)|>K\right)=0.$$
\item{(ii)} For any $\eta>0$ and $N>0$,
$$\lim_{\theta\downarrow 0}\limsup_{t\to+\infty}\pp_{\delta_{x}}\left(w'_{N}(Y^{t}_{2},\theta)\ge \eta\right)=0,$$
where for $\alpha\in \mathcal{D}[0,+\infty)$ and $\theta>0$, $w'_{N}(\alpha,\theta):=\inf\{\max_{1\le i\le k}\sup_{s,r\in [t_{i-1},t_{i}]}|\alpha(s)-\alpha(r)|:\ 0=t_{0}<t_{1}<\cdots<t_{k}=N,\mbox{ and } \inf_{1\le i\le k}|t_{i}-t_{i-1}|\ge \theta\}$.
\end{description}
Proof of (i): We note that
\begin{equation}\label{lem4.7.1}
\pp_{\delta_{x}}\left(\sup_{s\in [0,N]}|Y^{t}_{2}(s)|>K\right)=\pp_{\delta_{x}}\left(\sup_{s\in [0,N]}|\e^{\frac{\lambda_{1}}{2}s}Y^{t}_{0}|>K\right)
=\pp_{\delta_{x}}\left(|Y^{t}_{0}|>\e^{-\frac{\lambda_{1}}{2}N}K\right).
\end{equation}
Lemma \ref{lem:finite dimensional convergence} implies that $Y^{t}_{0}$ converges in distribution to some random variable as $t\to+\infty$.
Hence (i) follows immediately from \eqref{lem4.7.1} and the tightness of $\{Y^{t}_{0}:t\ge 0\}$.

\smallskip

\noindent Proof of (ii): We note that for $s,r\in [t_{i-1},t_{i}]\subset [0,N]$,
\begin{equation*}
\left|Y^{t}_{2}(s)-Y^{t}_{2}(r)\right|=\left|\e^{\frac{\lambda_{1}}{2}s}Y^{t}_{0}-\e^{\frac{\lambda_{1}}{2}r}Y^{t}_{0}\right|
=\e^{\frac{\lambda_{1}}{2}(s\wedge r)}\left(\e^{\frac{\lambda_{1}}{2}|s-r|}-1\right)|Y^{t}_{0}|
\le\e^{\frac{\lambda_{1}}{2}N}\left(\e^{\frac{\lambda_{1}}{2}(t_{i}-t_{i-1})}-1\right)|Y^{t}_{0}|.
\end{equation*}
Hence we have for $0=t_{0}<t_{1}<\cdots<t_{k}=N$,
$$\max_{1\le i\le k}\sup_{s,r\in [t_{i-1},t_{i}]}\left|Y^{t}_{2}(s)-Y^{t}_{2}(r)\right|\le \e^{\frac{\lambda_{1}}{2}N}\left(\e^{\frac{\lambda_{1}}{2}\max_{1\le i\le k}(t_{i}-t_{i-1})}-1\right)|Y^{t}_{0}|.$$
It follows that
$$w'_{N}(Y^{t}_{2},\theta)\le \e^{\frac{\lambda_{1}}{2}N}\left(\e^{\frac{\lambda_{1}}{2}\theta}-1\right)|Y^{t}_{0}|.$$
Thus
$$\pp_{\delta_{x}}\left(w'_{N}(Y^{t}_{2},\theta)\ge \eta\right)\le \pp_{\delta_{x}}\left(\e^{\frac{\lambda_{1}}{2}N}\left(\e^{\frac{\lambda_{1}}{2}\theta}-1\right)|Y^{t}_{0}|\ge \eta\right)
=\pp_{\delta_{x}}\left(|Y^{t}_{0}|\ge\e^{-\frac{\lambda_{1}}{2}N}\eta\left(\e^{\frac{\lambda_{1}}{2}\theta}-1\right)^{-1}\right).$$
Since $Y^{t}_{0}\stackrel{d}{\to}\xi$ in $\R$ as $t\to+\infty$ for some finite random variable $\xi$, we get that
\begin{align*}
\limsup_{t\to+\infty}\pp_{\delta_{x}}\left(w'_{N}(Y^{t}_{2},\theta)\ge \eta\right)&\le\limsup_{t\to+\infty}\pp_{\delta_{x}}\left(|Y^{t}_{0}|\ge\e^{-\frac{\lambda_{1}}{2}N}\eta\left(\e^{\frac{\lambda_{1}}{2}\theta}-1\right)^{-1}\right)\\
&\le\mathrm{P}\left(|\xi|\ge \e^{-\frac{\lambda_{1}}{2}N}\eta\left(\e^{\frac{\lambda_{1}}{2}\theta}-1\right)^{-1}\right)\to 0
\end{align*}
as $\theta\downarrow 0$. Hence we have proved that $\{Y^{t}_{2}:t\ge 0\}$ is tight.

Since $Y^{t}_{0}\stackrel{d}{\to}\xi$ in $\R$ as $t\to+\infty$,
we have for every $s_{1},s_{2},\cdots,s_{k}\in [0,+\infty)$,
$$(Y^{t}_{2}(s_{1}),Y^{t}_{2}(s_{2}),\cdots,Y^{t}_{2}(s_{k}))^{\mathrm{T}}
=Y^{t}_{0}\cdot(\e^{\frac{\lambda_{1}}{2}s_{1}},\e^{\frac{\lambda_{1}}{2}s_{2}},\cdots,\e^{\frac{\lambda_{1}}{2}s_{k}})^{\mathrm{T}}
\stackrel{d}{\to}
\xi\cdot(\e^{\frac{\lambda_{1}}{2}s_{1}},\e^{\frac{\lambda_{1}}{2}s_{2}},\cdots,\e^{\frac{\lambda_{1}}{2}s_{k}})^{\mathrm{T}}.$$
It thus by \cite[Lemma VI.3.19]{JS} that if a subsequence $\{Y^{t_{k}}_{2}:k\ge 1\}$ converges in distribution in $\mathcal{D}[0,+\infty)$, then the limit must have the same law as the continuous process $(\e^{\frac{\lambda_{1}}{2}s}\xi)_{s\ge 0}$. Therefore $\{Y^{t}_{2}:t\ge 0\}$ is C-tight.

\smallskip

\noindent{\it Step 2.} We shall prove $\{Y^{t}_{1}:=(Y^{t}_{1}(s))_{s\ge 0}:\ t\ge 0\}$ is tight.

For $t,s\ge 0$, let $\mathcal{F}^{t}_{s}$ be the $\sigma$-field generated by $\{Y^{t}_{1}(r):r\le s\}$. By Aldous' criterion for tightness (cf. \cite[Theorem VI.4.5]{JS}), it suffices to show the following hold for every $x\in E$.
\begin{description}
\item{(a)} For any $N>0$,
$$\lim_{K\to+\infty}\limsup_{t\to+\infty}\pp_{\delta_{x}}\left(\sup_{s\in [0,N]}|Y^{t}_{1}(s)|>K\right)=0.$$
\item{(b)} For any $\epsilon>0$ and $N>0$,
$$\lim_{\theta\downarrow 0}\limsup_{t\to+\infty}\sup_{s,\tau\in \mathcal{T}^{t}_{N}, 0\le \tau-s\le\theta}\pp_{\delta_{x}}\left(\left|Y^{t}_{1}(\tau)-Y^{t}_{1}(s)\right|\ge \epsilon\right)=0.$$
Here $\mathcal{T}^{t}_{N}$ denotes the set of $\mathcal{F}^{t}$-stopping times that are bounded by $N$.
\end{description}

\noindent Proof of (a): By Chebyshev's inequality we have
\begin{align*}
\pp_{\delta_{x}}\left(\sup_{s\in [0,N]}|Y^{t}_{1}(s)|>K\right)=&\pp_{\delta_{x}}\left(\sup_{s\in [0,N]}\e^{\frac{\lambda_{1}}{2}(t+s)}\left|{ W^{\varphi}_{t+s}}-{ W^{\varphi}_{t}}\right|>K\right)\\
\le&\pp_{\delta_{x}}\left(\e^{\frac{\lambda_{1}}{2}(t+N)}\sup_{s\in [0,N]}\left|{ W^{\varphi}_{t+s}}-{ W^{\varphi}_{t}}\right|>K\right)\\
\le&K^{-2}\e^{\lambda_{1}(t+N)}\pp_{\delta_{x}}\left[\sup_{s\in [0,N]}\left|{ W^{\varphi}_{t+s}}-{ W^{\varphi}_{t}}\right|^{2}\right].
\end{align*}
We can easily verify that $({ W^{\varphi}_{t+s}}-{ W^{\varphi}_{t}})_{s\ge 0}$ is a $\mathcal{F}^{t}$-martingale. Hence by the $L^{p}$-maximal inequality,
\begin{equation}
\pp_{\delta_{x}}\left(\sup_{s\in [0,N]}|Y^{t}_{1}(s)|>K\right)\le 4K^{-2}\e^{\lambda_{1}(t+N)}\pp_{\delta_{x}}\left[\left(W^{\varphi}_{t+N}-{ W^{\varphi}_{t}}\right)^{2}\right].\label{them1.1}
\end{equation}
By \eqref{eq:variance} we have for every $\mu\in\mf$ and $s\ge 0$,
$$\pp_{\mu}\left[\left(W^{\varphi}_{s}-{ W^{\varphi}_{0}}\right)^{2}\right]=\pp_{\mu}\left[\left(W^{\varphi}_{s}-\langle \varphi,\mu\rangle\right)^{2}\right]
=\mathrm{Var}_{\mu}[W^{\varphi}_{s}]=\langle f_{s},\mu\rangle,$$
Where
$$f_{s}(x):=\mathrm{Var}_{\delta_{x}}[W^{\varphi}_{s}]=\int_{0}^{s}\e^{-2\lambda_{1} r}T_{r}(\vartheta[\varphi])(x){ \ud s}\quad\forall x\in E.$$
Since $x\mapsto f_{s}(x)$ is a nonnegative bounded function on $E$, it follows by Proposition \ref{lem2.2}(i) that
$$\e^{-\lambda_{1} t}\langle f_{s},X_{t}\rangle \to \langle f_{s},\widetilde{\varphi}\rangle { W^{\varphi}_{\infty}}\mbox{ in }L^{2}(\pp_{\delta_{x}})\mbox{ as }t\to+\infty.$$
for every $x\in E$. Hence we have for every $s\ge 0$ and $x\in E$,
\begin{eqnarray}\label{them1.2}
\e^{\lambda_{1}(t+s)}\pp_{\delta_{x}}\left[\left({ W^{\varphi}_{t+s}}-{ W^{\varphi}_{t}}\right)^{2}\right]
&=&\e^{\lambda_{1} s-\lambda_{1} t}\pp_{\delta_{x}}\left[\pp_{X_{t}}\left[\left(W^{\varphi}_{s}-{ W^{\varphi}_{0}}\right)^{2}\right]\right]\nonumber\\
&=&\e^{\lambda_{1} s-\lambda_{1} t}\pp_{\delta_{x}}\left[\langle f_{s},X_{t}\rangle\right]\nonumber\\
&\to&\e^{\lambda_{1} s}\langle f_{s},\widetilde{\varphi}\rangle \varphi(x)={ \sigma^{2}_{\varphi}}\left(\e^{\lambda_{1} s}-1\right)\mbox{ as }t\to+\infty.
\end{eqnarray}
In the final equality we use the fact that
$$\langle f_{s},\widetilde{\varphi}\rangle=\int_{0}^{s}\e^{-2\lambda_{1} r}\langle T_{r}(\vartheta[\varphi]),\widetilde{\varphi}\rangle { \ud r}=\int_{0}^{s}\e^{-\lambda_{1} r}\langle \vartheta[\varphi],\widetilde{\varphi}\rangle { \ud r}={ \sigma^{2}_{\varphi}}(1-\e^{-\lambda_{1} s}).$$
\eqref{them1.2} together with \eqref{them1.1} yields that
$$\limsup_{t\to+\infty}\pp_{\delta_{x}}\left(\sup_{s\in [0,N]}|Y^{t}_{1}(s)|>K\right)\le 4K^{-2}\lim_{t\to+\infty}\e^{\lambda_{1}(t+N)}\pp_{\delta_{x}}\left[\left(W^{\varphi}_{t+N}-{ W^{\varphi}_{t}}\right)^{2}\right]\le 4K^{-2}{ \sigma^{2}_{\varphi}}(\e^{\lambda_{1} N}-1).$$
Hence we get (a) by letting $K\to+\infty$.

\smallskip

\noindent Proof of (b): For $s,\tau\in\mathcal{T}^{t}_{N}$,
\begin{eqnarray*}
Y^{t}_{1}(\tau)-Y^{t}_{1}(s)&=&\e^{\frac{\lambda_{1}}{2}(t+\tau)}\left(W^{\varphi}_{t+\tau}-{ W^{\varphi}_{t+s}}\right)+\e^{\frac{\lambda_{1}}{2}(t+s)}\left(\e^{\frac{\lambda_{1}}{2}(\tau-s)}-1\right)
\left({ W^{\varphi}_{t+s}}-{ W^{\varphi}_{t}}\right)\\
&=:&I(t,s,\tau)+II(t,s,\tau).
\end{eqnarray*}
Recall that $s,\tau$ are bounded by $N$ and that $0\le \tau-s\le \theta$. Let $\epsilon>0$. We have
\begin{eqnarray*}
\pp_{\delta_{x}}\left(|II(t,s,\tau)|>\frac{\epsilon}{2}\right)&\le&\pp_{\delta_{x}}\left(\e^{\frac{\lambda_{1}}{2}(t+s)}\left(\e^{\frac{\lambda_{1}}{2}\theta}-1\right)
\left|{ W^{\varphi}_{t+s}}-{ W^{\varphi}_{t}}\right|>\frac{\epsilon}{2}\right)\\
&\le&\pp_{\delta_{x}}\left(\sup_{r\in [0,N]}\e^{\frac{\lambda_{1}}{2}(t+r)}\left|W^{\varphi}_{t+r}-{ W^{\varphi}_{t}}\right|>\frac{\epsilon}{2}\left(\e^{\frac{\lambda_{1}}{2}\theta}-1\right)^{-1}\right)\\
&=&\pp_{\delta_{x}}\left(\sup_{r\in [0,N]}|Y^{t}_{1}(r)|>\frac{\epsilon}{2}\left(\e^{\frac{\lambda_{1}}{2}\theta}-1\right)^{-1}\right).
\end{eqnarray*}
Since $\left(\e^{\frac{\lambda_{1}}{2}\theta}-1\right)^{-1}\to +\infty$ as $\theta\downarrow 0$, it follows from the above inequality and (a) that
\begin{equation}\label{them1.3}
\lim_{\theta\downarrow 0}\limsup_{t\to+\infty}\pp_{\delta_{x}}\left(|II(t,s,\tau)|>\frac{\epsilon}{2}\right)=0.
\end{equation}
On the other hand, we have by Chebyshev's inequality
\begin{eqnarray*}
\pp_{\delta_{x}}\left(|I(t,s,\tau)|>\frac{\epsilon}{2}\right)&\le&\pp_{\delta_{x}}\left(\e^{\frac{\lambda_{1}}{2}(t+N)}\left|W^{\varphi}_{t+\tau}-{ W^{\varphi}_{t+s}}\right|>\frac{\epsilon}{2}\right)\\
&\le&\frac{4}{\epsilon^{2}}\e^{\lambda_{1}(t+N)}\pp_{\delta_{x}}\left[\left(W^{\varphi}_{t+\tau}-{ W^{\varphi}_{t+s}}\right)^{2}\right].
\end{eqnarray*}
It is easy to verify that both $t+\tau$ and $t+s$ are $\mathcal{F}$-stopping times. Hence by Lemma \ref{lem5.3} we have
\begin{eqnarray}
\pp_{\delta_{x}}\left(|I(t,s,\tau)|>\frac{\epsilon}{2}\right)&\le&\frac{4}{\epsilon^{2}}\e^{\lambda_{1}(t+N)}\pp_{\delta_{x}}\left[\int_{t+s}^{t+\tau}\e^{-2\lambda_{1} r}\langle \vartheta[\varphi],X_{r}\rangle { \ud r}\right]\nonumber\\
&=&\frac{4}{\epsilon^{2}}\e^{\lambda_{1}(t+N)}\pp_{\delta_{x}}\left[\int_{s}^{\tau}\e^{-2\lambda_{1}(t+r)}\langle \vartheta[\varphi],X_{t+r}\rangle { \ud r}\right]\nonumber\\
&\le&\frac{4}{\epsilon^{2}}\e^{\lambda_{1} N}\pp_{\delta_{x}}\left[\int_{s}^{\tau}\e^{-\lambda_{1}(t+r)}\langle \vartheta[\varphi],X_{t+r}\rangle { \ud r}\right].\label{them1.5}
\end{eqnarray}
Since $s,\tau\in\mathcal{T}^{t}_{N}$ and $0\le \tau-s\le\theta$,
\begin{align}
\pp_{\delta_{x}}\left[\int_{s}^{\tau}\e^{-\lambda_{1} (t+r)}\langle \vartheta[\varphi],X_{t+r}\rangle { \ud r}\right]
&=\pp_{\delta_{x}}\big[\int_{s}^{\tau}\left(\e^{-\lambda_{1} (t+r)}\langle \vartheta[\varphi],X_{t+r}\rangle -\langle \vartheta[\varphi],\widetilde{\varphi}\rangle { W^{\varphi}_{\infty}}\right){ \ud r}\nonumber\\
&\quad+\langle \vartheta[\varphi],\widetilde{\varphi}\rangle { W^{\varphi}_{\infty}}(\tau-s)\big]\nonumber\\
&\le\int_{0}^{N}\pp_{\delta_{x}}\big[\left|\e^{-\lambda_{1} (t+r)}\langle \vartheta[\varphi],X_{t+r}\rangle -\langle \vartheta[\varphi],\widetilde{\varphi}\rangle { W^{\varphi}_{\infty}}\right|\big]{ \ud r}\nonumber\\
&\quad+\langle \vartheta[\varphi],\widetilde{\varphi}\rangle \varphi(x)\theta.\label{them1.4}
\end{align}
By Proposition \ref{lem2.2}(i),
$\e^{-\lambda_{1} t}\langle\vartheta[\varphi],X_{t}\rangle\to \langle \vartheta[\varphi],\widetilde{\varphi}\rangle { W^{\varphi}_{\infty}}$ in $L^{2}(\pp_{\delta_{x}})$ as $t\to+\infty$, and thus
$$\pp_{\delta_{x}}\big[\left|\e^{-\lambda_{1} t}\langle \vartheta[\varphi],X_{t}\rangle -\langle \vartheta[\varphi],\widetilde{\varphi}\rangle { W^{\varphi}_{\infty}}\right|\big]\to 0\mbox{ as }t\to+\infty.$$ Hence by the dominated convergence theorem, the integral in the right hand side of \eqref{them1.4} converges to $0$ as $t\to+\infty$.
Thus we get by \eqref{them1.5} and \eqref{them1.4} that
\begin{equation*}
\limsup_{t\to+\infty}\pp_{\delta_{x}}\left(|I(t,s,\tau)|>\frac{\epsilon}{2}\right)\le \frac{4}{\epsilon^{2}}\e^{\lambda_{1} N}\langle \vartheta[\varphi],\widetilde{\varphi}\rangle \varphi(x)\theta\downarrow 0\quad\mbox{ as }\theta\downarrow 0.
\end{equation*}
This together with \eqref{them1.3} yields (b).

\smallskip

\noindent{\it Step 3.} Since $Y^{t}=Y^{t}_{1}+Y^{t}_{2}$ where $\{Y_{1}^{t}:t\ge 0\}$ is tight and $\{Y^{t}_{2}:t\ge 0\}$ is C-tight, it follows by \cite[Corollary VI.3.33]{JS} that $\{Y^{t}:t\ge 0\}$ is tight.\qed

\section{Rate of convergence in the law of large numbers}\label{sec4}

We begin this section by considering an arbitrary function $g\in\mathcal{B}_{b}(E)$ such that $\langle g,\widetilde{\varphi}\rangle=0$. Under Assumption \ref{AS1}, we have
$\sup_{x\in \R}|\varphi(x)^{-1}\e^{-\lambda_{1}t}T_{t}g(x)|\to 0$ as $t\to+\infty$.
In this case, the parameter $\epsilon(g)$, defined as
$$\epsilon(g)=-\limsup_{t\to+\infty}\frac{\log \sup_{x\in E}|\varphi(x)^{-1}\e^{-\lambda_{1}t}T_{t}g(x)|}{t},$$
characterizes the exponential decay rate of $\sup_{x\in \R}|\varphi(x)^{-1}\e^{-\lambda_{1}t}T_{t}g(x)|$ as $t\to+\infty$.
In the next section, we will show that if $\epsilon(g)>\lambda_{1}/2$, then the asymptotic growth rate of the variance of $\langle g,X_{t}\rangle$ depends on $\lambda_{1}$, whereas for $\epsilon(g)\le \lambda_{1}/2$, it depends on the growth rate of $T_{t}g$. The remainder of this section is divided into three parts, each corresponding to one of the cases outlined in Theorem \ref{them2}. The proof of Theorem \ref{them2} is given at the end of this section.

{\subsection{Asymptotic behaviour of the variance}}
\begin{lemma}\label{lem1case1}
Suppose $\epsilon(g)>\lambda_{1}/2$. Then we have
$$\varphi(x)^{-1}\e^{-\frac{\lambda_{1}}{2}t}T_{t}g(x)\rightrightarrows 0\mbox{ as }t\to+\infty,$$
and
$$\varphi(x)^{-1}\e^{-\lambda_{1} t}\mathrm{Var}_{\delta_{x}}[\langle g,X_{t}\rangle]\stackrel{b}{\to}{ \rho^{2}_{g}}\mbox{ as }t\to+\infty,$$
where ${ \rho^{2}_{g}}:=\int_{0}^{+\infty}\e^{-\lambda_{1} s}\langle \vartheta[T_{s}g],\widetilde{\varphi}\rangle { \ud s}\in [0,+\infty)$.
\end{lemma}
\begin{proof}
We have
\begin{align*}
\sup_{x\in E}\big|\varphi(x)^{-1}\e^{-\frac{\lambda_{1}}{2}t}T_{t}g(x)\big|=&\e^{\frac{\lambda_{1}}{2}t}\sup_{x\in E}\big|\varphi(x)^{-1}\e^{-\lambda_{1}t}T_{t}g(x)\big|\\
=&\exp\big\{t\big(\frac{\log\sup_{x\in E}\big|\varphi(x)^{-1}\e^{-\lambda_{1}t}T_{t}g(x)\big|}{t}+\frac{\lambda_{1}}{2}\big)\big\}=:&\exp\{t\chi(t)\}.
\end{align*}
Note that $\limsup_{t\to+\infty}\chi(t)=-\epsilon(g)+\lambda_{1}/2<0$. The first conclusion follows immediately.

Let $\epsilon\in (\frac{\lambda_{1}}{2},\epsilon(g)\wedge \lambda_{1})$. We have
\begin{align}
&\varphi(x)^{-1}\e^{-\lambda_{1} t}\mathrm{Var}_{\delta_{x}}[\langle g,X_{t}\rangle]\nonumber\\
=&\varphi(x)^{-1}\e^{-\lambda_{1} t}\int_{0}^{t}T_{t-s}(\vartheta[T_{s}g])(x){ \ud s}\nonumber\\
=&\int_{0}^{t}\e^{(\lambda_{1}-2\epsilon)s}\varphi(x)^{-1}\e^{-\lambda_{1}(t-s)}T_{t-s}\left(\vartheta[\e^{-(\lambda_{1}-\epsilon)s}T_{s}g]\right)(x){ \ud s}.\label{lem1case1.1}
\end{align}
Since
$$\limsup_{t\to+\infty}\frac{\log\sup_{x\in E}|\varphi(x)^{-1}\e^{-\lambda_{1} t}T_{t}g(x)|}{t}\le -\epsilon(g)< -\epsilon,$$
there is some $s_{0}>0$ such that for all $s\ge s_{0}$, $\sup_{x\in E}|\varphi(x)^{-1}\e^{-\lambda_{1} s}T_{s}g(x)|\le \e^{-\epsilon s}$, or equivalently,
$|\e^{-(\lambda_{1}-\epsilon)s}T_{s}g(x)|\le \varphi(x)$ for all $x\in E$. On the other hand, by \eqref{ieq:many-to-one} we have for $s\in [0,s_{0}]$ and $x\in E$,
$|\e^{-(\lambda_{1}-\epsilon)s}T_{s}g(x)|\le \e^{(c_{1}-\lambda_{1}+\epsilon)s}\|g\|_{\infty}\le \e^{|c_{1}-\lambda_{1}+\epsilon|t_{0}}\|g\|_{\infty}$.
Thus we have $M:=\sup_{s\ge 0,x\in E}|\e^{-(\lambda_{1}-\epsilon)s}T_{s}g(x)|\in [0,+\infty)$. It then follows by Assumption \ref{AS1} that
for all $0\le s\le t<+\infty$ and $x\in E$,
\begin{align*}
\left|\e^{-\lambda_{1}(t-s)}T_{t-s}(\vartheta[\e^{-(\lambda_{1}-\epsilon)s}T_{s}g])(x)\right|\le&M^{2}\e^{-\lambda_{1}(t-s)}T_{t-s}(\vartheta[1])(x)\\
\le&M^{2}\varphi(x)\|\vartheta[1]\|_{\infty}\left(\sup_{r\ge 0}\triangle_{r}+\langle 1,\widetilde{\varphi}\rangle\right),
\end{align*}
and
$$\varphi(x)^{-1}\e^{-\lambda_{1}(t-s)}T_{t-s}\left(\vartheta[\e^{-(\lambda_{1}-\epsilon)s}T_{s}g]\right)(x)\stackrel{b}{\to}\e^{-2(\lambda_{1}-\epsilon)s}\langle \vartheta[T_{s}g],\widetilde{\varphi}\rangle\mbox{ as }t\to+\infty.$$
Since $\lambda_{1}-2\epsilon<0$, it then follows by the bounded convergence theorem that the integral in the right hand side of \eqref{lem1case1.1} converges boundedly and pointwise to
\[\int_{0}^{+\infty}\e^{(\lambda_{1}-2\epsilon)s}\e^{-2(\lambda_{1}-\epsilon)s}\langle\vartheta[T_{s}g],\widetilde{\varphi}\rangle { \ud s}={ \rho^{2}_{g}}\in [0,+\infty)\mbox{ as }t\to+\infty.\]
\end{proof}
Now we consider the case when $\epsilon(g)\le \lambda_{1}/2$.

\begin{lemma}\label{lem:newcase2}
Suppose $g\in\mathcal{B}_{b}(E)$ satisfies that
\begin{equation}\label{lemnew.1}
\langle g,\widetilde{\varphi}\rangle=0\mbox{ and }\frac{\e^{-(\lambda_{1}-\epsilon)t}T_{t}g}{t^{r}}\rightrightarrows g^{*}\mbox{ as }t\to+\infty,
\end{equation}
for some $g^{*}\in\mathcal{B}_{b}(E)$, $r\in [0,+\infty)$ and $\epsilon\in [0,\lambda_{1}/2]$.
\begin{description}
\item{(i)} If $\epsilon=\lambda_{1}/2$, then
$$\frac{1}{t^{1+2r}}\e^{-\lambda_{1} t}\varphi(x)^{-1}\mathrm{Var}_{\delta_{x}}[\langle g,X_{t}\rangle]\stackrel{b}{\to}\frac{\varrho^{2}_{g^{*}}}{1+2r}\mbox{ as }t\to+\infty,$$
where $\varrho^{2}_{g^{*}}=\langle \vartheta[g^{*}],\widetilde{\varphi}\rangle$.
\item{(ii)} If $\epsilon\in [0,\lambda_{1}/2)$, then
\begin{equation}\nonumber
\frac{1}{t^{2 r}}\e^{-2(\lambda_{1}-\epsilon)t}\mathrm{Var}_{\delta_{x}}[\langle g,X_{t}\rangle]\stackrel{b}{\to}\delta^{2}_{g^{*}}(x)\mbox{ as }t\to+\infty,
\end{equation}
where for every $x\in E$, $\delta^{2}_{g^{*}}(x):=\int_{0}^{+\infty}\e^{-2(\lambda_{1}-\epsilon)s}T_{s}(\vartheta[g^{*}])(x){ \ud s}\in [0,+\infty)$.
\end{description}
\end{lemma}

\begin{remark}\label{rm1}
Obviously, when \eqref{lemnew.1} holds for $g^{*}\not\equiv 0$, we have $\epsilon(g)=\epsilon$.
Moreover, by the dominated convergence theorem we have $\langle g^{*},\widetilde{\varphi}\rangle=\lim_{t\to+\infty}t^{-r}\e^{-(\lambda_{1}-\epsilon)t}\langle T_{t}g,\widetilde{\varphi}\rangle=0$, and for any $s\ge 0$, $\e^{-(\lambda_{1}-\epsilon)s}T_{s}g^{*}=\lim_{t\to+\infty}t^{-r}\e^{-(\lambda_{1}-\epsilon)(t+s)}T_{t+s}g=g^{*}$.
This implies that $\e^{(\lambda_{1}-\epsilon)t}$ is an eigenvalue of $T_{t}$ with corresponding eigenfunction $g^{*}$.
\end{remark}

\noindent\textit{Proof of Lemma \ref{lem:newcase2}:}
(i) Suppose $\epsilon=\lambda_{1}/2$. We have
\begin{align}
&\frac{1}{t^{1+2r}}\e^{-\lambda_{1} t}\varphi(x)^{-1}\mathrm{Var}_{\delta_{x}}[\langle g,X_{t}\rangle]\nonumber\\
=&\frac{1}{t}\int_{0}^{t}\varphi(x)^{-1}\e^{-\lambda_{1}(t-s)}T_{t-s}\left(\vartheta\left[\frac{\e^{-\frac{\lambda_{1}}{2}s}T_{s}g}{t^{r}}\right]\right)(x){ \ud s}\nonumber\\
=&\int_{0}^{1}\varphi(x)^{-1}\e^{-\lambda_{1} t(1-u)}T_{t(1-u)}\left(\vartheta\left[\frac{\e^{-\frac{\lambda_{1}}{2}tu}T_{tu}g}{t^{r}}\right]\right)(x){ \ud u}\nonumber\\
=&\int_{0}^{1}u^{2r}\varphi(x)^{-1}\e^{-\lambda_{1} t(1-u)}T_{t(1-u)}\left(\vartheta\left[\frac{\e^{-\frac{\lambda_{1}}{2}tu}T_{tu}g}{(tu)^{r}}\right]\right)(x){ \ud u}.\label{lem1case2.1}
\end{align}
The second equality is from change of variables.
Given that $t^{-r}\e^{-\frac{\lambda_{1}}{2}t}T_{t}g\rightrightarrows g^{*}$ as $t\to+\infty$, we can apply similar argument as in the proof of Lemma \ref{lem5.1}(i) to show that for every $u\in (0,1)$,
$$\vartheta[\frac{\e^{-\frac{\lambda_{1}}{2}tu}T_{tu}g}{(tu)^{r}}]\stackrel{b}{\to}\vartheta[g^{*}]\mbox{ as }t\to+\infty.$$
Then by Assumption \ref{AS1} and the dominated convergence theorem, we have
$$\lim_{t\to+\infty}\varphi(x)^{-1}\e^{-\lambda_{1} t(1-u)}T_{t(1-u)}\left(\vartheta\left[\frac{\e^{-\frac{\lambda_{1}}{2}tu}T_{tu}g}{(tu)^{r}}\right]\right)(x)=\langle \vartheta[g^{*}],\widetilde{\varphi}\rangle=\varrho^{2}_{g^{*}}$$
for every $u\in (0,1)$ and $x\in E$.
So, the integral in the right hand side of \eqref{lem1case2.1} converges boundedly and pointwise to $\int_{0}^{1}u^{2r}\varrho^{2}_{g^{*}}{ \ud u}=\varrho^{2}_{g^{*}}/(1+2r)$ by the dominated convergence theorem, once, there is some nonnegative integrable function $f(u)$ on $(0,1)$ such that for every $x\in E$ and $t$ sufficiently large,
\begin{equation}\label{lem1case2.2}
u^{2r}\varphi(x)^{-1}\e^{-\lambda_{1} t(1-u)}T_{t(1-u)}\left(\vartheta\left[\frac{\e^{-\frac{\lambda_{1}}{2}tu}T_{tu}g}{(tu)^{r}}\right]\right)(x)\le f(u)\quad\forall u\in (0,1).
\end{equation}
We will prove \eqref{lem1case2.2} in the rest of this proof.
Our assumption implies that for any $\delta>0$, there is $t_{0}>0$ such that for $t\ge 2t_{0}$ and $u\ge t_{0}/t$, $\|(tu)^{-r}\e^{-\lambda_{1} tu}T_{tu}g-g^{*}\|_{\infty}\le \delta$. In this case, we have by \eqref{eq1.9}
\begin{eqnarray*}
&&\varphi(x)^{-1}\e^{-\lambda_{1} t(1-u)}T_{t(1-u)}\left(\vartheta\left[\frac{\e^{-\frac{\lambda_{1}}{2}tu}T_{tu}g}{(tu)^{r}}\right]\right)(x)\\
&\le& \|\vartheta\left[\frac{\e^{-\lambda_{1} tu}T_{tu}g}{(tu)^{r}}\right]\|_{\infty}(\triangle_{t(1-u)}+\langle 1,\widetilde{\varphi}\rangle)\\
&\le& (\delta+\|g^{*}\|_{\infty})^{2}\|\vartheta[1]\|_{\infty}(\sup_{s\ge 0}\triangle_{s}+\langle 1,\widetilde{\varphi}\rangle)=:c_{1}.
\end{eqnarray*}
On the other hand, by \eqref{ieq:many-to-one} we have $\|T_{tu}g\|_{\infty}\le \e^{c_{2}t_{0}}\|g\|_{\infty}$ for $t\ge 2t_{0}$ and $u\in [0,t_{0}/t]$, and thus by \eqref{eq1.9},
\begin{eqnarray*}
&&u^{2r}\varphi(x)^{-1}\e^{-\lambda_{1} t(1-u)}T_{t(1-u)}\left(\vartheta\left[\frac{\e^{-\frac{\lambda_{1}}{2}tu}T_{tu}g}{(tu)^{r}}\right]\right)(x)\\
&=&\frac{1}{t^{2r}}\e^{-\lambda_{1} tu}\varphi(x)^{-1}\e^{-\lambda_{1} t(1-u)}T_{t(1-u)}\left(\vartheta\left[T_{tu}g\right]\right)(x)\\
&\le&\frac{1}{(2t_{0})^{2r}}\|\vartheta[T_{tu}g]\|_{\infty}(\sup_{s\ge 0}\triangle_{s}+\langle 1,\widetilde{\varphi}\rangle)\\
&\le&\frac{1}{(2t_{0})^{2r}}\e^{2c_{2}t_{0}}\|g\|^{2}_{\infty}\|\vartheta[1]\|_{\infty}(\sup_{s\ge 0}\triangle_{s}+\langle 1,\widetilde{\varphi}\rangle)=:c_{3}.
\end{eqnarray*}
Hence \eqref{lem1case2.2} holds for $f(u)=c_{1}u^{2r}+c_{3}$.

\smallskip

\noindent(ii) Suppose $\epsilon\in [0,\lambda_{1}/2)$. We have
\begin{align*}
\frac{\e^{-2(\lambda_{1}-\epsilon)t}}{t^{2r}}\mathrm{Var}_{\delta_{x}}[\langle g,X_{t}\rangle]&=\frac{\e^{-2(\lambda_{1}-\epsilon)t}}{t^{2r}}\int_{0}^{t}T_{s}(\vartheta[T_{t-s}g])(x){ \ud s}\nonumber\\
&=\int_{0}^{t}\e^{-(\lambda_{1}-2\epsilon)s}\e^{-\lambda_{1} s}T_{s}\left(\vartheta\left[\left(1-\frac{s}{t}\right)^{r}\e^{-(\lambda_{1}-\epsilon)(t-s)}\frac{T_{t-s}g}{(t-s)^{r}}\right]\right)(x){ \ud s}.
\end{align*}
Using the fact that $\frac{\e^{-(\lambda_{1}-\epsilon)t}T_{t}g(x)}{t^{r}}\rightrightarrows g^{*}\mbox{ as }t\to+\infty$, one can easily prove that for a fixed $s\ge 0$,
\begin{equation}\label{lem2case3.1}
\vartheta\left[\e^{-(\lambda_{1}-\epsilon)(t-s)}\frac{T_{t-s}g}{(t-s)^{r}}\right]\stackrel{b}{\to} \vartheta[g^{*}]\mbox{ as }t\to +\infty.
\end{equation}
Hence by the dominated convergence theorem, we have for every $s\ge 0$,
\begin{align*}
T_{s}\left(\vartheta\left[\left(1-\frac{s}{t}\right)^{r}\e^{-(\lambda_{1}-\epsilon)(t-s)}\frac{T_{t-s}g}{(t-s)^{r}}\right]\right)
=&\left(1-\frac{s}{t}\right)^{2r}T_{s}\left(\vartheta\left[\e^{-(\lambda_{1}-\epsilon)(t-s)}\frac{T_{t-s}g}{(t-s)^{r}}\right]\right)\\
\stackrel{b}{\to}&T_{s}(\vartheta[g^{*}])(x)
\end{align*}
as $t\to+\infty$. Thus, the conclusion of (ii) follows from the dominated convergence theorem, once we prove that, there exists some nonnegative integrable function $h(s)$ on $(0,+\infty)$ such that for $t$ sufficiently large, $x\in E$ and $s>0$,
\begin{equation}\label{lem2case3.2}
h_{t,x}(s):=\e^{-(\lambda_{1}-2\epsilon)s}\e^{-\lambda_{1} s}T_{s}\left(\vartheta\left[\left(1-\frac{s}{t}\right)^{r}\e^{-(\lambda_{1}-\epsilon)(t-s)}\frac{T_{t-s}g}{(t-s)^{r}}\right]\right)(x)1_{\{s\le t\}}\le h(s).
\end{equation}
By \eqref{eq1.9}, we have for $0\le s\le t<+\infty$ and $x\in E$,
\begin{eqnarray}
&&\e^{-\lambda_{1} s}T_{s}\left(\vartheta\left[\left(1-\frac{s}{t}\right)^{r}\e^{-(\lambda_{1}-\epsilon)(t-s)}\frac{T_{t-s}g}{(t-s)^{r}}\right]\right)(x)\nonumber\\
&\le&\varphi(x)\|\vartheta\left[\e^{-(\lambda_{1}-\epsilon)(t-s)}\frac{T_{t-s}g}{(t-s)^{r}}\right]\|_{\infty}\left(\sup_{r\ge 0}\triangle_{r}+\langle 1,\widetilde{\varphi}\rangle\right).\label{eq4.25}
\end{eqnarray}
It follows by \eqref{lem2case3.1} that for every $\delta>0$, there is some $t_{0}>1$, such that for $t\ge t_{0}$ and $s\in [0,t-t_{0}]$,
$$\left\|\vartheta\left[\e^{-(\lambda_{1}-\epsilon)(t-s)}\frac{T_{t-s}g}{(t-s)^{r}}\right]-\vartheta[g^{*}]\right\|\le\delta.$$
In this case,
$$h_{t,x}(s)\le\e^{-(\lambda_{1}-2\epsilon)s}\varphi(x)(\|\vartheta[g^{*}]\|_{\infty}+\delta)(\sup_{r\ge 0}\triangle_{r}+\langle 1,\widetilde{\varphi}\rangle)=:c_{1}\e^{-(\lambda_{1}-2\epsilon)s}\varphi(x).$$
On the other hand, for $t\ge t_{0}$, $s\in (t-t_{0},t]$ and $y\in E$,
\begin{eqnarray*}\left|\left(1-\frac{s}{t}\right)^{r}\e^{-(\lambda_{1}-\epsilon)(t-s)}\frac{T_{t-s}g}{(t-s)^{r}}(y)\right|&=&\frac{1}{t^{r}}\e^{-(\lambda_{1}-\epsilon)(t-s)}|T_{t-s}g(y)|\\
&\le&
\frac{1}{t^{r}_{0}}\e^{-(\lambda_{1}-\epsilon)(t-s)}
T_{t-s}|g|(y)\le \e^{c_{2}t_{0}}\|g\|_{\infty},
\end{eqnarray*}
and consequently,
$$\|\vartheta\left[\left(1-\frac{s}{t}\right)^{r}\e^{-(\lambda_{1}-\epsilon)(t-s)}\frac{T_{t-s}g}{(t-s)^{r}}\right]\|_{\infty}\le \e^{2c_{2}t_{0}}\|g\|^{2}_{\infty}\|\vartheta[1]\|_{\infty}.$$
This together with \eqref{eq4.25} yields that
$h_{t,x}(s)\le c_{3}\e^{-(\lambda_{1}-2\epsilon)s}\varphi(x)$
for some $c_{3}>0$.
Hence we prove \eqref{lem2case3.2} holds for $h(s)=(c_{1}+c_{3})\e^{-(\lambda_{1}-2\epsilon)s}\|\varphi\|_{\infty}$.\qed

{\subsection{Small branching case $\epsilon(g)>\lambda_{1}/2$}}
\begin{lemma}\label{lem4.1}
Suppose $\{f_{t}:t\ge 0\}\subset \mathcal{B}_{b}(E)$ satisfies that $f_{t}\stackrel{b}{\to} f$ as $t\to+\infty$ for some $f\in\mathcal{B}_{b}(E)$, then for every $x\in E$ and $s\ge 0$, under $\pp_{\delta_{x}}$,
\begin{align*}
&\big(\e^{-\frac{\lambda_{1}}{2}(t+s)}\left(\langle f_{t+s},X_{t+s}\rangle-\pp_{\delta_{x}}\left[\langle f_{t+s},X_{t+s}\rangle |\mathcal{F}_{t}\right]\right),\e^{\frac{\lambda_{1}}{2}(t+s)}(W^{\varphi}_{t+s}-W^{\varphi}_{\infty})\big)
^{\mathrm{T}}\\
\stackrel{d}{\longrightarrow}&\,(\rho_{s}(f)\sqrt{{ W^{\varphi}_{\infty}}}N_{1},\sigma_{\varphi}\sqrt{W^{\varphi}_{\infty}}N_{2})^{\mathrm{T}} \mbox{ as }t\to+\infty,
\end{align*}
where $\rho^{2}_{s}(f):=\int_{0}^{s}\e^{-\lambda_{1} r}\langle \vartheta[T_{r}f],\widetilde{\varphi}\rangle { \ud r}$, and $N_{1},N_{2}$ are two independent standard normal random variables that are independent of ${ W^{\varphi}_{\infty}}$.
\end{lemma}
\begin{proof}
Fix arbitrary $s\ge 0$ and $x\in E$. Let $K_{t}(\theta_{1},\theta_{2})$ be the characteristic function of the random vector $\big(\e^{-\frac{\lambda_{1}}{2}(t+s)}\left(\langle f_{t+s},X_{t+s}\rangle-\pp_{\delta_{x}}\left[\langle f_{t+s},X_{t+s}\rangle |\mathcal{F}_{t}\right]\right),\e^{\frac{\lambda_{1}}{2}(t+s)}(W^{\varphi}_{t+s}-W^{\varphi}_{\infty})\big)^{\mathrm{T}}$. It suffices to prove that for every $\theta_{1},\theta_{2}\in\R$,
\begin{equation}\label{lem7.3.2}
\lim_{t\to+\infty}K_{t}(\theta_{1},\theta_{2})=\pp_{\delta_{x}}\left[\exp\{-\frac{1}{2}W^{\varphi}_{\infty}\left(\theta^{2}_{1}\rho^{2}_{s}(f)+\theta^{2}_{2}\sigma^{2}_{\varphi}\right)\}\right].
\end{equation}
We have
\begin{align*}
&K_{t}(\theta_{1},\theta_{2})\\
=&\pp_{\delta_{x}}\left[\exp\left\{i\theta_{1}\e^{-\frac{\lambda_{1}}{2}(t+s)}\left(\langle f_{t+s},X_{t+s}\rangle-\langle T_{s}f_{t+s},X_{t}\rangle\right)+i\theta_{2}\e^{\frac{\lambda_{1}}{2}(t+s)}({ W^{\varphi}_{t+s}}-{ W^{\varphi}_{\infty}})\right\}\right]\\
=&\pp_{\delta_{x}}\left[\exp\left\{-i\theta_{1}\e^{-\frac{\lambda_{1}}{2}(t+s)}\langle T_{s}f_{t+s},X_{t}\rangle+i\e^{-\frac{\lambda_{1}}{2}(t+s)}\langle\theta_{1}f_{t+s}+\theta_{2}\varphi,X_{t+s} \rangle\right\}\pp_{X_{t+s}}\left[\e^{-i\theta_{2}\e^{-\frac{\lambda_{1}}{2}(t+s)}{ W^{\varphi}_{\infty}}}\right]\right]\\
=&\pp_{\delta_{x}}\left[\exp\left\{-i\theta_{1}\e^{-\frac{\lambda_{1}}{2}(t+s)}\langle T_{s}f_{t+s},X_{t}\rangle+\langle i\theta_{1}\e^{-\frac{\lambda_{1}}{2}(t+s)}f_{t+s}-\frac{1}{2}\theta^{2}_{2}\e^{-\lambda_{1}(t+s)}\Theta+I_{\theta_{2}}(t+s,\cdot),X_{t+s}\rangle\right\}\right],
\end{align*}
where for $y\in E$,
$$I_{\theta_{2}}(t+s,y):=\int_{(0,+\infty)}\left(\exp\{-i\theta_{2}\e^{-\frac{\lambda_{1}}{2}(t+s)r}\}-1+i\theta_{2}\e^{-\frac{\lambda_{1}}{2}(t+s)r}+\frac{1}{2}\theta^{2}_{2}\e^{-\lambda_{1}(t+s)}r^{2}\right)\Upsilon(y,{ \ud r}).$$
The second equality follows from the Markov property, and the third equality is from \eqref{eq:charac}.
Then by the Markov property and \eqref{eq3}, we have
\begin{align}
&K_{t}(\theta_{1},\theta_{2})\nonumber\\
=&\pp_{\delta_{x}}\Big[\exp\big\{-i\theta_{1}\e^{-\frac{\lambda_{1}}{2}(t+s)}\langle T_{s}f_{t+s},X_{t}\rangle\big\}\nonumber\\
&\quad\pp_{X_{t}}\big[\exp\big\{\langle i\theta_{1}\e^{-\frac{\lambda_{1}}{2}(t+s)}f_{t+s}-\frac{1}{2}\theta^{2}_{2}\e^{-\lambda_{1}(t+s)}\Theta+I_{\theta_{2}}(t+s,\cdot),X_{s}\rangle\big\}\big]
\Big]\nonumber\\
=&\pp_{\delta_{x}}\big[\exp\big\{\e^{-\lambda_{1} t}\big\langle-T_{s}\big(\frac{1}{2}\theta^{2}_{2}\e^{-\lambda_{1} s}\Theta-\e^{\lambda_{1} t}I_{\theta_{2}}(t+s,\cdot)\big)\nonumber\\
&\quad+
\frac{1}{2}\mathrm{Var}_{\delta_{\cdot}}\big[\langle\frac{1}{2}\theta^{2}_{2}\e^{-\frac{\lambda_{1}}{2}t-\lambda_{1} s}\Theta-\e^{\frac{\lambda_{1}}{2}t}I_{\theta_{2}}(t+s,\cdot)-i\theta_{1}\e^{-\frac{\lambda_{1}}{2}s}f_{t+s},X_{s}\rangle\big]\nonumber\\
&\quad+\e^{\lambda_{1} t}\mathcal{L}_{s}\big[\e^{-\frac{\lambda_{1}}{2}t}\big(\frac{1}{2}\theta^{2}_{2}\e^{-\frac{\lambda_{1}}{2}t-\lambda_{1} s}\Theta-\e^{\frac{\lambda_{1}}{2}t}I_{\theta_{2}}(t+s,\cdot)-i\theta_{1}\e^{-\frac{\lambda_{1}}{2}s}f_{t+s}\big)\big],X_{t}\big\rangle\big\}\big]\nonumber\\
=:&\pp_{\delta_{x}}\big[\exp\{\e^{-\lambda_{1} t}\langle V(s,t,\theta_{1},\theta_{2}),X_{t}\rangle\}\big].\label{lem2case1.1}
\end{align}
We note that by \eqref{eq:inequality2},
$$\left|\e^{\lambda_{1} t}I_{\theta_{2}}(t+s,y)\right|\le\theta^{2}_{2}\e^{-\lambda_{1} s}\int_{(0,+\infty)}r^{2}\left(1\wedge \frac{|\theta_{2}|r\e^{-\frac{\lambda_{1}}{2}(t+s)}}{6}\right)\Upsilon(y,{ \ud r})\stackrel{b}{\to}0\mbox{ as }t\to+\infty.$$
Thus we have for every $s\ge 0$,
$$T_{s}\left(\frac{1}{2}\theta^{2}_{2}\e^{-\lambda_{1} s}\Theta-\e^{\lambda_{1} t}I_{\theta_{2}}(t+s,\cdot)\right)\stackrel{b}{\to}\frac{1}{2}\theta^{2}_{2}\e^{-\lambda_{1} s}T_{s}\Theta\mbox{ as }t\to+\infty,$$
and
\begin{equation}\label{eq4.4}
\frac{1}{2}\theta^{2}_{2}\e^{-\frac{\lambda_{1}}{2}t-\lambda_{1} s}\Theta-\e^{\frac{\lambda_{1}}{2}t}I_{\theta_{2}}(t+s,\cdot)-i\theta_{1}\e^{-\frac{\lambda_{1}}{2}s}f_{t+s}\stackrel{b}{\to}
-i\theta_{1}\e^{-\frac{\lambda_{1}}{2}s}f\mbox{ as }t\to+\infty.
\end{equation}
It then follows from \eqref{eq4.4}, Lemma \ref{lem5.1}(i) and (ii) that for $y\in E$,
\begin{eqnarray*}
&&\mathrm{Var}_{\delta_{y}}\left[\langle\frac{1}{2}\theta^{2}_{2}\e^{-\frac{\lambda_{1}}{2}t-\lambda_{1} s}\Theta-\e^{\frac{\lambda_{1}}{2}t}I_{\theta_{2}}(t+s,\cdot)-i\theta_{1}\e^{-\frac{\lambda_{1}}{2}s}f_{t+s},X_{s}\rangle\right]\\
&=&\int_{0}^{s}T_{s-r}\left(\vartheta\left[T_{r}\left(\frac{1}{2}\theta^{2}_{2}\e^{-\frac{\lambda_{1}}{2}t-\lambda_{1} s}\Theta-\e^{\frac{\lambda_{1}}{2}t}I_{\theta_{2}}(t+s,\cdot)-i\theta_{1}\e^{-\frac{\lambda_{1}}{2}s}f_{t+s}\right)\right]\right)(y){ \ud r}\\
&\stackrel{b}{\to}&\int_{0}^{s}T_{s-r}\left(\vartheta\left[T_{r}\left(-i\theta_{1}\e^{-\frac{\lambda_{1}}{2}s}f\right)\right]\right)(y){ \ud r}
=-\theta^{2}_{1}\e^{-\lambda_{1} s}\mathrm{Var}_{\delta_{y}}[\langle f,X_{s}\rangle],
\end{eqnarray*}
and
$$\e^{\lambda_{1} t}\mathcal{L}_{s}\left[\e^{-\frac{\lambda_{1}}{2}t}\left(\frac{1}{2}\theta^{2}_{2}\e^{-\frac{\lambda_{1}}{2}t-\lambda_{1} s}\Theta-\e^{\frac{\lambda_{1}}{2}t}I_{\theta_{2}}(t+s,\cdot)-i\theta_{1}\e^{-\frac{\lambda_{1}}{2}s}f_{t+s}\right)\right]\stackrel{b}{\to}0$$
as $t\to+\infty$. Hence by Proposition \ref{lem2.2}(i),
$$\e^{-\lambda_{1} t}\langle V(s,t,\theta_{1},\theta_{2}),X_{t}\rangle\to \langle-\frac{1}{2}\theta^{2}_{2}\e^{-\lambda_{1} s}T_{s}\Theta-\frac{1}{2}\theta^{2}_{1}\e^{-\lambda_{1} s}\mathrm{Var}_{\delta_{\cdot}}[\langle f,X_{s}\rangle],\widetilde{\varphi}\rangle { W^{\varphi}_{\infty}}\mbox{ in }L^{2}(\pp_{\delta_{x}})$$
as $t\to+\infty$, and consequently by \eqref{lem2case1.1}
\begin{equation}\label{lem2case1.2}
\lim_{t\to+\infty}K_{t}(\theta_{1},\theta_{2})=\pp_{\delta_{x}}\left[\exp\{\langle-\frac{1}{2}\theta^{2}_{2}\e^{-\lambda_{1} s}T_{s}\Theta-\frac{1}{2}\theta^{2}_{1}\e^{-\lambda_{1} s}\mathrm{Var}_{\delta_{\cdot}}[\langle f,X_{s}\rangle],\widetilde{\varphi}\rangle { W^{\varphi}_{\infty}}\}\right].
\end{equation}
By Fubini's theorem,
$T_{s}\Theta=\int_{0}^{+\infty}\e^{-2\lambda_{1} t}T_{t+s}(\vartheta[\varphi]){ \ud t}$. Thus we have
$$\langle \e^{-\lambda_{1} s}T_{s}\Theta,\widetilde{\varphi}\rangle=\e^{-\lambda_{1} s}\int_{0}^{+\infty}\e^{-2\lambda_{1} t+\lambda_{1}(t+s)}\langle\vartheta[\varphi],\widetilde{\varphi}\rangle { \ud t}={ \sigma^{2}_{\varphi}},$$
and
$$\langle \e^{-\lambda_{1} s}\mathrm{Var}_{\delta_{\cdot}}[\langle f,X_{s}\rangle],\widetilde{\varphi}\rangle=\e^{-\lambda_{1} s}\int_{0}^{s}\langle T_{s-r}(\vartheta[T_{r}f]),\widetilde{\varphi}\rangle { \ud r}=\int_{0}^{s}\e^{-\lambda_{1} r}\langle \vartheta[T_{r}f],\widetilde{\varphi}\rangle { \ud r}=\rho^{2}_{s}(f).$$
Inserting these to \eqref{lem2case1.2} yields \eqref{lem7.3.2}.
\end{proof}

\begin{proposition}\label{lem2case1}
Suppose $\epsilon(g)>\lambda_{1}/2$. Then for every $x\in E$, under $\pp_{\delta_{x}}$,
$$\big(\e^{-\frac{\lambda_{1}}{2}t}\langle g,X_{t}\rangle,\e^{\frac{\lambda_{1}}{2}t}({ W^{\varphi}_{t}}-{ W^{\varphi}_{\infty}})\big)^{\mathrm{T}}
\stackrel{d}{\to}\big({\rho_{g}}\sqrt{{ W^{\varphi}_{\infty}}}N_{1},
\sigma_{\varphi}\sqrt{{ W^{\varphi}_{\infty}}} N_{2}\big)^{\mathrm{T}}\mbox{ as }t\to+\infty,$$
where $N_{1},N_{2}$ are two independent standard normal random variables, and $N_{1},N_{2}$ are independent of ${ W^{\varphi}_{\infty}}$.
\end{proposition}

\begin{proof}
For $t,s\ge 0$, define
\begin{eqnarray}
U(s+t)&:=&\big(\e^{-\frac{\lambda_{1}}{2}(t+s)}\langle g,X_{t+s}\rangle,\e^{\frac{\lambda_{1}}{2}(t+s)}({ W^{\varphi}_{t+s}}-{ W^{\varphi}_{\infty}})\big)^{\mathrm{T}}\nonumber\\
&=&\big(\e^{-\frac{\lambda_{1}}{2}(t+s)}\big(\langle g,X_{t+s}\rangle-\langle T_{s}g,X_{t}\rangle\big),\e^{\frac{\lambda_{1}}{2}(t+s)}({ W^{\varphi}_{t+s}}-{ W^{\varphi}_{\infty}})\big)^{\mathrm{T}}\nonumber\\
&&\quad+\big(\e^{-\frac{\lambda_{1}}{2}(t+s)}\langle T_{s}g,X_{t}\rangle,0\big)^{\mathrm{T}}\nonumber\\
&=:&U_{1}(s,t)+U_{2}(s,t).\label{lem4.2.1}
\end{eqnarray}
We observe that the double limit, first as $t\to+\infty$ and then $s\to+\infty$, of the right hand side of \eqref{lem4.2.1}, if exists, is equal to the limit of $U(t)$ as $t\to+\infty$.
By Lemma \ref{lem4.1}, we have for every $s\ge 0$, $U_{1}(s,t)\stackrel{d}{\to} (\sqrt{W^{\varphi}_{\infty}}\rho_{s}(g)N_{1},\sqrt{W^{\varphi}_{\infty}}\sigma_{\varphi}N_{2})^{\mathrm{T}}$ as $t\to+\infty$. We note that
$$\rho^{2}_{s}(g)=\int_{0}^{s}\e^{-\lambda_{1} r}\langle\vartheta[T_{r}g],\widetilde{\varphi}\rangle { \ud r}\uparrow \int_{0}^{+\infty}\e^{-\lambda_{1} r}\langle\vartheta[T_{r}g],\widetilde{\varphi}\rangle { \ud r}=\rho^{2}_{g}\mbox{ as }s\to+\infty.$$
Hence $U_{1}(s,t)$ converges in distribution to $(\sqrt{W^{\varphi}_{\infty}}\rho_{g}N_{1},\sqrt{W^{\varphi}_{\infty}}\sigma_{\varphi}N_{2})^{\mathrm{T}}$, as first $t\to+\infty$ and then $s\to+\infty$.

Next we deal with $U_{2}(s,t)$. Let $h_{s}(y):=\e^{-\frac{\lambda_{1}}{2}s}T_{s}g(y)$ for $s\ge 0$ and $y\in E$. Then we have $\langle h_{s},\widetilde{\varphi}\rangle=\e^{-\frac{\lambda_{1}}{2}s}\langle T_{s}g,\widetilde{\varphi}\rangle=\e^{\frac{\lambda_{1}}{2}s}\langle g,\widetilde{\varphi}\rangle=0$, and
$\e^{-\lambda_{1} t}T_{t}h_{s}=\e^{\frac{\lambda_{1}}{2}s}\e^{-\lambda_{1}(t+s)}T_{t+s}g$. Immediately, we have $\epsilon(h_{s})=\epsilon(g)>\lambda_{1}/2$. Hence by Lemma \ref{lem1case1},  $\e^{-\frac{\lambda_{1}}{2}t}T_{t}h_{s}(x)\to 0$ and $\e^{-\lambda_{1} t}\mathrm{Var}_{\delta_{x}}[\langle h_{s},X_{t}\rangle]\to { \rho^{2}_{h_{s}}}\varphi(x)$ as $t\to+\infty$, where
$${ \rho^{2}_{h_{s}}}=\int_{0}^{+\infty}\e^{-\lambda_{1} r}\langle \vartheta[T_{r}h_{s}],\widetilde{\varphi}\rangle { \ud r}=\int_{s}^{+\infty}\e^{-\lambda_{1} r}\langle \vartheta[T_{r}g],\widetilde{\varphi}\rangle { \ud r}.$$
Since
$$\pp_{\delta_{x}}\left[\langle\e^{-\frac{\lambda_{1}}{2}(t+s)}T_{s}g,X_{t}\rangle^{2}\right]
=\e^{-\lambda_{1} t}\pp_{\delta_{x}}\left[\langle h_{s},X_{t}\rangle^{2}\right]=\e^{-\lambda_{1} t}\mathrm{Var}_{\delta_{x}}[\langle h_{s},X_{t}\rangle]+\left(\e^{-\frac{\lambda_{1}}{2}t}T_{t}h_{s}(x)\right)^{2},$$
we get
$$\lim_{s\to+\infty}\lim_{t\to+\infty}\pp_{\delta_{x}}\left[\langle\e^{-\frac{\lambda_{1}}{2}(t+s)}T_{s}g,X_{t}\rangle^{2}\right]
=\lim_{s\to+\infty}{ \rho^{2}_{h_{s}}}\varphi(x)=0.$$
This implies that $\lim_{s\to+\infty}\lim_{t\to+\infty}U_{2}(s,t)=0$ in $L^{2}(\pp_{\delta_{x}})$.
Hence by Slutsky theorem, $U(s+t)=U_{1}(s,t)+U_{2}(s,t)$ converges in distribution to $({\rho_{g}}\sqrt{{ W^{\varphi}_{\infty}}}N_{1},\sigma_{\varphi}\sqrt{{ W^{\varphi}_{\infty}}} N_{2})^{\mathrm{T}}$, as first $t\to+\infty$ and then $s\to+\infty$.
Thus we have proved this lemma.
\end{proof}

{\subsection{Critical branching case $\epsilon(g)=\lambda_{1}/2$}}
\begin{proposition}\label{lem2case2}
Suppose $\pp_{\delta_{x}}\big(X_{t}=0\mbox{ for some }t\in (0,+\infty)\big)>0$ and \eqref{lemnew.1} holds for $\epsilon=\lambda_{1}/2$.
Let $N$ be a standard normal random variable independent of $W^{\varphi}_{\infty}$. Then the following are true.
\begin{description}
\item{(i)} Under $\pp_{\delta_{x}}$,
$$\frac{1}{\sqrt{t}}\e^{-\frac{\lambda_{1}}{2}t}\langle g^{*},X_{t}\rangle\stackrel{d}{\to}
{\varrho_{g^{*}}}\sqrt{{ W^{\varphi}_{\infty}}}N\mbox{ as }t\to+\infty.$$
Moreover, if \eqref{lemnew.1} holds for $r=0$, then under $\pp_{\delta_{x}}$,
$$\frac{1}{\sqrt{t}}\e^{-\frac{\lambda_{1}}{2}t}\langle g,X_{t}\rangle\stackrel{d}{\to}{\varrho_{g^{*}}}\sqrt{{ W^{\varphi}_{\infty}}}N\mbox{ as }t\to+\infty.$$
\item{(ii)} If \eqref{lemnew.1} holds for $r>0$, and for all $d,\delta>0$,
\begin{equation}\label{lem2case2.0}
\lim_{t\to+\infty}\pp_{\delta_{x}}\left[\left|\frac{1}{t^{\frac{1}{2}+r}}\e^{-\frac{\lambda_{1}}{2}t}\langle g,X_{t}\rangle\right|^{2},\ \left|\frac{1}{t^{\frac{1}{2}+r}}\e^{-\frac{\lambda_{1}}{2}t}\langle g,X_{t}\rangle\right|>d\e^{\delta t}\right]=0,
\end{equation}
then under $\pp_{\delta_{x}}$,
$$\frac{1}{t^{\frac{1}{2}+r}}\e^{-\frac{\lambda_{1}}{2}t}\langle g,X_{t}\rangle\stackrel{d}{\to}\frac{\varrho_{g^{*}}}{\sqrt{1+2r}}\sqrt{{ W^{\varphi}_{\infty}}}N\mbox{ as }t\to+\infty.$$
\end{description}
\end{proposition}

We note that by Chebyshev's inequality and Lemma \ref{lem:newcase2}(i), for any $d,\delta>0$,
\begin{align*}
\pp_{\delta_{x}}\left(\left|\frac{1}{t^{1/2+r}}\e^{-\frac{\lambda_{1}}{2}t}\langle g,X_{t}\rangle\right|>d\e^{\delta t}\right)&\le d^{-2}\e^{-2\delta t}\frac{1}{t^{1+2r}}\e^{-\lambda_{1} t}\pp_{\delta_{x}}\left[\langle g,X_{t}\rangle^{2}\right]\\
&=d^{-2}\e^{-2\delta t}\frac{1}{t^{1+2r}}\e^{-\lambda_{1} t}\left(\mathrm{Var}_{\delta_{x}}[\langle g,X_{t}\rangle]+(T_{t}g(x))^{2}\right)\to 0
\end{align*}
as $t\to+\infty$. Hence a sufficient condition for \eqref{lem2case2.0} to hold is that the family $\{\frac{1}{t^{1/2+r}}\e^{-\frac{\lambda_{1}}{2}t}\langle g,X_{t}\rangle:\ t>0\}$ is uniformly integrable.
By Lemma \ref{lemA.4}, this holds if $\sup_{y\in E}\int_{\mfo }{\langle 1,\nu\rangle}^{4}H(y,{ \ud \nu})<+\infty$.
Thus we obtain the following result.

\begin{corollary}\label{cor4.7}
Suppose \[\sup_{y\in E}\int_{\mfo}{\langle 1,\nu\rangle}^{4}H(y,{ \ud \nu})<+\infty,\]
and that $\pp_{\delta_{x}}\left(X_{t}=0\mbox{ for some }t\in (0,+\infty)\right)>0$. If \eqref{lemnew.1} holds for $\epsilon=\lambda_{1}/2$, then
under $\pp_{\delta_{x}}$,
$$\frac{1}{t^{\frac{1}{2}+r}}\e^{-\frac{\lambda_{1}}{2}t}\langle g,X_{t}\rangle\stackrel{d}{\to}\frac{\varrho_{g^{*}}}{\sqrt{1+2r}}\sqrt{{ W^{\varphi}_{\infty}}}N\mbox{ as }t\to+\infty,$$
where $N$ is a standard normal random variable independent of $W^{\varphi}_{\infty}$.
\end{corollary}

The following two lemmas are needed for the proofs of Proposition \ref{lem2case2}(i) and (ii), respectively.

\begin{lemma}\label{lem4case2}
Suppose the assumptions of Proposition \ref{lem2case2} hold.
Then for every $d,\delta>0$ and $x\in E$,
\begin{equation}\label{lem4case2.0}
\lim_{t\to+\infty}\pp_{\delta_{x}}\left[\left|\frac{1}{\sqrt{t}}\e^{-\frac{\lambda_{1}}{2}t}\langle g^{*},X_{t}\rangle\right|^{2},\ \left|\frac{1}{\sqrt{t}}\e^{-\frac{\lambda_{1}}{2}t}\langle g^{*},X_{t}\rangle\right|>d\e^{\delta t}\right]=0.
\end{equation}
\end{lemma}
\begin{proof}
{Recall from Remark \ref{rm1} that $T_{t}g^{*}=\e^{\frac{\lambda_{1}}{2}t}g^{*}$ for $t\ge 0$.}
For $x\in E$ and $t>0$, let
$$S_{t}g^{*}(x):=\frac{1}{\sqrt{t}}\e^{-\frac{\lambda_{1}}{2}t}\langle g^{*},X_{t}\rangle-\frac{1}{\sqrt{t}}g^{*}(x)=\frac{1}{\sqrt{t}}\e^{-\frac{\lambda_{1}}{2}t}\left(\langle g^{*},X_{t}\rangle-T_{t}g^{*}(x)\right).$$
By Lemma \ref{lem:newcase2}(i) we have
\begin{equation}\label{lem4case2.1}
\pp_{\delta_{x}}\left[|S_{t}g^{*}(x)|^{2}\right]=\frac{1}{t}\e^{-\lambda_{1} t}\mathrm{Var}_{\delta_{x}}[\langle g^{*},X_{t}\rangle]\stackrel{b}{\to}\varrho^{2}_{g^{*}}\varphi(x)\mbox{ as }t\to+\infty.
\end{equation}
Let
$$R(t,g^{*}):=\e^{-\frac{\lambda_{1}}{2}(t+1)}\left(\langle g^{*},X_{t+1}\rangle-\pp_{\delta_{x}}\left[\langle g^{*},X_{t+1}\rangle |\mathcal{F}_{t}\right]\right)=\e^{-\frac{\lambda_{1}}{2}(t+1)}\left(\langle g^{*},X_{t+1}\rangle-\langle T_{1}g^{*},X_{t}\rangle\right).$$
Then by the Markov property,
\begin{eqnarray*}
\pp_{\delta_{x}}\left[R(t,g^{*})^{2}\right]&=&\pp_{\delta_{x}}\left[\mathrm{Var}\left[\e^{-\frac{\lambda_{1}}{2}(t+1)}\langle g^{*},X_{t+1}\rangle |\mathcal{F}_{t}\right]\right]\\
&=&\e^{-\lambda_{1}(t+1)}\pp_{\delta_{x}}\left[\mathrm{Var}_{X_{t}}[\langle g^{*},X_{1}\rangle]\right]\\
&=&\e^{-\lambda_{1}(t+1)}T_{t}(\mathrm{Var}_{\delta_{\cdot}}[\langle g^{*},X_{1}\rangle])(x).
\end{eqnarray*}
Thus by Assumption \ref{AS1} we have
\begin{equation}\label{lem4case2.2}
\pp_{\delta_{x}}\left[R(t,g^{*})^{2}\right]\stackrel{b}{\to}\varphi(x)\langle \e^{-\lambda_{1}}\mathrm{Var}_{\delta_{\cdot}}[\langle g^{*},X_{1}\rangle],\widetilde{\varphi}\rangle\mbox{ as }t\to+\infty.
\end{equation}
Moreover, we have by Lemma \ref{lem4.1} that $R(t,g^{*})\stackrel{d}{\to}\sqrt{{ W^{\varphi}_{\infty}}}\rho_{1}(g^{*})N$ as $t\to+\infty$ where $N$ is a standard normal random variable independent of ${ W^{\varphi}_{\infty}}$. Based on this, \eqref{lem4case2.1} and \eqref{lem4case2.2}, we can apply similar argument as in the proof of \cite[Lemma 3.2]{RSZ14a} (with (2.15) and (3.44) there replaced by \eqref{lem4case2.1} and \eqref{lem4case2.2}, respectively) to show that for any $d,\delta>0$,
$$\lim_{t\to+\infty}\pp_{\delta_{x}}\left(|S_{t}g^{*}(x)|^{2},\ |S_{t}g^{*}(x)|>d\e^{\delta t}\right)=0.$$
Hence \eqref{lem4case2.0} follows immediately from the above convergence and the fact that
\[\frac{1}{\sqrt{t}}\e^{-\frac{\lambda_{1}}{2}t}\langle g^{*},X_{t}\rangle=S_{t}g^{*}(x)+\frac{1}{\sqrt{t}}g^{*}(x)=S_{t}g^{*}(x)+O(t^{-1/2}).\]
\end{proof}

\begin{lemma}\label{lem3case2}
Suppose the assumptions of Proposition \ref{lem2case2} hold. Then condition \eqref{lem2case2.0} implies that for any $d,\delta>0$,
\begin{equation}\label{lem3case2.0}
\frac{1}{t^{1+2r}}\e^{-\lambda_{1} t}\int_{\frac{1}{t^{1/2+r}}\e^{-\frac{\lambda_{1}}{2}t}|{\langle g,\nu\rangle}|>d\e^{\delta t}}{\langle g,\nu\rangle}^{2}L_{t}(x,{ \ud \nu})\stackrel{b}{\to}0\mbox{ as }t\to+\infty.
\end{equation}
\end{lemma}

\begin{proof}
The boundedness of the left hand side of \eqref{lem3case2.0} for large $t$ is obvious, since for $x\in E$,
$$\frac{1}{t^{1+2r}}\e^{-\lambda_{1} t}\int_{\mfo}{\langle g,\nu\rangle}^{2}L_{t}(x,{ \ud \nu})=\frac{1}{t^{1+2r}}\e^{-\lambda_{1} t}\mathrm{Var}_{\delta_{x}}[\langle g,X_{t}\rangle],$$
and the right hand side is bounded from above for large $t$ by Lemma \ref{lem:newcase2}(i). We only need to show the pointwise convergence in \eqref{lem3case2.0}.
For $t>0$, let $g_{t}:=g/t^{r}$ and $Y^{(t,x)}$ be a random variable with the same law as $(\frac{1}{\sqrt{t}}\e^{-\frac{\lambda_{1}}{2}t}\langle g_{t},X_{t}\rangle,\pp_{\delta_{x}})$.
Then $Y^{(t,x)}$ is an infinitely divisible random variable, and by \eqref{eq2} its characteristic function can be represented by
$$\mathrm{E}\left[\e^{i\theta Y^{(t,x)}}\right]=\exp\{i\delta^{(t,x)}\theta-\int_{\R}\left(1-\e^{i\theta y}\right)\pi^{(t,x)}({ \ud y})\}\quad\forall \theta\in\R,$$
where $\delta^{(t,x)}:=-\frac{1}{\sqrt{t}}\e^{-\frac{\lambda_{1}}{2}t}\Lambda(x,g_{t})$, and $\pi^{(t,x)}({ \ud y})$ is the L\'{e}vy measure on $\R$ satisfying that
$$\pi^{(t,x)}(A)=\int_{\mfo }1_{\{\frac{1}{\sqrt{t}}\e^{-\frac{\lambda_{1}}{2}t}\nu(g_{t})\in A\}}L_{t}(x,{ \ud \nu})\quad\forall A\in\mathcal{B}(\R).$$
Thus it suffices to prove that for any $d,\delta>0$,
\begin{equation}\label{lem3case2.1}
\lim_{t\to+\infty}\int_{|y|>d\e^{\delta t}}y^{2}\pi^{(t,x)}({ \ud y})=0.
\end{equation}
We note that for $\theta,t\ge 0$,
\begin{eqnarray*}
\pp_{\delta_{x}}\left[\e^{-\theta W^{\varphi}_{t}}\right]&=&\exp\{-\theta\e^{-\lambda_{1} t}\Lambda_{t}(x,\varphi)-\int_{\mfo }\left(1-\e^{-\theta\e^{-\lambda_{1} t}{\langle \varphi,\nu\rangle}}\right)L_{t}(x,{ \ud \nu})\}\\
&\le&\exp\{-\int_{\mfo }\left(1-\e^{-\theta\e^{-\lambda_{1} t}{\langle \varphi,\nu\rangle}}\right)L_{t}(x,{ \ud \nu})\}.
\end{eqnarray*}
Letting $\theta\to +\infty$, we have $\exp\{-\int_{\mfo }L_{t}(x,{ \ud \nu})\}\ge \pp_{\delta_{x}}\left(W_{t}=0\right)=\pp_{\delta_{x}}\left(X_{t}=0\right)$. Consequently,
$$\limsup_{t\to+\infty}\int_{\mfo }L_{t}(x,{ \ud \nu})\le\limsup_{t\to+\infty}-\log\pp_{\delta_{x}}\left(X_{t}=0\right)=-\log\pp_{\delta_{x}}\left(X_{t}=0\mbox{ for some }t\in (0,+\infty)\right).$$
Thus, given that $\pp_{\delta_{x}}\left(X_{t}=0\mbox{ for some }t\in (0,+\infty)\right)>0$, there is some $t_{0}>0$ such that $\sup_{t\ge t_{0}}\pi^{(t,x)}(\R)=\sup_{t\ge t_{0}}\int_{\mfo }L_{t}(x,{ \ud \nu})<+\infty$. Since the L\'{e}vy measure $\pi^{(t,x)}$ is finite for $t\ge t_{0}$, we can identify $Y^{(t,x)}$ in the form
$$Y^{(t,x)}=\delta^{(t,x)}+\sum_{i=1}^{N}\xi_{i}$$
where $N$ is a Poisson random variable with parameter $\pi^{(t,x)}(\R)$, and $\{\xi_{i}:i\ge 1\}$ is an independent sequence of i.i.d. $\R$-valued random variables with law $\pi^{(t,x)}({ \ud y})/\pi^{(t,x)}(\R)$. Condition \eqref{lem2case2.0} yields that
\begin{equation}\label{lem3case2.2}
\lim_{t\to+\infty}\mathrm{P}\left(|Y^{(t,x)}|^{2},\ |Y^{(t,x)}|>d\e^{\delta t}\right)=0
\end{equation}
for all $d,\delta>0$. Note that
\begin{align*}
\mathrm{P}\left(|Y^{(t,x)}|^{2},\ |Y^{(t,x)}|>d\e^{\delta t}\right)=&\mathrm{P}\left(|\sum_{i=1}^{N}\xi_{i}+\delta^{(t,x)}|^{2},\ |\sum_{i=1}^{N}\xi_{i}+\delta^{(t,x)}|>d\e^{\delta t}\right)\\
\ge&\mathrm{P}\left(|\xi_{1}+\delta^{(t,x)}|^{2},\ N=1,\ |\xi_{1}+\delta^{(t,x)}|>d\e^{\delta t}\right)\\
=&\pp(N=1)\mathrm{P}\left(|\xi_{1}+\delta^{(t,x)}|^{2},\ |\xi_{1}+\delta^{(t,x)}|>d\e^{\delta t}\right)\\
=&\e^{-\pi^{(t,x)}(\R)}\int_{|y+\delta^{(t,x)}|>d\e^{\delta t}}|y+\delta^{(t,x)}|^{2}\pi^{(t,x)}({ \ud y}).
\end{align*}
By \eqref{lem3case2.2} and the fact that $\sup_{t\ge t_{0}}\pi^{(t,x)}(\R)<+\infty$, we get
\begin{equation}\label{lem3case2.3}
\lim_{t\to+\infty}\int_{|y+\delta^{(t,x)}|>d\e^{\delta t}}|y+\delta^{(t,x)}|^{2}\pi^{(t,x)}({ \ud y})=0.
\end{equation}
We note that $\delta^{(t,x)}=\mathrm{E}[Y^{(t,x)}]-\int_{\R}y\pi^{(t,x)}({ \ud y})$, where
$\mathrm{E}[Y^{(t,x)}]=\frac{1}{t^{1/2+r}}\e^{-\frac{\lambda_{1}}{2}t}T_{t}g(x)\to 0$ as $t\to+\infty$, and
$$|\int_{\R}y\pi^{(t,x)}({ \ud y})|\le \frac{1}{2}\int_{\R}(1+y^{2})\pi^{(t,x)}({ \ud y})=\frac{1}{2}(\pi^{(t,x)}(\R)+\frac{1}{t^{1+2r}}\e^{-\lambda_{1} t}\mathrm{Var}_{\delta_{x}}[\langle g,X_{t}\rangle]).$$
Lemma \ref{lem:newcase2}(i) implies that $\delta^{(t,x)}=O(1)$ as $t\to+\infty$. Therefore, \eqref{lem3case2.1} follows directly from \eqref{lem3case2.3}.
\end{proof}

\noindent\textbf{Proof of Proposition \ref{lem2case2}:}
We shall prove first (ii) and then (i).
For $r\not=0$, $n\in\mathbb{N}$ and $t>0$, we may write
\begin{eqnarray}
\frac{1}{(nt)^{1/2+r}}\e^{-\frac{\lambda_{1}}{2}nt}\langle g,X_{nt}\rangle
&=&\frac{1}{(nt)^{1/2+r}}\e^{-\frac{\lambda_{1}}{2}nt}\left(\langle g,X_{nt}\rangle
-\langle T_{(n-1)t}g,X_{t}\rangle\right)\nonumber\\
&&\quad+\frac{1}{(nt)^{1/2+r}}\e^{-\frac{\lambda_{1}}{2}nt}\langle T_{(n-1)t}g,X_{t}\rangle\nonumber\\
&=:&U_{1}(n,t)+U_{2}(n,t).\label{eq4.18}
\end{eqnarray}
We observe that the double limit, first as $t\to+\infty$ and then $n\to+\infty$, of the right hand side of \eqref{eq4.18}, if exists, is equal to the limit of $\frac{1}{t^{1/2+r}}\e^{-\frac{\lambda_{1}}{2}t}\langle g,X_{t}\rangle$ as $t\to+\infty$.
By the Markov property and \eqref{eq3}, we have for $\theta\in\R$,
\begin{align}
&\pp_{\delta_{x}}\left[\exp\left\{i\theta U_{1}(n,t)\right\}\right]\nonumber\\
=&\pp_{\delta_{x}}\left[\exp\left\{-i\theta\frac{1}{(nt)^{1/2+r}}\e^{-\frac{\lambda_{1}}{2}nt}\langle T_{(n-1)t}g,X_{t}\rangle\right\}\pp_{X_{t}}\left[\exp\left\{i\theta\frac{1}{(nt)^{1/2+r}}\e^{-\frac{\lambda_{1}}{2}nt}\langle g,X_{(n-1)t}\rangle\right\}\right]\right]\nonumber\\
=&\pp_{\delta_{x}}\big[\exp\big\{\e^{-\lambda_{1} t}\langle -\frac{1}{2}\theta^{2}\left(1-\frac{1}{n}\right)^{1+2r}\frac{1}{((n-1)t)^{1+2r}}\e^{-\lambda_{1}(n-1)t}
\mathrm{Var}_{\delta_{\cdot}}[\langle g,X_{(n-1)t}\rangle]\nonumber\\
&\quad\quad-\e^{\lambda_{1} t}\mathcal{L}_{(n-1)t}[-i\theta\frac{1}{(nt)^{1/2+r}}\e^{-\frac{\lambda_{1}}{2}nt}g],X_{t}\rangle\big\}\big]\nonumber\\
=:&\pp_{\delta_{x}}\left[\exp\{\e^{-\lambda_{1} t}\langle V_{1}(\theta,n,t),X_{t}\rangle\}\right].\label{eq4.19}
\end{align}
Here for $y\in E$,
\begin{align*}
&\mathcal{L}_{(n-1)t}[-i\theta\frac{1}{(nt)^{1/2+r}}\e^{-\frac{\lambda_{1}}{2}nt}g](y)\\
=&\int_{\mfo }\big(1-\exp\{i\theta{\langle\frac{1}{(nt)^{1/2+r}}\e^{-\frac{\lambda_{1}}{2}nt}g,\nu\rangle}\}+i\theta{\langle\frac{1}{(nt)^{1/2+r}}\e^{-\frac{\lambda_{1}}{2}nt}g,\nu\rangle}\\
&-\frac{1}{2}\theta^{2}{\langle\frac{1}{(nt)^{1/2+r}}\e^{-\frac{\lambda_{1}}{2}nt}g,\nu\rangle}^{2}\big)L_{(n-1)t}(y,{ \ud \nu}).
\end{align*}
It then follows by \eqref{eq:inequality2} that for any $\epsilon\in (0,1)$,
\begin{align*}
&\left|\mathcal{L}_{(n-1)t}[-i\theta\frac{1}{(nt)^{1/2+r}}\e^{-\frac{\lambda_{1}}{2}nt}g](y)\right|\\
\le&\frac{1}{(nt)^{1+2r}}\e^{-\lambda_{1} nt}\theta^{2}\int_{\mfo }{\langle g,\nu\rangle}^{2}\left(1\wedge\frac{|\theta|\frac{1}{(nt)^{1/2+r}}\e^{-\frac{\lambda_{1}}{2}nt}|{\langle g,\nu\rangle}|}{6}\right)L_{(n-1)t}(y,{ \ud \nu})\\
\le&\frac{1}{(nt)^{1+2r}}\e^{-\lambda_{1} nt}\theta^{2}\big[\frac{\epsilon}{6}\int_{\frac{1}{(nt)^{1/2+r}}\e^{-\frac{\lambda_{1}}{2}nt}|\theta||{\langle g,\nu\rangle}|\le\epsilon}{\langle g,\nu\rangle}^{2}L_{(n-1)t}(y,{ \ud \nu})\\
&\quad+\int_{\frac{1}{(nt)^{1/2+r}}\e^{-\frac{\lambda_{1}}{2}nt}|\theta||{\langle g,\nu\rangle}|>\epsilon}{\langle g,\nu\rangle}^{2}L_{(n-1)t}(y,{ \ud \nu})\big]\\
\le&\frac{1}{(nt)^{1+2r}}\e^{-\lambda_{1} nt}\theta^{2}\big[\frac{\epsilon}{6}\int_{\mfo }{\langle g,\nu\rangle}^{2}L_{(n-1)t}(y,{ \ud \nu})+\int_{\frac{1}{(nt)^{1/2+r}}\e^{-\frac{\lambda_{1}}{2}nt}|\theta||{\langle g,\nu\rangle}|>\epsilon}{\langle g,\nu\rangle}^{2}L_{(n-1)t}(y,{ \ud \nu})\big]\\
=&\frac{1}{(nt)^{1+2r}}\e^{-\lambda_{1} nt}\theta^{2}\big[\frac{\epsilon}{6}\mathrm{Var}_{\delta_{y}}[\langle g,X_{(n-1)t}\rangle]+\int_{\frac{1}{(nt)^{1/2+r}}\e^{-\frac{\lambda_{1}}{2}nt}|\theta||{\langle g,\nu\rangle}|>\epsilon}{\langle g,\nu\rangle}^{2}L_{(n-1)t}(y,{ \ud \nu})\big].
\end{align*}
Hence we have
\begin{eqnarray*}
&&\left|\e^{\lambda_{1} t}\mathcal{L}_{(n-1)t}[-i\theta\frac{1}{(nt)^{1/2+r}}\e^{-\frac{\lambda_{1}}{2}nt}g](y)\right|\\
&\le&\theta^{2}\left(1-\frac{1}{n}\right)^{1+2r}\frac{1}{((n-1)t)^{1+2r}}\e^{-\lambda_{1} (n-1)t}\big[\frac{\epsilon}{6}\mathrm{Var}_{\delta_{y}}[\langle g,X_{(n-1)t}\rangle]\\
&&\quad+\int_{\frac{1}{((n-1)t)^{1/2+r}}\e^{-\frac{\lambda_{1}}{2}(n-1)t}|\theta||{\langle g,\nu\rangle}|>\epsilon\left(\frac{n}{n-1}\right)^{1/2+r}\e^{\frac{\lambda_{1}}{2}t}}{\langle g,\nu\rangle}^{2}L_{(n-1)t}(y,{ \ud \nu})\big].
\end{eqnarray*}
Letting $t\to+\infty$ and then $\epsilon\to 0$, we get by Lemma \ref{lem:newcase2}(i) and \eqref{lem3case2.0} that
$$\e^{\lambda_{1} t}\mathcal{L}_{(n-1)t}[-i\theta\frac{1}{(nt)^{1/2+r}}\e^{-\frac{\lambda_{1}}{2}nt}g]\stackrel{b}{\to}0\mbox{ as }t\to+\infty.$$
It then follows by Lemma \ref{lem:newcase2}(i) and Proposition \ref{lem2.2}(i) that
\begin{align*}
&\lim_{t\to+\infty}\e^{-\lambda_{1} t}\langle V_{1}(\theta,n,t),X_{t}\rangle\\
=&\lim_{t\to+\infty}\langle -\frac{1}{2}\theta^{2}\left(1-\frac{1}{n}\right)^{1+2r}\frac{1}{((n-1)t)^{1+2r}}\e^{-\lambda_{1}(n-1)t}
\mathrm{Var}_{\delta_{\cdot}}[\langle g,X_{(n-1)t}\rangle],X_{t}\rangle\\
=&-\frac{1}{2}\theta^{2}\left(1-\frac{1}{n}\right)^{1+2r}\frac{\varrho^{2}_{g^{*}}}{1+2r}W^{\varphi}_{\infty}\quad\mbox{ in }L^{2}(\pp_{\delta_{x}}).
\end{align*}
Thus we get by \eqref{eq4.19} that
\begin{eqnarray*}
\lim_{n\to+\infty}\lim_{t\to+\infty}\pp_{\delta_{x}}\left[\exp\{i\theta U_{1}(n,t)\}\right]
&=&\lim_{n\to+\infty}\pp_{\delta_{x}}\left[\exp\left\{-\frac{1}{2}\theta^{2}\left(1-\frac{1}{n}\right)^{1+2r}\frac{\varrho^{2}_{g^{*}}}{1+2r}W^{\varphi}_{\infty}\right\}\right]\\
&=&\pp_{\delta_{x}}\left[\exp\left\{-\frac{1}{2}\theta^{2}\frac{\varrho^{2}_{g^{*}}}{1+2r}W^{\varphi}_{\infty}\right\}\right].
\end{eqnarray*}
This implies that $U_{1}(n,t)$ converges in distribution to $\frac{\varrho_{g^{*}}}{\sqrt{1+2r}}\sqrt{{ W^{\varphi}_{\infty}}}N$  as $t\to+\infty$ and then $n\to+\infty$.

Next we deal with $U_{2}(n,t)$. We note that
\begin{align*}
\mathrm{Var}_{\delta_{x}}[\langle T_{(n-1)t}g,X_{t}\rangle]
=&\int_{0}^{t}T_{t-s}\left(\vartheta[T_{(n-1)t+s}g]\right)(x){ \ud s}\\
=&\int_{(n-1)t}^{nt}T_{nt-u}\left(T_{u}g\right)(x){ \ud u}\\
=&\int_{0}^{nt}T_{nt-u}\left(T_{u}g\right)(x){ \ud u}-\int_{0}^{(n-1)t}T_{t}\circ T_{(n-1)t-u}\left(T_{u}g\right)(x){ \ud u}\\
=&\mathrm{Var}_{\delta_{x}}[\langle g,X_{nt}\rangle]-T_{t}\left(\mathrm{Var}_{\delta_{\cdot}}[\langle g,X_{(n-1)t}\rangle]\right)(x).
\end{align*}
The final equality is from Fubini's theorem. Thus we have
\begin{align}\label{lem2case2.1}
\mathrm{Var}_{\delta_{x}}[U_{2}(n,t)]=&\frac{1}{(nt)^{1+2r}}\e^{-\lambda_{1} nt}\mathrm{Var}_{\delta_{x}}[\langle T_{(n-1)t}g,X_{t}\rangle]\nonumber\\
=&\frac{1}{(nt)^{1+2r}}\e^{-\lambda_{1} nt}\mathrm{Var}_{\delta_{x}}[\langle g,X_{nt}\rangle]\nonumber\\
&-\left(1-\frac{1}{n}\right)^{1+2r}\e^{-\lambda_{1} t}T_{t}\left(\frac{1}{((n-1)t)^{1+2r}}\e^{-\lambda_{1}(n-1)t}\mathrm{Var}_{\delta_{\cdot}}[\langle g,X_{(n-1)t}\rangle]\right)(x).\nonumber\\
\end{align}
Recall that $\frac{1}{t^{1+2r}}\e^{-\lambda_{1} t}\varphi(y)^{-1}\mathrm{Var}_{\delta_{y}}[\langle g,X_{t}\rangle]\stackrel{b}{\to}\frac{\varrho^{2}_{g^{*}}}{1+2r}$ as $t\to+\infty$. It follows by Assumption \ref{AS1} and the dominated convergence theorem that
\begin{align*}
&\lim_{t\to+\infty}\e^{-\lambda_{1} t}T_{t}\left(\frac{1}{((n-1)t)^{1+2r}}\e^{-\lambda_{1}(n-1)t}\mathrm{Var}_{\delta_{\cdot}}[\langle g,X_{(n-1)t}\rangle]\right)(x)\\
=&\varphi(x)\langle \frac{\varrho^{2}_{g^{*}}}{1+2r}\varphi,\widetilde{\varphi}\rangle=\frac{\varrho^{2}_{g^{*}}}{1+2r}\varphi(x).
\end{align*}
Hence we get from \eqref{lem2case2.1} that
\begin{equation}\nonumber
\lim_{t\to+\infty}\mathrm{Var}_{\delta_{x}}[U_{2}(n,t)]=\left(1-\left(1-\frac{1}{n}\right)^{1+2r}\right)\frac{\varrho^{2}_{g^{*}}}{1+2r}\varphi(x).
\end{equation}
On the other hand, under our assumption,
$$\pp_{\delta_{x}}\left[U_{2}(n,t)\right]=\frac{1}{(nt)^{1/2+r}}\e^{-\frac{\lambda_{1}}{2}n t}T_{nt}g(x)\to 0\mbox{ as }t\to+\infty.$$
Thus we get $\lim_{n\to+\infty}\lim_{t\to+\infty}\pp_{\delta_{x}}[U_{2}(n,t)^{2}]=\lim_{n\to+\infty}\lim_{t\to+\infty}
\mathrm{Var}_{\delta_{x}}[U_{2}(n,t)]=0$, which yields that $\lim_{n\to+\infty}\lim_{t\to+\infty}U_{2}(n,t)=0$ in $L^{2}(\pp_{\delta_{x}})$.
By Slutsky theorem, \[\frac{1}{(nt)^{1/2+r}}\e^{-\frac{\lambda_{1}}{2}nt}\langle g,X_{nt}\rangle=U_{1}(n,t)+U_{2}(n,t)\] converges in distribution to $\frac{\varrho_{g^{*}}}{\sqrt{1+2r}}\sqrt{{ W^{\varphi}_{\infty}}}N,\mbox{ as first }t\to+\infty$ and then $n\to+\infty$. Hence we have proved (ii).

Applying similar argument as in the proof of Lemma \ref{lem3case2} (with \eqref{lem2case2.0} replaced by \eqref{lem4case2.0}), we can show that for all $d,\delta>0$,
\begin{equation}\label{eq4.21}
\frac{1}{t}\e^{-\lambda_{1} t}\int_{\frac{1}{\sqrt{t}}\e^{-\frac{\lambda_{1}}{2}t}|{\langle g^{*},\nu\rangle}|>d\e^{\delta t}}{\langle g^{*},\nu\rangle}^{2}L_{t}(y,{ \ud \nu})\stackrel{b}{\to}0\mbox{ as }t\to+\infty.
\end{equation}
Then the first conclusion of (i) follows from \eqref{eq4.21} (used in place of \eqref{lem3case2.0}) in the same way as for (ii). Now suppose \eqref{lemnew.1} holds for $r=0$. In this case, we have $\e^{-\frac{\lambda_{1}}{2}t}T_{t}(g-g^{*})\rightrightarrows 0$ as $t\to+\infty$. It then follows by Lemma \ref{lem:newcase2}(i) that
$$\frac{1}{t}\e^{-\lambda_{1} t}\varphi(x)^{-1}\mathrm{Var}_{\delta_{x}}[\langle g-g^{*},X_{t}\rangle]\to 0\mbox{ as }t\to+\infty.$$
We also note that
$$\frac{1}{\sqrt{t}}\e^{-\frac{\lambda_{1}}{2}t}\pp_{\delta_{x}}\left[\langle g-g^{*},X_{t}\rangle\right]=\frac{1}{\sqrt{t}}\e^{-\frac{\lambda_{1}}{2}t}T_{t}(g-g^{*})(x)\to 0\mbox{ as }t\to+\infty.$$
Thus $\frac{1}{t}\e^{-\lambda_{1} t}\pp_{\delta_{x}}\left[\langle g-g^{*},X_{t}\rangle^{2}\right]\to 0$ as $t\to+\infty$. This implies that $\lim_{t\to+\infty}\frac{1}{\sqrt{t}}\e^{-\frac{\lambda_{1}}{2} t}\langle g-g^{*},X_{t}\rangle=0$ in $L^{2}(\pp_{\delta_{x}})$. Hence the second conclusion of (i) follows from the first conclusion and the Slutsky theorem.\qed

{\subsection{Large branching case $\epsilon(t)<\lambda_{1}/2$}}

\begin{proposition}\label{lem3case3}
Suppose \eqref{lemnew.1} holds for $\epsilon\in [0,\lambda_{1}/2)$.
Define $ W^{g^{*}}_{t}:=\e^{-(\lambda_{1}-\epsilon)t}\langle g^{*},X_{t}\rangle$ for $t\ge 0$. Then the following are true.
\begin{description}
\item{(i)} For every $\mu\in\mf$, $\{ W^{g^{*}}_{t}:t\ge 0\}$ is an $L^{2}(\pp_{\mu})$-bounded $\mathcal{F}_{t}$-martingale, hence converges $\pp_{\mu}$-a.s. and in $L^{2}(\pp_{\mu})$ to a limit $ W^{g^{*}}_{\infty}$ as $t\to+\infty$.
\item{(ii)} $ W^{g^{*}}_{\infty}$ is non-degenerate for
all
$\mu\in\mfo$ if and only if $\varrho^{2}_{g^{*}}>0$, where $\varrho^{2}_{g^{*}}=\langle \vartheta[g^{*}],\widetilde{\varphi}\rangle$.
\item{(iii)} For every $x\in E$, under $\pp_{\delta_{x}}$,
$$\e^{(\frac{\lambda_{1}}{2}-\epsilon)t}\left( W^{g^{*}}_{t}- W^{g^{*}}_{\infty}\right)
\stackrel{d}{\to}\frac{\varrho_{g^{*}}}{\sqrt{\lambda_{1}-2\epsilon}}\sqrt{W^{\varphi}_{\infty}}N\mbox{ as }t\to+\infty,$$
where $N$ is a standard normal random variable independent of ${ W^{\varphi}_{\infty}}$.
\item{(iv)} For every $x\in E$,
$$\lim_{t\to+\infty}\frac{1}{t^{r}}\e^{-(\lambda_{1}-\epsilon)t}\langle g,X_{t}\rangle= W^{g^{*}}_{\infty}\quad\mbox{ in }L^{2}(\pp_{\delta_{x}}).$$
\end{description}
\end{proposition}

\begin{proof}
(i) We have proved in Remark \ref{rm1} that $\pp_{\delta_{x}}\left[W^{g^{*}}_{t}\right]=\e^{-(\lambda_{1}-\epsilon)t}T_{t}g^{*}(x)=g^{*}(x)$ for all $t\ge 0$ and $x\in E$. Using this and the Markov property, we can show that $\{ W^{g^{*}}_{t}:t\ge 0\}$ is a $\mathcal{F}_{t}$-martingale. Moreover, by the bounded convergence theorem, we have $\lim_{t\to+\infty}\mathrm{Var}_{\mu}[ W^{g^{*}}_{t}]=\lim_{t\to+\infty}\langle \mathrm{Var}_{\delta_{\cdot}}[ W^{g^{*}}_{t}],\mu\rangle=\lim_{t\to+\infty}\langle \e^{-2(\lambda_{1}-\epsilon)t}\mathrm{Var}_{\delta_{\cdot}}[\langle g^{*},X_{t}\rangle],\mu\rangle=\langle \delta^{2}_{g^{*}},\mu\rangle$. It follows that $\pp_{\mu}\left[( W^{g^{*}}_{t})^{2}\right]=\mathrm{Var}_{\mu}[ W^{g^{*}}_{t}]+\langle g^{*},\mu\rangle^{2}$ is bounded from above for large $t$.

\smallskip

\noindent(ii) It suffices to prove that $\mathrm{Var}_{\mu}[ W^{g^{*}}_{\infty}]>0$ for all $\mu\in\mfo$ if and only if $\varrho^{2}_{g^{*}}>0$.

``$\Rightarrow$" is obvious since $\mathrm{Var}_{\widetilde{\varphi}}[ W^{g^{*}}_{\infty}]=\langle \delta^{2}_{g^{*}},\widetilde{\varphi}\rangle=\frac{1}{\lambda_{1}-2\epsilon}\varrho^{2}_{g^{*}}.$

$\Leftarrow$: Suppose $\varrho^{2}_{g^{*}}>0$. Since $\lim_{s\to+\infty}\triangle_{s}=0$, there is some $s_{0}>0$ such that $\triangle_{s}\le \frac{1}{2}\frac{\varrho^{2}_{g^{*}}}{\|\vartheta[g^{*}]\|_{\infty}}$ for all $s\ge s_{0}$. Then by \eqref{ieq:H3}, we have for every $x\in E$ and $s\ge s_{0}$,
$$\varphi(x)^{-1}\e^{-\lambda_{1} s}T_{s}(\vartheta[g^{*}])(x)\ge\langle \vartheta[g^{*}],\widetilde{\varphi}\rangle-\|\vartheta[g^{*}]\|_{\infty}\triangle_{s}\ge \frac{1}{2}\varrho^{2}_{g^{*}}>0.$$
Thus we get
$$\mathrm{Var}_{\delta_{x}}[ W^{g^{*}}_{\infty}]=\delta^{2}_{g^{*}}(x)\ge \int_{s_{0}}^{+\infty}\e^{-(\lambda_{1}-2\epsilon)s}\e^{-\lambda_{1} s}T_{s}(\vartheta[g^{*}])(x){ \ud s}\ge \frac{1}{2(\lambda_{1}-2\epsilon)}\varrho^{2}_{g^{*}}\varphi(x)\e^{-(\lambda_{1}-2\epsilon)s_{0}}$$
for all $x\in E$. Immediately, $\mathrm{Var}_{\mu}[ W^{g^{*}}_{\infty}]=\langle\mathrm{Var}_{\delta_{\cdot}}[ W^{g^{*}}_{\infty}],\mu\rangle>0$ for all $\mu\in\mfo$.

\smallskip

\noindent(iii) The proof of this assertion is similar to that of Lemma \ref{prop1}. For the reader's convenience we also give details here.
We observe that the law of $ W^{g^{*}}_{\infty}$ under $\pp_{\delta_{x}}$ is infinitely divisible with mean $g^{*}(x)$ and variance $\delta^{2}_{g^{*}}(x)$.
Let $\Phi^{*}_{x}(\theta)$ be the corresponding characteristic exponent, that is,
$$\e^{-\Phi^{*}_{x}(\theta)}=\pp_{\delta_{x}}\left[\e^{i\theta  W^{g^{*}}_{\infty}}\right]\quad\forall \theta\in\R.$$
We can represent $\Phi^{*}_{x}(\theta)$ as
$$\Phi^{*}_{x}(\theta)=-i\theta g^{*}(x)+\frac{1}{2}\delta^{2}_{g^{*}}(x)\theta^{2}+\int_{\R}\left(1-\e^{i\theta r}+i\theta r-\frac{1}{2}\theta^{2}r^{2}\right)\Upsilon^{*}(x,{ \ud r}),$$
where $r^{2}\Upsilon^{*}(x,{ \ud r})$ is a bounded measure on $\R$ such that $\int_{\R}r^{2}\Upsilon^{*}(x,{ \ud r})=\delta^{2}_{g^{*}}(x)$. By the Markov property we have
\begin{align}
&\pp_{\delta_{x}}\big[\exp\{i\theta \e^{(\frac{\lambda_{1}}{2}-\epsilon)t}( W^{g^{*}}_{t}- W^{g^{*}}_{\infty})\}\big]\nonumber\\
=&\pp_{\delta_{x}}\big[\exp\{i\theta\e^{-\frac{\lambda_{1}}{2}t}\langle g^{*},X_{t}\rangle\}\pp_{X_{t}}\big[\exp\{-i\theta\e^{-\frac{\lambda_{1}}{2}t} W^{g^{*}}_{\infty}\}\big]\big]\nonumber\\
=&\pp_{\delta_{x}}\big[\exp\{-\e^{-\lambda_{1} t}\langle\frac{1}{2}\theta^{2}\delta^{2}_{g^{*}}\nonumber\\
&\quad+\e^{\lambda_{1} t}
\int_{\R}\big(1-\exp\{-i\theta\e^{-\frac{\lambda_{1}}{2}t}r\}-i\theta\e^{-\frac{\lambda_{1}}{2}t}r-\frac{1}{2}\theta^{2}\e^{-\lambda_{1} t}r^{2}\big)\Upsilon^{*}(\cdot,{ \ud r}),X_{t}\rangle\}\big]\nonumber\\
=:&\pp_{\delta_{x}}\big[\exp\{-\e^{-\lambda_{1} t}\langle U(\theta,t),X_{t}\rangle\}\big].\label{eq4.26}
\end{align}
We note that for every $y\in E$,
\begin{eqnarray}\label{lem3case3.1}
&&\left|\int_{\R}\left(1-\exp\{-i\theta\e^{-\frac{\lambda_{1}}{2}t}r\}-i\theta\e^{-\frac{\lambda_{1}}{2}t}r-\frac{1}{2}\theta^{2}\e^{-\lambda_{1} t}r^{2}\right)\Upsilon^{*}(y,{ \ud r})\right|\nonumber\\
&\le& |\theta|^{2}\e^{-\lambda_{1} t}\int_{\R}r^{2}\left(1\wedge \frac{|\theta|\e^{-\frac{\lambda_{1}}{2}t}|r|}{6}\right)\Upsilon^{*}(y,{ \ud r}).
\end{eqnarray}
Since $\int_{\R}r^{2}\Upsilon^{*}(y,{ \ud r})=\delta^{2}_{g^{*}}(y)$ is a bounded function on $E$, \eqref{lem3case3.1} yields that
$$\e^{\lambda_{1} t}\int_{\R}\left(1-\exp\{-i\theta\e^{-\frac{\lambda_{1}}{2}t}r\}-i\theta\e^{-\frac{\lambda_{1}}{2}t}r-\frac{1}{2}\theta^{2}\e^{-\lambda_{1} t}r^{2}\right)\Upsilon^{*}(y,{ \ud r})\stackrel{b}{\to}0\mbox{ as }t\to+\infty.$$
It then follows by Proposition \ref{lem2.2}(i) that
$$\lim_{t\to+\infty}\e^{-\lambda_{1} t}\langle U(\theta,t),X_{t}\rangle= \langle \frac{1}{2}\theta^{2}\delta^{2}_{g^{*}},\widetilde{\varphi}\rangle { W^{\varphi}_{\infty}}=\frac{1}{2}\theta^{2}\frac{\varrho^{2}_{g^{*}}}{\lambda_{1}-2\epsilon}W^{\varphi}_{\infty}\mbox{ in }L^{2}(\pp_{\delta_{x}}).$$
It thus follows from \eqref{eq4.26} that
$$\lim_{t\to+\infty}\pp_{\delta_{x}}\left[\exp\{i\theta\e^{(\frac{\lambda_{1}}{2}-\epsilon)t}( W^{g^{*}}_{t}- W^{g^{*}}_{\infty})\}\right]=\pp_{\delta_{x}}\left[\exp\{-\frac{1}{2}\theta^{2}
\frac{\varrho^{2}_{g^{*}}}{\lambda_{1}-2\epsilon}W^{\varphi}_{\infty}\}\right].$$
Hence we prove (iii).

\smallskip

\noindent (iv) For $t>0$, let $\tilde{g}_{t}:=t^{-r}g-g^{*}$, then $\e^{-(\lambda_{1}-\epsilon)t}\langle\tilde{g}_{t},X_{t}\rangle=t^{-r}\e^{-(\lambda_{1}-\epsilon)t}\langle g,X_{t}\rangle- W^{g^{*}}_{t}$.
It suffices to show that
\begin{equation}\nonumber
\lim_{t\to+\infty}\e^{-(\lambda_{1}-\epsilon)t}\langle\tilde{g}_{t},X_{t}\rangle=0\mbox{ in }L^{2}(\pp_{\delta_{x}}).
\end{equation}
We note that
$$\pp_{\delta_{x}}\left[\e^{-(\lambda_{1}-\epsilon)t}\langle\tilde{g}_{t},X_{t}\rangle\right]=\e^{-(\lambda_{1}-\epsilon)t}T_{t}\tilde{g}_{t}(x)=\frac{1}{t^{r}}\e^{-(\lambda_{1}-\epsilon)t}T_{t}g(x)-g^{*}(x)
\to 0\mbox{ as }t\to+\infty.$$
Hence we only need to show that
\begin{equation}
\lim_{t\to+\infty}\mathrm{Var}_{\delta_{x}}\left[\e^{-(\lambda_{1}-\epsilon)t}\langle\tilde{g},X_{t}\rangle\right]=0.\label{eq4.29}
\end{equation}
In fact, we have
\begin{align}
&\mathrm{Var}_{\delta_{x}}\left[\e^{-(\lambda_{1}-\epsilon)t}\langle\tilde{g}_{t},X_{t}\rangle\right]\nonumber\\
=&\e^{-2(\lambda_{1}-\epsilon)t}\int_{0}^{t}T_{s}\left(\vartheta[T_{t-s}\tilde{g}_{t}]\right)(x){ \ud s}\nonumber\\
=&\int_{0}^{t}\e^{-(\lambda_{1}-2\epsilon)s}\e^{-\lambda_{1} s}T_{s}\left(\vartheta\left[\left(1-\frac{s}{t}\right)^{r}\frac{\e^{-(\lambda_{1}-\epsilon)(t-s)}T_{t-s}g}{(t-s)^{r}}-g^{*}\right]\right)(x){ \ud s}.\label{lem3case3.3}
\end{align}
Given that $t^{-r}\e^{-(\lambda_{1}-\epsilon)t}T_{t}g-g^{*}\rightrightarrows 0$, we can easily show that for every $s>0$,
$$\left(1-\frac{s}{t}\right)^{r}\frac{\e^{-(\lambda_{1}-\epsilon)(t-s)}T_{t-s}g}{(t-s)^{r}}-g^{*}\rightrightarrows 0
\mbox{ and }\vartheta\Big[\left(1-\frac{s}{t}\right)^{r}\frac{\e^{-(\lambda_{1}-\epsilon)(t-s)}T_{t-s}g}{(t-s)^{r}}-g^{*}\Big]\stackrel{b}{\to}0,$$
as $t\to+\infty$. Thus, by the dominated convergence theorem, the integrand in the right hand side of \eqref{lem3case3.3} converges to $0$ as $t\to+\infty$.
Moreover, we can apply similar arguments (with minor modification) as in the proof of Lemma \ref{lem:newcase2}(ii) to show that, for large $t$, the integrand in the right hand side of \eqref{lem3case3.3} is bounded from above by some nonnegative integrable function on $(0,+\infty)$. Hence, by the dominated convergence theorem, the corresponding integral converges to $0$ as $t\to+\infty$, and \eqref{eq4.29} follows.
\end{proof}

\noindent\textbf{Proof of Theorem \ref{them2}:} (i) Obviously, the function $\widehat{f}=f-\langle f,\widetilde{\varphi}\rangle \varphi$ satisfies that $\langle \widehat{f},\widetilde{\varphi}\rangle=0$ and $\epsilon(\widehat{f})=\epsilon(f)$. Moreover, we have
\begin{equation}\label{eq4.32}
\e^{-\lambda_{1} t}\langle f,X_{t}\rangle-\langle f,\widetilde{\varphi}\rangle W^{\varphi}_{\infty}=\e^{-\lambda_{1} t}\langle \widehat{f},X_{t}\rangle+\langle f,\widetilde{\varphi}\rangle \left(W^{\varphi}_{t}-W^{\varphi}_{\infty}\right).
\end{equation}
Hence (i) is a directly result of Proposition \ref{lem2case1}.

\smallskip

\noindent(ii) Recall from Theorem \ref{them1} that $\e^{\frac{\lambda_{1}}{2}t}(W^{\varphi}_{t}-W^{\varphi}_{\infty})$ converges in distribution to some finite random variable as $t\to+\infty$. Hence for any $r\ge 0$,
$\frac{1}{t^{1/2+r}}\e^{\frac{\lambda_{1}}{2}t}(W^{\varphi}_{t}-W^{\varphi}_{\infty})$ converges in probability to $0$ as $t\to+\infty$. Thus, (ii) follows from \eqref{eq4.32}, Proposition \ref{lem2case2} and the Slutsky theorem.

\smallskip

\noindent(iii) Given that $t^{-r}\e^{-(\lambda_{1}-\epsilon)t}T_{t}\widehat{f}\rightrightarrows f^{*}$ for $\epsilon\in [0,\lambda_{1}/2)$, we have by Proposition \ref{lem3case3} that $t^{-r}\e^{-(\lambda_{1}-\epsilon)t}\langle\widehat{f},X_{t}\rangle\to W^{f^{*}}_{\infty}$ as $t\to+\infty$ in $L^{2}(\pp_{\delta_{x}})$.
In view of \eqref{eq4.32}, to prove the first conclusion of (iii), we only need to show that
\begin{equation}\label{eq4.33}
\lim_{t\to+\infty}\frac{\e^{\epsilon t}}{t^{r}}\left(W^{\varphi}_{t}-W^{\varphi}_{\infty}\right)=0\mbox{ in }L^{2}(\pp_{\delta_{x}}).
\end{equation}
Since $\lim_{t\to+\infty}W^{\varphi}_{t}=W^{\varphi}_{\infty}$ in $L^{2}(\pp_{\delta_{x}})$, we have $W^{\varphi}_{t}=\pp_{\delta_{x}}\left[W^{\varphi}_{\infty}|\mathcal{F}_{t}\right]$ $\pp_{\delta_{x}}$-a.s., and $$\pp_{\delta_{x}}\left[W^{\varphi}_{t}W^{\varphi}_{\infty}\right]=\pp_{\delta_{x}}\left[W^{\varphi}_{t}\pp_{\delta_{x}}\left[W^{\varphi}_{\infty}|\mathcal{F}_{t}\right]\right]=\pp_{\delta_{x}}
\left[(W^{\varphi}_{t})^{2}\right].$$
Thus
\begin{align*}
\pp_{\delta_{x}}\left[\left(W^{\varphi}_{t}-W^{\varphi}_{\infty}\right)^{2}\right]=&\pp_{\delta_{x}}\left[(W^{\varphi}_{\infty})^{2}\right]-\pp_{\delta_{x}}\left[(W^{\varphi}_{t})^{2}\right]\\
=&\mathrm{Var}_{\delta_{x}}\left[W^{\varphi}_{\infty}\right]-\mathrm{Var}_{\delta_{x}}\left[W^{\varphi}_{t}\right]=\int_{t}^{+\infty}\e^{-2\lambda_{1} s}T_{s}(\vartheta[\varphi])(x){ \ud s}.
\end{align*}
It follows by \eqref{eq1.9} that there is some $c_{1}>0$ such that $\e^{-\lambda_{1} s}T_{s}(\vartheta[\varphi])(x)\le c_{1}\varphi(x)$ for all $x\in E$ and $s\ge 0$.
Thus $\int_{t}^{+\infty}\e^{-2\lambda_{1} s}T_{s}(\vartheta[\varphi])(x){ \ud s}\le c_{1}\varphi(x)\int_{t}^{+\infty}\e^{-\lambda_{1} s}{ \ud s}=\frac{c_{1}}{\lambda_{1}}\e^{-\lambda_{1} t}\varphi(x)$.
Hence we get
$$\frac{e^{2\epsilon t}}{t^{2r}}\pp_{\delta_{x}}\left[\left(W^{\varphi}_{t}-W^{\varphi}_{\infty}\right)^{2}\right]\le \frac{c_{1}}{\lambda_{1}}\frac{\e^{(2\epsilon-\lambda_{1})t}}{t^{2r}}\varphi(x)\to 0\mbox{ as }t\to+\infty,$$
and \eqref{eq4.33} follows immediately. The second and third conclusions of (iii) follow directly from Proposition \ref{lem3case3}(ii) and (iii).\qed

\section{Proofs of the propositions in Section \ref{sec1.3}}\label{sec5}

\noindent\textbf{Proof of Proposition \ref{prop:a}:}
{ Recall the definitions of $\alpha(f),\gamma(f),s(f)$ and $\mathfrak{I}(f)$ given in Section \ref{sec1.3}.
For $j\in \mathfrak{I}(f)$, we further define $\mathfrak{M}_{j}(f):=\{n:\ 1\le n\le n_{j}\mbox{ and }r^{(j)}_{n}(f)=\gamma(f)\}$.}
In view of \eqref{eg1.5}, we may write $T_{t}f$ as
\begin{equation}\label{eg1.6}
T_{t}f=I_{t}(f)+II_{t}(f)+III_{t}(f)+T_{t}\widetilde{f}_{s(f)},
\end{equation}
where
\begin{align}
I_{t}(f):=&\sum_{j\in \mathfrak{I}(f)}\e^{\lambda_{j}t}\sum_{n\in \mathfrak{M}_{j}(f)}\big(\langle f,\widehat{\phi}^{(j)}_{n}\rangle_{m}+t\langle f,\widehat{\phi}^{(j)}_{n+1}\rangle_{m}+\cdots+
\frac{t^{\gamma(f)}}{\gamma(f)!}\langle f,\widehat{\phi}^{(j)}_{n+\gamma(f)}\rangle_{m}\big)\phi^{(j)}_{n},\nonumber\\
II_{t}(f):=&\sum_{j\in \mathfrak{I}(f)}\e^{\lambda_{j}t}\sum_{n\in \{1,\cdots,n_{j}\}\setminus \mathfrak{M}_{j}(f)}\big(\langle f,\widehat{\phi}^{(j)}_{n}\rangle_{m}+t\langle f,\widehat{\phi}^{(j)}_{n+1}\rangle_{m}+\cdots+
\frac{t^{\gamma^{(j)}_{n}(f)}}{\gamma^{(j)}_{n}(f)!}\langle f,\widehat{\phi}^{(j)}_{n+\gamma^{(j)}_{n}(f)}\rangle_{m}\big)\phi^{(j)}_{n},\nonumber\\
III_{t}(f):=&\sum_{j\in \{2,\cdots,s(f)\}\setminus \mathfrak{I}(f)}\e^{\lambda_{j}t}\Phi^{\mathrm{T}}_{j}D_{j}(t)\langle f,\widehat{\Phi}_{j}\rangle_{m}.\nonumber
\end{align}
We have $\re \lambda_{j}=\alpha(f)$ for $j\in \mathfrak{I}(f)$ and $r^{(j)}_{n}(f)<\gamma(f)$ for $n\in \{1,\cdots,n_{j}\}\setminus \mathfrak{M}_{j}(f)$.
Using this and the boundedness of $\phi^{(j)}_{n}$, we can easily prove that
$$\frac{\e^{-\alpha(f)t}}{t^{\gamma(f)}}(I_{t}(f)+II_{t}f)-\frac{1}{\gamma(f)!}\sum_{j\in \mathfrak{I}(f)}\e^{i t\im \lambda_{j}}F_{j}(f)\rightrightarrows 0\mbox{ on $E$, as }t\to+\infty,$$
where $F_{j}(f)=\sum_{n\in \mathfrak{M}_{j}(f)}\langle f,\widehat{\phi}^{(j)}_{n+\gamma(f)}\rangle_{m}\phi^{(j)}_{n}$.
For $j\in \{2,\cdots,s(f)\}\setminus \mathfrak{I}(f)$, we have either $\re \lambda_{j}<\alpha(f)$, or $\re\lambda_{j}=\alpha(f)$ and $\max_{1\le n\le n_{j}}r^{(j)}_{n}(f)<\gamma(f)$. In either case,
$\frac{\e^{-\alpha(f)t}}{t^{\gamma(f)}}\e^{\lambda_{j}t}\Phi^{\mathrm{T}}_{j}D_{j}(t)\langle f,\widehat{\Phi}_{j}\rangle_{m}\rightrightarrows 0\mbox{ on $E$ as }t\to+\infty.$ Thus $\frac{\e^{-\alpha(f)t}}{t^{\gamma(f)}}III_{t}(f)\rightrightarrows 0$ as $t\to+\infty$. Therefore, to prove Proposition \ref{prop:a}, it suffices to prove that
\begin{equation}\label{a.1}
\frac{\e^{-\alpha(f)t}T_{t}\widetilde{f}_{s(f)}}{t^{\gamma(f)}}\rightrightarrows 0\mbox{ as }t\to+\infty.
\end{equation}
We note that for each $k\ge 1$, $\widehat{\mathcal{N}}^{\perp}_{k}$ is an invariant subspace of $T_{t}$. In fact, for $g\in\widehat{\mathcal{N}}^{\perp}_{k}$, we have $\langle g,\widehat{\Phi}_{j}\rangle_{m}=0$ for all $1\le j\le k$, and
thus $\langle T_{t}g,\widehat{\Phi}_{j}\rangle_{m}=\langle g,\widehat{T}_{t}\widehat{\Phi}_{j}\rangle_{m}=\e^{\overline{\lambda_{j}}t}\langle g,D_{j}(t)\widehat{\Phi}_{j}\rangle_{m}=0$, which yields that $T_{t}g\in\widehat{\mathcal{N}}^{\perp}_{k}$.
Thus $\{\left.T_{t}\right|_{\widehat{\mathcal{N}}^{\perp}_{k}}:t\ge 0\}$ is a semigroup on $\widehat{\mathcal{N}}^{\perp}_{k}$. By \cite[Theorem 6.7.5]{BP}, we have $\sigma(\left.T_{t}\right|_{\widehat{\mathcal{N}}^{\perp}_{k}})\setminus\{0\}=\{\e^{\lambda_{j}t}:\ j\ge k+1\}$.
So the spectral radius of $\left.T_{t}\right|_{\widehat{\mathcal{N}}^{\perp}_{k}}$ is $r(\left.T_{t}\right|_{\widehat{\mathcal{N}}^{\perp}_{k}})=\e^{t\re \lambda_{k+1}}$.
$T_{t}$ is a compact operator in $L^{2}(E,m,\mathbb{C})$, and hence is compact in $\widehat{\mathcal{N}}^{\perp}_{k}$. Thus the operator norm $\|\left.T_{t}\right|_{\widehat{\mathcal{N}}^{\perp}_{k}}\|_{2}$ is equal to $r(\left.T_{t}\right|_{\widehat{\mathcal{N}}^{\perp}_{k}})$.
By H\"{o}lder's inequality and \eqref{eg1.2}, we have for any $g\in \widehat{\mathcal{N}}^{\perp}_{k}$, $t,s>0$ and any $u=(i,x)\in E$, there is a constant $c_{1}=c_{1}(s)>0$ such that,
\begin{align}\label{a.2}
|T_{t+s}g(u)|=&|\int_{E}\mathfrak{p}(s,u,v)T_{t}g(v)m({ \ud v})|\nonumber\\
\le&\|T_{t}g\|_{2}\left(\int_{E}\mathfrak{p}(s,u,v)^{2}m({ \ud v})\right)^{1/2}
\le c_{1}\delta_{D}(x)\|g\|_{2}\e^{t\re \lambda_{k+1}}.
\end{align}
This implies that for any $a>\re\lambda_{k+1}$, $\e^{-a t}T_{t}g\rightrightarrows 0$ as $t\to+\infty$. Thus \eqref{a.1} follows, since $\alpha(f)>\re\lambda_{s(f)+1}$, $\widetilde{f}_{s(f)}\in \widehat{\mathcal{N}}^{\perp}_{s(f)}$ and $\gamma(f)\ge 0$.\qed

\noindent\textbf{Proof of Proposition \ref{prop:d}(i):} Suppose $\alpha(f)=\lambda_{1}/2$. To simplify the notation, we write $\gamma$ for $\gamma(f)$. For every $u\in E$,
\begin{align*}
\frac{\e^{-\lambda_{1}t}}{t^{1+2\gamma}}\varphi(u)^{-1}\mathrm{Var}_{\delta_{u}}[\langle f,X_{t}\rangle]
=&\frac{\e^{-\lambda_{1}t}}{t^{1+2\gamma}}\varphi(u)^{-1}\int_{0}^{t}T_{t-s}(\vartheta[T_{s}f])(u){ \ud s}\\
=&\varphi(u)^{-1}\int_{0}^{1}r^{2\gamma}\e^{-\lambda_{1}t(1-r)}T_{t(1-r)}\left(\vartheta\left[\frac{\e^{-\frac{\lambda_{1}}{2}t r}T_{t r}f}{(t r)^{\gamma}}\right]\right)(u){ \ud r}.
\end{align*}
The second equality follows from change of variables.
For $t>0$, let
$$\varepsilon_{t}(f):=\frac{\e^{-\frac{\lambda_{1}}{2}t}T_{t}f}{t^{\gamma}}-\frac{1}{\gamma!}\sum_{j\in \mathfrak{I}(f)}\e^{i t\im \lambda_{j}}F_{j}(f).$$
We observe that $j'\in \mathfrak{I}(f)$ if $j\in \mathfrak{I}(f)$, and $\im\lambda_{j'}=-\im\lambda_{j}$, $F_{j'}(f)=\overline{F_{j}(f)}$. So $\sum_{j\in\mathfrak{I}(f)}\e^{it\im \lambda_{j}}F_{j}(f)$ is a real-valued function.
 Since $\alpha(f)=\lambda_{1}/2$, Proposition \ref{prop:a} yields that $\varepsilon_{t}(f)\rightrightarrows 0$ as $t\to+\infty$.
 Using the fact that $\vartheta[h+g]=\vartheta[h]+2\vartheta[h,g]+\vartheta[g]$ for $h,g\in\mathcal{B}_{b}(E)$,
 we can rewrite $\frac{\e^{-\lambda_{1}t}}{t^{1+2\gamma}}\varphi(u)^{-1}\mathrm{Var}_{\delta_{u}}[\langle f,X_{t}\rangle]$ as $I_{t}^{f}(u)+II_{t}^{f}(u)+III_{t}^{f}(u)$, where
\begin{align*}
I_{t}^{f}(u):=&\int_{0}^{1}r^{2\gamma}\varphi(u)^{-1}\e^{-\lambda_{1}t(1-r)}T_{t(1-r)}\left(\vartheta[\varepsilon_{tr}(f)]\right)(u){ \ud r},\\
II_{t}^{f}(u):=&2\int_{0}^{1}r^{2\gamma}\varphi(u)^{-1}\e^{-\lambda_{1}t(1-r)}T_{t(1-r)}\Big(\vartheta\big[\frac{1}{\gamma!}\sum_{j\in \mathfrak{I}(f)}\e^{i t r\im \lambda_{j}}F_{j}(f),\varepsilon_{tr}(f)\big]\Big)(u){ \ud r},\\
III_{t}^{f}(u):=&\int_{0}^{1}r^{2\gamma}\varphi(u)^{-1}\e^{-\lambda_{1}t(1-r)}T_{t(1-r)}\Big(\vartheta\big[\frac{1}{\gamma!}\sum_{j\in \mathfrak{I}(f)}\e^{i t r\im \lambda_{j}}F_{j}(f)\big]\Big)(u){ \ud r}.
\end{align*}
For $I^{f}_{t}(u)$, we can apply similar argument as in the proof of Lemma \ref{lem:newcase2}(i) (with $\e^{-\lambda_{1}tu}T_{tu}g/(tu)^{r}$ replaced by $\varepsilon_{tu}(f)$) to show that $I^{f}_{t}(u)\stackrel{b}{\to}0$ as $t\to+\infty$.

It is easy to see from the definition of $\vartheta[\cdot,\cdot]$ that for any $h,g\in\mathcal{B}_{b}(E)$ and $u\in E$,
$$\vartheta[h,g](u)\le \|h\|_{\infty}\|g\|_{\infty}\vartheta[1](u).$$
Thus we have
$$\vartheta\big[\frac{1}{\gamma!}\sum_{j\in \mathfrak{I}(f)}\e^{i t r\im \lambda_{j}}F_{j}(f),\varepsilon_{tr}(f)\big](u)
\le \big(\frac{1}{\gamma !}\sum_{j\in \mathfrak{I}(f)}\|F_{j}(f)\|_{\infty}\big)\|\varepsilon_{tr}(f)\|_{\infty}\vartheta[1](u).$$
Since $\vartheta[1](u)$ is a bounded function on $E$, we have $\vartheta\big[\frac{1}{\gamma!}\sum_{j\in \mathfrak{I}(f)}\e^{i t r\im \lambda_{j}}F_{j}(f),\varepsilon_{tr}(f)\big]\rightrightarrows 0$ as $t\to+\infty$. Then we can show that $II^{f}_{t}(u)\stackrel{b}{\to} 0$ as $t\to+\infty$ similarly as we did for $I^{f}_{t}(u)$.
Next we deal with $III^{f}_{t}(u)$.
By the definition of $\vartheta[\cdot,\cdot]$ and the fact that $\sum_{k\in \mathfrak{I}}\e^{itr\im\lambda_{k}}F_{k}(j)=\sum_{k\in \mathfrak{I}}\e^{-itr\im\lambda_{k}}\overline{F_{k}(j)}$, we have
\begin{align*}
&\vartheta\big[\frac{1}{\gamma!}\sum_{j\in \mathfrak{I}(f)}\e^{i t r\im \lambda_{j}}F_{j}(f)\big](u)\\
&=\frac{1}{(\gamma!)^{2}}\sum_{j\in\mathfrak{I}(f)}\sum_{k\in\mathfrak{I}(f)}\e^{i t r(\im \lambda_{j}-\im \lambda_{k})}{ \vartheta[F_{j}(f),\overline{F_{k}(f)}]}(u)\\
&=\frac{1}{(\gamma!)^{2}}\big(\sum_{j\in \mathfrak{I}(f)}{ \vartheta[F_{j}(f),\overline{F_{j}(f)}]}(u)+\sum_{j\in \mathfrak{I}(f)}\sum_{k\in \mathfrak{I}(f)\setminus\{j\}}\e^{i t r(\im\lambda_{j}-\im \lambda_{k})}{ \vartheta[F_{j}(f),\overline{F_{k}(f)}]}(u)\big).
\end{align*}
Hence we have $III^{f}_{t}(u)=III^{f,1}_{t}(u)+III^{f,2}_{t}(u)$,
where
\begin{align*}
III^{f,1}_{t}(u):=&\frac{1}{(\gamma !)^{2}}\sum_{j\in \mathfrak{I}(f)}\int_{0}^{1}r^{2\gamma}\varphi(u)^{-1}
\e^{-\lambda_{1}t(1-r)}T_{t(1-r)}({ \vartheta[F_{j}(f),\overline{F_{j}(f)}]})(u){ \ud r},\\
III^{f,2}_{t}(u):=&\frac{1}{(\gamma !)^{2}}\sum_{j\in \mathfrak{I}(f)}\sum_{k\in \mathfrak{I}(f)\setminus\{j\}}\int_{0}^{1}r^{2\gamma}\e^{i t r(\im \lambda_{j}-\im \lambda_{k})}\varphi(u)^{-1}
\e^{-\lambda_{1}t(1-r)}T_{t(1-r)}({ \vartheta[F_{j}(f),\overline{F_{k}(f)}]})(u){ \ud r}.
\end{align*}
Using the dominated convergence theorem and the fact that
$$\varphi(u)^{-1}
\e^{-\lambda_{1}t}T_{t}({ \vartheta[F_{j}(f),\overline{F_{j}(f)}]})(u)\stackrel{b}{\to}\langle { \vartheta[F_{j}(f),\overline{F_{j}(f)}]},\widehat{\varphi}\rangle_{m}\mbox{ as }t\to+\infty,$$
we can show that
$III^{f,1}_{t}(u)\stackrel{b}{\to}(1+2\gamma)^{-1}(\gamma !)^{-2}\sum_{j\in \mathfrak{I}(f)}\langle { \vartheta[F_{j}(f),\overline{F_{j}(f)}]},\widehat{\varphi}\rangle_{m}$ as $t\to+\infty$.
On the other hand, by the dominated convergence theorem and the fact that $\varphi(u)^{-1}
\e^{-\lambda_{1}t}T_{t}({ \vartheta[F_{j}(f),\overline{F_{k}(f)}]})(u)\stackrel{b}{\to}\langle { \vartheta[F_{j}(f),\overline{F_{k}(f)}]},\widehat{\varphi}\rangle_{m}$ as $t\to+\infty$, we have
\begin{align}\label{d.1.1}
&\big|\int_{0}^{1}r^{2\gamma}\e^{i t r(\im \lambda_{j}-\im \lambda_{k})}\varphi(u)^{-1}
\e^{-\lambda_{1}t(1-r)}T_{t(1-r)}({ \vartheta[F_{j}(f),\overline{F_{k}(f)}]})(u){ \ud r}\nonumber\\
&-\langle { \vartheta[F_{j}(f),\overline{F_{k}(f)}]},\widehat{\varphi}\rangle_{m}\int_{0}^{1}r^{2\gamma}\e^{i t r(\im \lambda_{j}-\im \lambda_{k})}{ \ud r}
\big|\nonumber\\
\le&\int_{0}^{1}r^{2\gamma}\big|\varphi(u)^{-1}
\e^{-\lambda_{1}t(1-r)}T_{t(1-r)}({ \vartheta[F_{j}(f),\overline{F_{k}(f)}]})(u)-\langle { \vartheta[F_{j}(f),\overline{F_{k}(f)}]},\widehat{\varphi}\rangle_{m}\big|{ \ud r}\stackrel{b}{\to}0
\end{align}
as $t\to+\infty$.
By the integration by parts formula we have
for any $0\not=\theta\in \R$ and $t>0$,
$$\big|\int_{0}^{1}r^{2\gamma}\e^{i\theta r}{ \ud r}\big|=\Big|\frac{\left.r^{2\gamma}\e^{i\theta t r}\right|^{1}_{0}-\int_{0}^{1}2\gamma r^{2\gamma-1}\e^{i\theta t r}{ \ud r}}{i\theta t}\Big|\le \frac{2}{|\theta t|}\to 0\mbox{ as }t\to+\infty.$$
This together with \eqref{d.1.1} yields that $$\int_{0}^{1}r^{2\gamma}\e^{i t r(\im \lambda_{j}-\im \lambda_{k})}\varphi(u)^{-1}
\e^{-\lambda_{1}t(1-r)}T_{t(1-r)}({ \vartheta[F_{j}(f),\overline{F_{k}(f)}]})(u){ \ud r}\stackrel{b}{\to}0\mbox{ as }t\to+\infty,$$ for $j\in \mathfrak{I}(f)$ and $k\in \mathfrak{I}(f)\setminus\{j\}$. Consequently, we have $III^{f,2}_{t}(u)\stackrel{b}{\to}0$. Hence we have proved that
\begin{equation}\label{d.1.2}
\frac{\e^{-\lambda_{1}t}}{t^{1+2\gamma}}\varphi(u)^{-1}\mathrm{Var}_{\delta_{u}}[\langle f,X_{t}\rangle]=I^{f}_{t}(u)+II^{f}_{t}(u)+III^{f,1}_{t}(u)+III^{f,2}_{t}(u)\stackrel{b}{\to}\frac{\varrho^{2}_{f}}{1+2\gamma}\mbox{ as }t\to+\infty,
\end{equation}
where $\varrho^{2}_{f}=(\gamma!)^{-2}\sum_{j\in \mathfrak{I}(f)}\langle { \vartheta[F_{j}(f),\overline{F_{j}(f)}]},\widehat{\varphi}\rangle_{m}$.

Proposition \ref{prop:a} also implies that for every $u\in E$,
$$\frac{\e^{-\frac{\lambda_{1}}{2}t}T_{t}f}{t^{\gamma+\frac{1}{2}}}(u)\to 0\mbox{ as }t\to+\infty.$$
Using this and \eqref{d.1.2}, we can repeat the proof of Proposition \ref{lem2case2}(ii) and Corollary \ref{cor4.7} to prove the second conclusion of Proposition \ref{prop:d}(i).\qed

\noindent\textbf{Proof of Proposition \ref{prop:d}(ii):}
Suppose $\alpha(f)>\lambda_{1}/2$. Fix an arbitrary $j\in\mathfrak{I}(f)$.
To prove $W^{(j)}_{t}(f)$ is a martingale, it suffices to prove that
$\pp_{\delta_{u}}\big[W^{(j)}_{t}(f)\big]=\e^{-\lambda_{j}t}T_{t}(F_{j}(f))(u)=F_{j}(f)(u)$ for all $t\ge 0$ and $u\in E$. Since $\mathcal{N}_{j,v_{j}}$ is an invariant subspace of $T_{t}$,
we have $T_{t}(F_{j}(f))\in\mathcal{N}_{j,v_{j}}$ for $F_{j}(f)=\sum_{n\in \mathfrak{M}_{j}(f)}\langle f,\widehat{\phi}^{(j)}_{n+\gamma(f)}\rangle_{m}\phi^{(j)}_{n}\in \mathcal{N}_{j,v_{j}}$.
Hence we only need to prove that for the basis $\{\phi^{(j)}_{k}:1\le k\le n_{j}\}$ of $\mathcal{N}_{j,v_{j}}$,
\begin{equation}\label{d.2.1}
\langle \e^{-\lambda_{j}t} T_{t}(F_{j}(f)),\phi^{(j)}_{k}\rangle_{m}=\langle F_{j}(f),\phi^{(j)}_{k}\rangle_{m}\quad\forall t\ge 0,\ 1\le k\le n_{j}.
\end{equation}
In fact, we have $T_{t}(F_{j}(f))=\sum_{n\in \mathfrak{M}_{j}(f)}\langle f,\widehat{\phi}^{(j)}_{n+\gamma(f)}\rangle_{m}T_{t}\phi^{(j)}_{n}$, where $T_{t}\phi^{(j)}_{n}$ is the $n$-th component of the vector $(T_{t}\Phi_{j})^{\mathrm{T}}=\e^{\lambda_{j}t}\Phi^{\mathrm{T}}_{j}D_{j}(t)$.
For $1\le l\le r_{j}$ and $n\in [\sum_{i=0}^{l-1}d_{j,i}+1,\sum_{i=0}^{l}d_{j,i}]$ (here we take the convention that $d_{j,0}=0$), $T_{t}\phi^{(j)}_{n}$ is given by
$$\e^{\lambda_{j}t}\Big(\frac{t^{n-\sum_{i=0}^{l-1}d_{j,i}-1}}{(n-\sum_{i=0}^{l-1}d_{j,i}-1)!}\phi^{(j)}_{\sum_{i=0}^{l-1}d_{j,i}+1}+\frac{t^{n-\sum_{i=0}^{l-1}d_{j,i}-2}}{(n-\sum_{i=0}^{l-1}d_{j,i}-2)!}
\phi^{(j)}_{\sum_{i=0}^{l-1}d_{j,i}+2}+\cdots+\phi^{(j)}_{n}\Big)$$
Thus, for $1\le k\le n_{j}$ and $n\in [\sum_{i=0}^{l-1}d_{j,i}+1,\sum_{i=0}^{l}d_{j,i}]$, we have
\begin{align*}
\langle f,\widehat{\phi}^{(j)}_{n+\gamma(f)}\rangle_{m}\langle T_{t}\phi^{(j)}_{n},\phi^{(j)}_{k}\rangle_{m}
=&1_{\{\sum_{i=0}^{l-1}d_{j,i}+1\le k\le n\}}\e^{\lambda_{j}t}\frac{t^{n-k}}{(n-k)!}\langle f,\widehat{\phi}^{(j)}_{n+\gamma(f)}\rangle_{m}\\
=&1_{\{k=n\}}\e^{\lambda_{j}t}\langle f,\widehat{\phi}^{(j)}_{n+\gamma(f)}\rangle_{m}.
\end{align*}
The final equality is from the fact that $\langle f,\widehat{\phi}^{(j)}_{n+\gamma(f)}\rangle_{m}=\langle f,\widehat{\phi}^{(j)}_{k+r^{(j)}_{k}(f)+(n+\gamma(f)-k-r^{(j)}_{k}(f))}\rangle_{m}=0$ for $k<n$.
Hence we get
$$\e^{-\lambda_{j}t}\langle T_{t}(F_{j}(f)),\phi^{(j)}_{k}\rangle_{m}=\e^{-\lambda_{j}t}\sum_{n\in \mathfrak{M}_{j}(f)}\langle f,\widehat{\phi}^{(j)}_{n+\gamma(f)}\rangle_{m}\langle T_{t}\phi^{(j)}_{n},\phi^{(j)}_{k}\rangle_{m}=1_{\{k\in \mathfrak{M}_{j}(f)\}}\langle f,\widehat{\phi}^{(j)}_{k+\gamma(f)}\rangle_{m}$$
for all $t\ge 0$,
and \eqref{d.2.1} follows immediately.

Next we prove that $\{W^{(j)}_{t}(f):t\ge 0\}$ is $L^{2}(\pp_{\delta_{u}})$-bounded for every $u\in E$. It follows from the second moment formula of superprocess that
\begin{align}
\pp_{\delta_{u}}\left[|W^{(j)}_{t}(f)|^{2}\right]
=&\big|\pp_{\delta_{u}}\left[W^{(j)}_{t}(f)\right]\big|^{2}+\mathrm{Cov}(W^{(j)}_{t}(f),\overline{W^{(j)}_{t}(f)})\nonumber\\
=&|F_{j}(f)(u)|^{2}+\int_{0}^{t}T_{s}\big(\vartheta[\e^{-\lambda_{j}t}T_{t-s}(F_{j}(f)),\overline{\e^{-\lambda_{j}t}T_{t-s}(F_{j}(f))}]\big)(u){ \ud s}\nonumber\\
=&|F_{j}(f)(u)|^{2}+\int_{0}^{t}\e^{-(2\alpha(f)-\lambda_{1})s}\cdot\e^{-\lambda_{1}s}T_{s}({ \vartheta[F_{j}(f),\overline{F_{j}(f)}]})(u){ \ud s}.\label{5.6}
\end{align}
The final equality is because $T_{s}(F_{j}(f))=\e^{\lambda_{j}s}F_{j}(f)$ for all $s\ge 0$ and $\re\lambda_{j}=\alpha(f)$ for $j\in\mathfrak{I}$.
Under Assumption \ref{AS1}, the function $u\mapsto\e^{-\lambda_{1}s}T_{s}({ \vartheta[F_{j}(f),\overline{F_{j}(f)}]})(u)$ is bounded from above for large $s$.
 Since $2\alpha(f)-\lambda_{1}>0$, by the dominated convergence theorem, the final integral in the right hand side of \eqref{5.6} converges to $\int_{0}^{+\infty}\e^{-2\alpha(f)s}T_{s}({ \vartheta[F_{j}(f),\overline{F_{j}(f)}]})(u){ \ud s}\in [0,+\infty)$ as $t\to+\infty$. Thus we get $\sup_{t\ge 0}\pp_{\delta_{u}}\left[|W^{(j)}_{t}(f)|^{2}\right]<+\infty$. Consequently, the martingale $W^{(j)}_{t}(f)$ converges $\pp_{\delta_{u}}$-a.s. and in $L^{2}(\pp_{\delta_{u}})$ to a limit random variable $W^{(j)}_{\infty}(f)$.

To prove the second conclusion of Proposition \ref{prop:d}(ii), we only need to show that
\begin{equation}\label{d.2.2}
\frac{\e^{-\alpha(f)t}}{t^{\gamma(f)}}\langle f,X_{t}\rangle-\frac{1}{\gamma(f)!}\sum_{j\in \mathfrak{I}(f)}\e^{i t\im \lambda_{j}}W^{(j)}_{t}(f)\to 0\mbox{ in }L^{2}(\pp_{\delta_{u}})\mbox{ as }t\to+\infty.
\end{equation}
Since $\re\lambda_{j}=\alpha(f)$ for $j\in\mathfrak{I}(f)$, we have
\begin{align}
&\mathrm{Var}_{\delta_{u}}\big[\frac{\e^{-\alpha(f)t}}{t^{\gamma(f)}}\langle f,X_{t}\rangle-\frac{1}{\gamma(f)!}\sum_{j\in \mathfrak{I}(f)}\e^{i t\im \lambda_{j}}W^{(j)}_{t}(f)\big]\nonumber\\
=&\mathrm{Var}_{\delta_{u}}\big[\e^{-\alpha(f)t}\langle \frac{f}{t^{\gamma(f)}}-\frac{1}{\gamma(f)!}\sum_{j\in \mathfrak{I}(f)}F_{j}(f),X_{t}\rangle\big]\nonumber\\
=&\e^{-2\alpha(f)t}\int_{0}^{t}T_{s}\big(\vartheta\big[\frac{T_{t-s}f}{t^{\gamma(f)}}-\frac{1}{\gamma(f)!}\sum_{j\in \mathfrak{I}(f)}\e^{\lambda_{j}(t-s)}F_{j}(f)\big]\big)(u){ \ud s}\nonumber\\
=&\int_{0}^{t}\e^{-(2\alpha(f)-\lambda_{1})s}\e^{-\lambda_{1}s}T_{s}\big(\vartheta\big[\left(1-\frac{s}{t}\right)^{\gamma(f)}\cdot\frac{\e^{-\alpha(f)(t-s)}T_{t-s}f}{(t-s)^{\gamma(f)}}\nonumber\\
&\quad\quad-\frac{1}{\gamma(f)!}\sum_{j\in \mathfrak{I}(f)}\e^{i(t-s)\im\lambda_{j}}F_{j}(f)\big]\big)(u){ \ud s}.\label{d.2.3}
\end{align}
By Proposition \ref{prop:a}, we can easily show that for every $s>0$,
\begin{eqnarray*}
&&\left(1-\frac{s}{t}\right)^{\gamma(f)}\cdot\frac{\e^{-\alpha(f)(t-s)}T_{t-s}f}{(t-s)^{\gamma(f)}}
-\frac{1}{\gamma(f)!}\sum_{j\in \mathfrak{I}(f)}\e^{i(t-s)\im\lambda_{j}}F_{j}(f)\\
&=&\left(1-\frac{s}{t}\right)^{\gamma(f)}\left(\frac{\e^{-\alpha(f)(t-s)}T_{t-s}f}{(t-s)^{\gamma(f)}}
-\frac{1}{\gamma(f)!}\sum_{j\in \mathfrak{I}(f)}\e^{i(t-s)\im\lambda_{j}}F_{j}(f)\right)\\
&&-\left(1-\left(1-\frac{s}{t}\right)^{\gamma(f)}\right)\frac{1}{\gamma(f)!}\sum_{j\in \mathfrak{I}(f)}\e^{i(t-s)\im\lambda_{j}}F_{j}(f)\\
&\rightrightarrows& 0 \mbox{ as }t\to+\infty.
\end{eqnarray*}
Then applying similar argument (with minor modification) as in the proof of Lemma \ref{lem:newcase2}(ii), we can show that the integral in the right hand side of \eqref{d.2.3} converges to $0$ as $t\to+\infty$. Again, by Proposition \ref{prop:a}, for every $u\in E$,
$$\pp_{\delta_{u}}\big[\frac{\e^{-\alpha(f)t}}{t^{\gamma(f)}}\langle f,X_{t}\rangle-\frac{1}{\gamma(f)!}\sum_{j\in \mathfrak{I}(f)}\e^{i t\im \lambda_{j}}W^{(j)}_{t}(f)\big]
=\frac{\e^{-\alpha(f)t}}{t^{\gamma(f)}}T_{t}f(u)-\frac{1}{\gamma(f)!}\sum_{j\in \mathfrak{I}(f)}\e^{it\im \lambda_{j}}F_{j}(f)(u)$$
converges to $0$,
as $t\to+\infty$. Hence we get $\pp_{\delta_{u}}\big[\big|\frac{\e^{-\alpha(f)t}}{t^{\gamma(f)}}\langle f,X_{t}\rangle-\frac{1}{\gamma(f)!}\sum_{j\in \mathfrak{I}(f)}\e^{i t\im \lambda_{j}}W^{(j)}_{t}(f)\big|^{2}\big]$ $\to 0$ as $t\to+\infty$, which yields \eqref{d.2.2}.\qed

%%%%%%%%%%%%%%%%%%%%%%%%%%%%%%%%%%%%%%%%%%%%%%%%%%%%%%%%%%%%%%%%%%%
%%                                                               %%
%% Supplementary Material, if any, should be provided in         %%
%% {supplement} environment  with title and short description.   %%
%%                                                               %%
%%%%%%%%%%%%%%%%%%%%%%%%%%%%%%%%%%%%%%%%%%%%%%%%%%%%%%%%%%%%%%%%%%%

%\begin{supplement}
%\stitle{Appendix A}
%\sdescription{Campbell's formula with complex-valued functions for Poisson random measure.}
%\end{supplement}
%\begin{supplement}
%\stitle{Appendix B}
%\sdescription{Short description of Supplement B.}
%\end{supplement}
\appendix
\section{Appendix}
\subsection{Campbell's formula with complex-valued functions for Poisson random measure}
Suppose $(S,\mathcal{S},\eta)$ is a $\sigma$-finite measure space and $N$ is a Poisson random measure on $(S,\mathcal{S})$ with intensity $\eta$. Let $g:S\to\R$ be a measurable function. Then $\int_{S}g(x)N({ \ud x})$ is a.s. absolutely finite if and only if
\begin{equation}\label{A1}
\int_{S}1\wedge |g(x)|\eta({ \ud x})<+\infty.
\end{equation}
Moreover, when \eqref{A1} holds,
$$\mathrm{E}\left[\e^{i\theta\int_{S}g(x)N({ \ud x})}\right]=\exp\big\{\int_{S}\left(\e^{i\theta g(x)}-1\right)\eta({ \ud x})\big\}\quad\forall \theta\in\R.$$
The above result is known as the \textit{Campbell's formula}. In the following we shall give a version of the Campbell's formula with complex valued functions.

\begin{lemma}\label{lemA.1}
Suppose $(S,\mathcal{S},\eta)$ is a $\sigma$-finite measure space and $(\Omega,\mathcal{F},\mathrm{P})$ is a probability space where $N$ is a Poisson random measure on $(S,\mathcal{S})$ with intensity $\eta$.
If $F$ is a $\mathbb{C}^{+}$-valued measurable function, then $\int_{S}F(x)N({ \ud x})$ is a.s. absolutely finite if and only if
$$\int_{S}\left(1\wedge |\re F(x)|+1\wedge |\im F(x)|\right)\eta({ \ud x})<+\infty.$$
In this case,
$$\mathrm{E}\left[\e^{-\int_{S}F(x)N({ \ud x})}\right]=\exp\big\{\int_{S}\left(\e^{-F(x)}-1\right)\eta(\ud x)\big\}.$$
Here $\mathrm{E}$ denotes the expectation operator corresponding to $\mathrm{P}$.
\end{lemma}

The following lemma is needed to prove Lemma \ref{lemA.1}.
\begin{lemma}\label{lemA.2}
Suppose the assumptions of Lemma \ref{lemA.1} hold.
Let $a\ge 0$ and $A$ be an arbitrary set in $\mathcal{S}$ satisfying that $\eta(A)<+\infty$.
Define a probability measure $\mathrm{P}^{a}_{A}$ on $\Omega$ by
$$\left.\frac{d\mathrm{P}^{a}_{A}}{d\mathrm{P}}\right|_{\mathcal{F}}=\frac{\e^{-aN(A)}}{\mathrm{E}\left[\e^{-aN(A)}\right]}.$$
Then under $\mathrm{P}^{a}_{A}$, $N$ is a Poisson random measure on $(S,\mathcal{S})$ with intensity $\eta^{a}_{A}({ \ud x}):=\e^{-a1_{A}(x)}\eta({ \ud x})$.
\end{lemma}
\begin{proof}
First we recall the Laplace transform of a Poisson distributed random variable: Suppose $\xi$ is a random variable. Then $\xi$ is Poisson distributed with parameter $\delta>0$ if and only if its Laplace transform is given by
\begin{equation}\label{lemA.2.1}
\mathrm{E}\left[\e^{-\theta \xi}\right]=\exp\{(\e^{-\theta}-1)\delta\}\quad\forall \theta\ge 0.
\end{equation}
We use $\mathrm{E}^{a}_{A}$ to denote the expectation operator corresponding to $\mathrm{P}^{a}_{A}$. Take an arbitrary $B\in\mathcal{S}$. We have for every $\theta\ge 0$,
\begin{eqnarray}
\mathrm{E}^{a}_{A}\left[\e^{-\theta N(B)}\right]&=&\frac{\mathrm{E}\left[\e^{-aN(A)-\theta N(B)}\right]}{\mathrm{E}\left[\e^{-aN(A)}\right]}\nonumber\\
&=&\frac{\mathrm{E}\left[\e^{-aN(A\setminus B)-(a+\theta)N(AB)-\theta N(B\setminus A)}\right]}{\mathrm{E}\left[\e^{-a N(A)}\right]}\nonumber\\
&=&\exp\{(1-\e^{-a})\eta(A)\}\mathrm{E}\left[\e^{-a N(A\setminus B)}\right]\mathrm{E}\left[\e^{-(a+\theta)N(AB)}\right]\mathrm{E}\left[\e^{-\theta N(B\setminus A)}\right]\nonumber\\
&=&\exp\{(\e^{-\theta}-1)\e^{-a}\eta(AB)+(\e^{-\theta}-1)\eta(B\setminus A)\}\nonumber\\
&=&\exp\{(e^{-\theta}-1)\eta^{a}_{A}(B)\}.\label{lemA.2.2}
\end{eqnarray}
The third equality is from the independence of Poisson random measure on disjoint sets and the fourth equality is from \eqref{lemA.2.1}.
The above equation implies that $(N(B),\mathrm{P}^{a}_{A})$ is a Poisson distributed random variable with parameter $\eta^{a}_{A}(B)$. Let $C\in \mathcal{S}$ be a disjoint set of $B$, and $D:=S\setminus (B\cup C)$.
For $\theta_{1},\theta_{2}\ge 0$, we have
\begin{eqnarray*}
\mathrm{E}^{a}_{A}\left[\e^{-\theta_{1}N(B)-\theta_{2}N(C)}\right]&=&\frac{\mathrm{E}\left[\e^{-aN(A)-\theta_{1}N(B)-\theta_{2}N(C)}\right]}{\mathrm{E}\left[\e^{-aN(A)}\right]}\\
&=&\frac{\mathrm{E}\left[\e^{-aN(AB)-\theta_{1}N(B)}\right]\mathrm{E}\left[\e^{-aN(AC)-\theta_{2}N(C)}\right]\mathrm{E}\left[\e^{-aN(AD)}\right]}
{\mathrm{E}\left[\e^{-aN(AB)}\right]\mathrm{E}\left[\e^{-aN(AC)}\right]\mathrm{E}\left[\e^{-aN(AD)}\right]}\\
&=&\mathrm{E}^{a}_{AB}\left[\e^{-\theta_{1}N(B)}\right]\mathrm{E}^{a}_{AC}\left[\e^{-\theta_{2}N(C)}\right]\\
&=&\exp\{(\e^{-\theta_{1}}-1)\eta^{a}_{AB}(B)\}\exp\{(\e^{-\theta_{2}}-1)\eta^{a}_{AC}(C)\}\\
&=&\mathrm{E}^{a}_{A}\left[\e^{-\theta_{1}N(B)}\right]\mathrm{E}^{a}_{A}\left[\e^{-\theta_{2}N(C)}\right].
\end{eqnarray*}
The fourth equality is from \eqref{lemA.2.2} and
the final equality is from the facts that $\eta^{a}_{AB}(B)=\eta^{a}_{A}(B)$ and $\eta^{a}_{AC}(C)=\eta^{a}_{A}(C)$. Hence we prove $N(B)$ and $N(C)$ are independent under $\mathrm{P}^{a}_{A}$. Applying similar argument, one can further prove that if $A_{1},\cdot,A_{n}$ are mutually disjoint sets in $\mathcal{S}$, then $N(A_{1}),\cdots,N(A_{n})$ are independent under $\mathrm{P}^{a}_{A}$.
Therefore under $\mathrm{P}^{a}_{A}$, $N$ is a Poisson random measure with intensity $\eta^{a}_{A}$.
\end{proof}

It follows by Lemma \ref{lemA.2} and the classic Campbell's formula that if $g:\ S\to \R$ is a measurable function satisfying \eqref{A1}, then
\begin{align}
\mathrm{E}\left[\e^{-a N(A)+i\int_{S}g(x)N({ \ud x})}\right]=&\mathrm{E}\left[\e^{-a N(A)}\right]\mathrm{E}^{a}_{A}\left[\e^{i\int_{S}g(x)N({ \ud x})}\right]\nonumber\\
=&\exp\{(\e^{-a}-1)\eta(A)\}\exp\{\int_{S}\left(\e^{ig(x)}-1\right)\eta^{a}_{A}({ \ud x})\}\nonumber\\
=&\exp\{(\e^{-a}-1)\eta(A)\}\exp\{\int_{S}\left(\e^{ig(x)}-1\right)\e^{-a 1_{A}(x)}\eta({ \ud x})\}\nonumber\\
=&\exp\{\int_{S}\left(\e^{-a 1_{A}(x)+ig(x)}-1\right)\eta({ \ud x})\}.\label{A2}
\end{align}

\noindent\textbf{Proof of Lemma \ref{lemA.1}:} We note that $\int_{S}F(x)N({ \ud x})=\int_{S}\re F(x)N({ \ud x})+i\int_{S}\im F(x)N({ \ud x})$. Thus the first assertion of this lemma follows directly from the classic Campbell's formula. We only need to prove the second assertion. Let $f(x):=\re F(x)$ and $g(x)=-\im F(x)$. We have to show that for every $f:S\to \R^{+}$ and $g:S\to \R$ measurable functions such that $\int_{S}\left(1\wedge |f(x)|+1\wedge |g(x)|\right)\eta({ \ud x})<+\infty$,
\begin{equation}\label{lemA.1.1}
\mathrm{E}\left[\e^{-\int_{S}f(x)N({ \ud x})+i\int_{S}g(x)N({ \ud x})}\right]=\exp\{\int_{S}\left(\e^{-f(x)+ig(x)}-1\right)\eta({ \ud x})\}.
\end{equation}
The assumption that $f$ is a nonnegative function assures that both sides of the above equation are finite.

For $k\in\mathbb{N}$, let $\mathcal{A}_{k}$ be the set of nonnegative measurable functions on $S$ that can be represented by $\sum_{i=1}^{k}a_{i}1_{A_{i}}(x)$ where $a_{1},\cdots,a_{k}\in\R^{+}$ and $A_{1},\cdots,A_{k}$ are mutually disjoint sets in $\mathcal{S}$ with $\eta(A_{i})<+\infty$ for $1\le i\le k$.
First we claim that \eqref{lemA.1.1} holds for all functions $f\in \mathcal{A}_{k}$, $k\ge 1$. We shall prove this claim by induction.
It is easy to see form \eqref{A2} that this claim is valid for $k=1$. Now we suppose it is valid for $k=n$. Suppose $f(x)=\sum_{i=1}^{n+1}a_{i}1_{A_{i}}(x)\in \mathcal{A}_{n+1}$.
Let $f_{n}(x):=\sum_{i=1}^{n}a_{i}1_{A_{i}}(x)$. We have
\begin{align*}
\mathrm{E}\left[\e^{-\int_{S}f(x)N({ \ud x})+i\int_{S}g(x)N({ \ud x})}\right]%\\
=&\mathrm{E}\left[\e^{-a_{n+1}N(A_{n+1})-\int_{S}f_{n}(x)N({ \ud x})+i\int_{S}g(x)N({ \ud x})}\right]\\
=&\mathrm{E}\left[\e^{-a_{n+1}N(A_{n+1})}\right]\mathrm{E}^{a_{n+1}}_{A_{n+1}}\left[\e^{-\int_{S}f_{n}(x)N({ \ud x})+i\int_{S}g(x)N({ \ud x})}\right].
\end{align*}
By Lemma \ref{lemA.2}, $(N,\mathrm{P}^{a_{n+1}}_{A_{n+1}})$ is a Poisson random measure with intensity $\eta^{a_{n+1}}_{A_{n+1}}$. Applying our assumption for $k=n$ to the final expectation, we get that
\begin{eqnarray*}
&&\mathrm{E}\left[\e^{-\int_{S}f(x)N({ \ud x})+i\int_{S}g(x)N({ \ud x})}\right]\\
&=&\exp\{(\e^{-a_{n+1}}-1)\eta(A_{n+1})\}\exp\{\int_{S}\left(\e^{-f_{n}(x)+ig(x)}-1\right)\eta^{a_{n+1}}_{A_{n+1}}({ \ud x})\}\\
&=&\exp\{(\e^{-a_{n+1}}-1)\eta(A_{n+1})+\int_{S}\left(\e^{-f_{n}(x)+ig(x)}-1\right)\e^{-a_{n+1}1_{A_{n+1}}(x)}\eta({ \ud x})\}\\
&=&\exp\{\int_{S}\left(\e^{-f(x)+ig(x)}-1\right)\eta({ \ud x})\}.
\end{eqnarray*}
Hence we prove this claim for $k=n+1$.

Suppose $f$ is an arbitrary nonnegative function on $S$ satisfying that $\int_{S}1\wedge|f(x)|\eta({ \ud x})<+\infty$. One could select an increasing sequence of functions $\{f_{n}:n\ge 1\}$ from $\cup_{k\ge 1}\mathcal{A}_{k}$ such that $f_{n}$ converges pointwise to $f$. Then \eqref{lemA.1.1} is also valid for $f$ by the monotone convergence theorem. Hence we complete the proof.\qed

\subsection{An upperbound inequality for the fourth moment of superprocess}

\begin{lemma}\label{lemA.5}
Suppose
\begin{equation}\label{lemA.5.0}
\sup_{x\in E}\int_{\mfo }{\langle 1,\nu\rangle}^{4}H(x,{ \ud \nu})<+\infty.
\end{equation}
There exists a constant $c>0$ such that for every $g\in\mathcal{B}_{b}(E)$, $\tau>0$ and $x\in E$,
$$\pp_{\delta_{x}}\left[\langle g,X_{\tau}\rangle^{4}\right]\le c\left((T_{\tau}g(x))^{4}+\mathrm{Var}^{2}_{\delta_{x}}[\langle g,X_{\tau}\rangle]+\int_{0}^{\tau}T_{s}\left(\vartheta[\mathrm{Var}_{\delta_{\cdot}}[\langle g,X_{\tau-s}\rangle]]\right)(x){ \ud s}\right).$$
\end{lemma}

\begin{proof}
Fix arbitrary $g\in\mathcal{B}_{b}(E)$, $\tau>0$ and $x\in E$. The argument in Section \ref{sec:stochastic integral} tells us that $\{M^{\tau}_{t}:=\int_{0}^{t}\int_{E}T_{\tau-s}g(x)M({ \ud s},{ \ud x}):\ t\in [0,\tau]\}$ is a c\`{a}dl\`{a}g local martingale with quadratic variation
$$[M^{\tau}]_{t}=\int_{0}^{t}\langle 2b(T_{\tau-s}g)^{2},X_{s}\rangle { \ud s}+\sum_{s\le t}\langle T_{\tau-s}g,\triangle X_{s}\rangle^{2}1_{\{\triangle X_{s}\not=0\}}.$$
Moreover, by \eqref{eq:stochastic integral representation}, $\langle g,X_{\tau}\rangle$ can be represented by
$$\langle g,X_{\tau}\rangle=T_{\tau}g(x)+M^{\tau}_{\tau}\quad\pp_{\delta_{x}}\mbox{-a.s.}$$
Hence we have
\begin{equation}\label{lemA.5.1}
\pp_{\delta_{x}}\left[\langle g,X_{\tau}\rangle^{4}\right]\le 2^{3}\left((T_{\tau}g(x))^{4}+\pp_{\delta_{x}}[(M^{\tau}_{\tau})^{4}]\right).
\end{equation}
By the Burkholder-Davis-Gundy inequality, we have
\begin{align*}
\pp_{\delta_{x}}\left[|M^{\tau}_{\tau}|^{4}\right]&\le\pp_{\delta_{x}}\left[\sup_{t\in [0,\tau]}|M^{\tau}_{t}|^{4}\right]\lesssim\pp_{\delta_{x}}\left[[M^{\tau}]^{2}_{\tau}\right]\\
&\lesssim \pp_{\delta_{x}}\left[\left(\int_{0}^{\tau}\langle 2b(T_{\tau-s}g)^{2},X_{s}\rangle { \ud s}\right)^{2}\right]
+\pp_{\delta_{x}}\left[\left(\sum_{s\le \tau}\langle T_{\tau-s}g,\triangle X_{s}\rangle^{2}1_{\{\triangle X_{s}\not=0\}}\right)^{2}\right]\\
&=:I(\tau)+II(\tau).
\end{align*}
Let $f^{\tau}_{1,s}(y):=2b(y)(T_{\tau-s}g(y))^{2}$ for $s\in [0,\tau]$ and $y\in E$. Note that
$$\left(\int_{0}^{\tau}\langle f^{\tau}_{1,s},X_{s}\rangle { \ud s}\right)^{2}=\int_{0}^{\tau}\int_{0}^{\tau}\langle f^{\tau}_{1,s},X_{s}\rangle\langle f^{\tau}_{1,r},X_{r}\rangle \ud s\ud r=2\int_{0}^{\tau}\langle f^{\tau}_{1,s},X_{s}\rangle { \ud s}\int_{s}^{\tau}\langle f^{\tau}_{1,r},X_{r}\rangle { \ud r}.$$
Thus by the Fubini's theorem and Markov property we have
\begin{align}\label{lemA.4.1}
I(\tau)=&2\pp_{\delta_{x}}\left[\int_{0}^{\tau}\langle f^{\tau}_{1,s},X_{s}\rangle { \ud s}\int_{s}^{\tau}\langle f^{\tau}_{1,r},X_{r}\rangle { \ud r}\right]\nonumber\\
=&2\int_{0}^{\tau}\pp_{\delta_{x}}\left[\langle f^{\tau}_{1,s},X_{s}\rangle \int_{s}^{\tau}\langle T_{r-s}f^{\tau}_{1,r},X_{s}\rangle { \ud r}\right]{ \ud s}=2\int_{0}^{\tau}\pp_{\delta_{x}}
\left[\langle f^{\tau}_{1,s},X_{s}\rangle \langle g^{\tau}_{1,s},X_{s}\rangle\right]{ \ud s},
\end{align}
where $g^{\tau}_{1,s}(y):=\int_{s}^{\tau}T_{r-s}f^{\tau}_{1,r}(y){ \ud r}$ for $s\in [0,\tau]$ and $y\in E$.
By \eqref{eq:second moment} we have for any $f,h\in\mathcal{B}_{b}(E)$ and $t\ge 0$,
$$\pp_{\delta_{x}}\left[\langle f,X_{t}\rangle\langle h,X_{t}\rangle\right]=T_{t}f(x)\cdot T_{t}h(x)+\int_{0}^{t}T_{s}\left(\vartheta[T_{t-s}f,T_{t-s}g]\right)(x){ \ud s}.$$
Hence we get from \eqref{lemA.4.1} that
\begin{align*}
I(\tau)=&2\left(\int_{0}^{\tau}T_{s}f^{\tau}_{1,s}(x)T_{s}g^{\tau}_{1,s}(x){ \ud s}+\int_{0}^{\tau}{ \ud s}\int_{0}^{s}T_{r}\left(\vartheta[T_{s-r}f^{\tau}_{1,s},T_{s-r}g^{\tau}_{1,s}]\right)(x){ \ud r}\right)\\
=:&2(I_{1}(\tau)+I_{2}(\tau)).
\end{align*}
By Fubini's theorem,  $T_{s-r}g^{\tau}_{1,s}=T_{s-r}\left(\int_{s}^{\tau}T_{u-s}f^{\tau}_{1,u}{ \ud u}\right)=\int_{s}^{\tau}T_{u-r}f^{\tau}_{1,u}{ \ud u}$. In particular, $T_{s}g^{\tau}_{1,s}=\int_{s}^{\tau}T_{u}f^{\tau}_{1,u}{ \ud u}$. Thus we have
\begin{align}\label{lemA.4.2}
I_{1}(\tau)
=&\int_{0}^{\tau}T_{s}f^{\tau}_{1,s}(x){ \ud s}\int_{s}^{\tau}T_{u}f^{\tau}_{1,u}(x){ \ud u}
=\frac{1}{2}\left(\int_{0}^{\tau}T_{s}f^{\tau}_{1,s}(x){ \ud s}\right)^{2}\nonumber\\
=&\frac{1}{2}\left(\int_{0}^{\tau}T_{s}\left(2b(T_{\tau-s}g)^{2}\right)(x){ \ud s}\right)^{2}\le\frac{1}{2}\left(\int_{0}^{\tau}T_{s}(\vartheta[T_{\tau-s}g])(x){ \ud s}\right)^{2}=\frac{1}{2}\mathrm{Var}^{2}_{\delta_{x}}[\langle g,X_{\tau}\rangle].
\end{align}
On the other hand, by the definition of $\vartheta[\cdot,\cdot]$ and Fubini's theorem,
\begin{align*}
\int_{0}^{s}T_{r}\left(\vartheta[T_{s-r}f^{\tau}_{1,s},T_{s-r}g^{\tau}_{1,s}]\right)(x){ \ud r}=&\int_{0}^{s}T_{r}\left(\vartheta[T_{s-r}f^{\tau}_{1,s},\int_{s}^{\tau}T_{u-r}f^{\tau}_{1,u}{ \ud u}]\right)(x){ \ud r}\\
=&\int_{0}^{s}{ \ud r}\int_{s}^{\tau}T_{r}\left(\vartheta[T_{s-r}f^{\tau}_{1,s},T_{u-r}f^{\tau}_{1,u}]\right)(x){ \ud u}.
\end{align*}
Then we have
\begin{align}
I_{2}(\tau)=&\int_{0}^{\tau}{ \ud s}\int_{0}^{s}{ \ud r}\int_{s}^{\tau}T_{r}\left(\vartheta[T_{s-r}f^{\tau}_{1,s},T_{u-r}f^{\tau}_{1,u}]\right)(x){ \ud u}\nonumber\\
=&\int_{0}^{\tau}T_{r}\left(\int_{r}^{\tau}{ \ud s}\int_{s}^{\tau}\vartheta[T_{s-r}f^{\tau}_{1,s},T_{u-r}f^{\tau}_{1,u}]{ \ud u}\right)(x){ \ud r}\nonumber\\
=&\frac{1}{2}\int_{0}^{\tau}T_{r}\left(\vartheta\left[\int_{r}^{\tau}T_{s-r}f^{\tau}_{1,s}{ \ud s}\right]\right)(x){ \ud r}=\frac{1}{2}\int_{0}^{\tau}T_{r}\left(\vartheta\left[g^{\tau}_{1,r}\right]\right)(x){ \ud r}.\label{lemA.4.3}
\end{align}
In the third equality we use the fact that for any $h_{s}\in \mathcal{B}^{+}_{b}(E)$,
$$2\int_{r}^{\tau}{ \ud s}\int_{s}^{\tau}\vartheta[h_{s},h_{u}]{ \ud u}=\int_{r}^{\tau}{ \ud s}\int_{r}^{\tau}\vartheta[h_{s},h_{u}]{ \ud u}=\vartheta[\int_{r}^{\tau}h_{s}{ \ud s},\int_{r}^{\tau}h_{u}{ \ud u}]
=\vartheta[\int_{r}^{\tau}h_{s}{ \ud s}].$$
We note that $g^{\tau}_{1,r}(y)=\int_{0}^{\tau-r}T_{s}f^{\tau}_{s+r}(y){ \ud s}=\int_{0}^{\tau-r}T_{s}(2b(T_{\tau-r-s}g)^{2})(y){ \ud s}\le \mathrm{Var}_{\delta_{y}}[\langle g,X_{\tau-r}\rangle]$ for $y\in E$. Hence by \eqref{lemA.4.3} we get
\begin{equation}\label{lemA.5.2}
I_{2}(\tau)\le \frac{1}{2}\int_{0}^{\tau}T_{r}\left(\vartheta[\mathrm{Var}_{\delta_{\cdot}}[\langle g,X_{\tau-r}\rangle]]\right)(x){ \ud r}.
\end{equation}
Next we deal with $II(\tau)$. We note that
\begin{align}
II(\tau)=&\pp_{\delta_{x}}\left[\left(\sum_{s\le \tau}{\langle T_{\tau-s}g,\triangle X_{s}\rangle^{2}}\right)^{2}\right]\nonumber\\
=&2\pp_{\delta_{x}}\left[\sum_{s\le \tau}{\langle T_{\tau-s}g,\triangle X_{s}\rangle^{2}}\sum_{s\le r\le\tau}
{\langle T_{\tau-r}g,\triangle X_{r}\rangle^{2}}\right]-\pp_{\delta_{x}}\left[\sum_{s\le \tau}{\langle T_{\tau-s}g,\triangle X_{s}\rangle^{4}}\right]\nonumber\\
=&2\pp_{\delta_{x}}\left[\sum_{s\le \tau}{\langle T_{\tau-s}g,\triangle X_{s}\rangle^{2}}\sum_{s\le r\le\tau}{\langle T_{\tau-r}g,\triangle X_{r}\rangle^{2}}\right]\nonumber\\
&-\pp_{\delta_{x}}\left[\int_{0}^{\tau}\langle\int_{\mfo }{\langle T_{\tau-s}g,\nu\rangle}^{4}H(\cdot,{ \ud \nu}),X_{s}\rangle { \ud s}\right].\label{lemA.4.8}
\end{align}
We remark here that the final term of \eqref{lemA.4.8} is finite under our assumption \eqref{lemA.5.0}.
By the Markov property, the first term in the right hand side equals
\begin{align}
&2\pp_{\delta_{x}}\left[\sum_{s\le \tau}{\langle T_{\tau-s}g,\triangle X_{s}\rangle^{2}}1_{\{\triangle X_{s}\not=0\}}\pp_{X_{s}}\left[\sum_{r\le\tau-s}
{\langle T_{\tau-s-r}g,\triangle X_{r}\rangle^{2}}1_{\{\triangle X_{r}\not=0\}}\right]\right]\nonumber\\
=&2\pp_{\delta_{x}}\left[\sum_{s\le \tau}{\langle T_{\tau-s}g,\triangle X_{s}\rangle^{2}}1_{\{\triangle X_{s}\not=0\}}\pp_{X_{s}}\left[\int_{0}^{\tau-s}\langle\int_{\mfo }{\langle T_{\tau-s-r}g,\nu\rangle}^{2}H(\cdot,{ \ud \nu}),X_{r-}\rangle { \ud r}\right]\right].\label{lemA.4.7}
\end{align}
 Define $f^{\tau}_{2,s}(y):=\int_{\mfo }{\langle T_{\tau-s}g,\nu\rangle}^{2}H(y,{ \ud \nu})$ and $g^{\tau}_{2,s}(y):=\int_{0}^{\tau-s}T_{r}f^{\tau}_{2,s+r}(y){ \ud r}=\int_{s}^{\tau}T_{r-s}f^{\tau}_{2,r}(y){ \ud r}$ for $y\in E$ and $s\in [0,\tau]$. We can continue the above computation to get
\begin{align*}
&\mbox{RHS of }\eqref{lemA.4.7}\\
=&2\pp_{\delta_{x}}\left[\sum_{s\le \tau}{\langle T_{\tau-s}g,\triangle X_{s}\rangle^{2}}1_{\{\triangle X_{s}\not=0\}}\pp_{X_{s}}\left[\int_{0}^{\tau-s}\langle
f^{\tau}_{2,s+r},X_{r-}\rangle { \ud r}\right]\right]\\
=&2\pp_{\delta_{x}}\left[\sum_{s\le \tau}{\langle T_{\tau-s}g,\triangle X_{s}\rangle^{2}}1_{\{\triangle X_{s}\not=0\}}\langle g^{\tau}_{2,s},X_{s}\rangle\right]\\
=&2\pp_{\delta_{x}}\left[\sum_{s\le \tau}\langle g^{\tau}_{2,s},X_{s-}\rangle{\langle T_{\tau-s}g,\triangle X_{s}\rangle^{2}}1_{\{\triangle X_{s}\not=0\}}\right]+2\pp_{\delta_{x}}\left[\sum_{s\le \tau}{\langle g^{\tau}_{2,s},\triangle X_{s}\rangle}{\langle T_{\tau-s}g,\triangle X_{s}\rangle^{2}}1_{\{\triangle X_{s}\not=0\}}\right]\\
=&2\pp_{\delta_{x}}\left[\int_{0}^{\tau}\langle f^{\tau}_{2,s},X_{s}\rangle\langle g^{\tau}_{2,s},X_{s}\rangle { \ud s}\right]
+2\pp_{\delta_{x}}\left[\int_{0}^{\tau}\langle\int_{\mfo }{\langle T_{\tau-s}g,\nu\rangle}^{2}{\langle g^{\tau}_{2,s},\nu\rangle}H(\cdot,{ \ud \nu}),X_{s} \rangle { \ud s}\right]\\
\le&2\pp_{\delta_{x}}\left[\int_{0}^{\tau}\langle f^{\tau}_{2,s},X_{s}\rangle\langle g^{\tau}_{2,s},X_{s}\rangle { \ud s}\right]+\pp_{\delta_{x}}\left[\int_{0}^{\tau}\langle\int_{\mfo }{\langle T_{\tau-s}g,\nu\rangle}^{4}+{\langle g^{\tau}_{2,s},\nu\rangle}^{2}H(\cdot,{ \ud \nu}),X_{s} \rangle { \ud s}\right].
\end{align*}
Inserting this to the right hand side of \eqref{lemA.4.8} yields that
\begin{align*}
II(\tau)\le&2\pp_{\delta_{x}}\left[\int_{0}^{\tau}\langle f^{\tau}_{2,s},X_{s}\rangle\langle g^{\tau}_{2,s},X_{s}\rangle { \ud s}\right]+\pp_{\delta_{x}}\left[\int_{0}^{\tau}\langle\int_{\mathcal{M}^{0}(E)}{\langle g^{\tau}_{2,s},\nu\rangle}^{2}H(\cdot,{ \ud \nu}),X_{s}\rangle { \ud s}\right]\\
\le&2\int_{0}^{\tau}\pp_{\delta_{x}}\left[\langle f^{\tau}_{2,s},X_{s}\rangle\langle g^{\tau}_{2,s},X_{s}\rangle\right]{ \ud s}+\int_{0}^{\tau}T_{s}\left(\int_{\mfo }{\langle g^{\tau}_{2,s},\nu\rangle}^{2}H(\cdot,{ \ud \nu})\right)(x){ \ud s}\\
\le&2\int_{0}^{\tau}\pp_{\delta_{x}}\left[\langle f^{\tau}_{2,s},X_{s}\rangle\langle g^{\tau}_{2,s},X_{s}\rangle\right]{ \ud s}+\int_{0}^{\tau}T_{s}\left(\vartheta[g^{\tau}_{2,s}]\right)(x){ \ud s}\\
=&2\int_{0}^{\tau}T_{s}f^{\tau}_{2,s}(x)T_{s}g^{\tau}_{2,s}(x){ \ud s}+2\int_{0}^{\tau}{ \ud s}\int_{0}^{s}T_{s-r}\left(\vartheta[T_{r}f^{\tau}_{2,s},T_{r}g^{\tau}_{2,s}]\right)(x){ \ud r}
+\int_{0}^{\tau}T_{s}(\vartheta[g^{\tau}_{2,s}])(x){ \ud s}\\
=:&2II_{1}(\tau)+2II_{2}(\tau)+III_{3}(\tau).
\end{align*}
We note that for $s\in [0,\tau]$ and $y\in E$, $f^{\tau}_{2,s}(y)\le\vartheta[T_{\tau-s}g](y)$, and that
\begin{align*}
g^{\tau}_{2,s}(y)=&\int_{0}^{\tau-s}T_{r}\left(\int_{\mfo }{\langle T_{\tau-s-r}g,\nu\rangle}^{2}H(\cdot,{ \ud \nu})\right)(y){ \ud s}\\
\le& \int_{0}^{\tau-s}T_{r}(\vartheta[T_{\tau-s-r}g])(y){ \ud r}=\mathrm{Var}_{\delta_{y}}[\langle g,X_{\tau-s}\rangle].
\end{align*}
Similarly as in \eqref{lemA.4.2}, we can show that
\begin{equation}\label{lemA.5.3}
II_{1}(\tau)\le \frac{1}{2}\mathrm{Var}^{2}_{\delta_{x}}[\langle g,X_{\tau}\rangle].
\end{equation}
Similarly as in \eqref{lemA.4.3} and \eqref{lemA.5.2} we can show that
\begin{equation}\label{lemA.5.4}
II_{2}(\tau)=\frac{1}{2}III_{3}(\tau)\le\frac{1}{2}\int_{0}^{\tau}T_{s}\left(\vartheta[\mathrm{Var}_{\delta_{\cdot}}[\langle g,X_{\tau-s}\rangle]]\right)(x){ \ud s}.
\end{equation}
Combining \eqref{lemA.4.2}, \eqref{lemA.5.2}, \eqref{lemA.5.3} and \eqref{lemA.5.4}, we get that
$$\pp_{\delta_{x}}\left[(M^{\tau}_{\tau})^{4}\right]\le c_{1}\left(\mathrm{Var}^{2}_{\delta_{x}}[\langle g,X_{\tau}\rangle]+\int_{0}^{\tau}T_{s}\left(\vartheta[\mathrm{Var}_{\delta_{\cdot}}[\langle g,X_{\tau-s}\rangle]]\right)(x){ \ud s}\right)$$
for some constant $c_{1}>0$ independent of $g$, $x$ and $\tau$.
This together with \eqref{lemA.5.1} yields the desired result.
\end{proof}

\begin{lemma}\label{lemA.4}
Suppose \eqref{lemnew.1} holds for $\epsilon=\lambda_{1}/2$.
Then for every $x\in E$,
there exists $t_{0}>0$ such that
$$\sup_{t\ge t_{0}}\frac{\e^{-2\lambda_{1} t}}{t^{2+4r}}\pp_{\delta_{x}}\left[\langle g,X_{t}\rangle^{4}\right]<+\infty.$$
\end{lemma}

\begin{proof}
By Lemma \ref{lemA.5}, there is some constant $c_{1}>0$, such that,
$$\pp_{\delta_{x}}\left[\langle g,X_{t}\rangle^{4}\right]\le c_{1}\left((T_{t}g(x))^{4}+\mathrm{Var}^{2}_{\delta_{x}}[\langle g,X_{t}\rangle]+\int_{0}^{t}T_{s}\left(\vartheta[\mathrm{Var}_{\delta_{\cdot}}[\langle g,X_{t-s}\rangle]]\right)(x){ \ud s}\right)$$
for all $t>0$ and $x\in E$. Our assumption on $g$ implies that
\begin{equation}\label{LemA.4.1}
\frac{\e^{-2\lambda_{1} t}}{t^{2+4r}}(T_{t}g(x))^{4}\to 0\mbox{ as }t\to+\infty.
\end{equation}
Lemma \ref{lem:newcase2}(i) implies that there is $t_{1}\ge 1$ and $0<M<+\infty$ such that
\begin{equation}\label{lemA.4.6}
\frac{\e^{-\lambda_{1} t}}{t^{1+2r}}\mathrm{Var}_{\delta_{y}}[\langle g,X_{t}\rangle]\le M\quad \forall t\ge t_{1},\ y\in E.
\end{equation}
It thus follows by \eqref{eq1.9} that $T_{r}\left(\vartheta[\mathrm{Var}_{\delta_{\cdot}}[\langle g,X_{ t-r}\rangle]]\right)(x)\lesssim \e^{\lambda_{1} r}\varphi(x)\|\vartheta[\mathrm{Var}_{\delta_{\cdot}}[\langle g,X_{ t-r}\rangle]]\|_{\infty}$. Hence we have
\begin{align*}
\int_{0}^{t}T_{s}\left(\vartheta[\mathrm{Var}_{\delta_{\cdot}}[\langle g,X_{t-s}\rangle]]\right)(x){ \ud s}&\lesssim \int_{0}^{ t}\e^{\lambda_{1} r}\varphi(x)\|\vartheta[\mathrm{Var}_{\delta_{\cdot}}[\langle g,X_{ t-r}\rangle]]\|_{\infty}{ \ud r}\\
&\le \varphi(x)\|\vartheta[1]\|_{\infty}\int_{0}^{ t}\e^{\lambda_{1} r}\|\mathrm{Var}_{\delta_{\cdot}}[\langle g,X_{ t-r}\rangle]\|^{2}_{\infty}{ \ud r},
\end{align*}
and consequently,
\begin{equation}\label{lemA.4.4}
\frac{\e^{-2\lambda_{1} t}}{ t^{2+4r}}\int_{0}^{t}T_{s}\left(\vartheta[\mathrm{Var}_{\delta_{\cdot}}[\langle g,X_{t-s}\rangle]]\right)(x){ \ud s}\lesssim \varphi(x)\frac{\e^{-\lambda_{1} t}}{ t^{2+4r}}\int_{0}^{ t}\e^{-\lambda_{1} u}\|\mathrm{Var}_{\delta_{\cdot}}[\langle g,X_{u}\rangle]\|^{2}_{\infty}{ \ud u}.
\end{equation}
By \eqref{lemA.4.6}, for $ t\ge t_{1}$,
\begin{equation}\label{lemA.4.5}
\int_{t_{1}}^{ t}\e^{-\lambda_{1} u}\|\mathrm{Var}_{\delta_{\cdot}}[\langle g,X_{u}\rangle]\|^{2}_{\infty}{ \ud u}=
\int_{t_{1}}^{ t}\e^{\lambda_{1} u}u^{2+4r}\left\|\frac{e^{-\lambda_{1} u}}{u^{1+2r}}\mathrm{Var}_{\delta_{\cdot}}[\langle g,X_{u}\rangle]\right\|^{2}_{\infty}{ \ud u}\lesssim M^{2} t^{2+4r}\e^{\lambda_{1} t}.
\end{equation}
On the other hand, by \eqref{ieq:many-to-one} we have for $u\ge 0$ and $y\in E$,
\begin{align*}
\mathrm{Var}_{\delta_{y}}[\langle g,X_{u}\rangle]=\int_{0}^{u}T_{u-s}(\vartheta[T_{s}g])(y){ \ud s}\le& \int_{0}^{u}\e^{c_{1}(u-s)}\|\vartheta[T_{s}g]\|_{\infty}{ \ud s}\\
\le&\int_{0}^{u}\e^{c_{1}(u+s)}\|g\|^{2}_{\infty}\|\vartheta[1]\|_{\infty}{ \ud s}\lesssim \e^{2c_{1}u},
\end{align*}
and thus $\int_{0}^{t_{1}}\e^{-\lambda_{1} u}\|\mathrm{Var}_{\delta_{\cdot}}[\langle g,X_{u}\rangle]\|^{2}_{\infty}{ \ud u}\lesssim \e^{2c_{1}t_{1}}.$
This together with \eqref{lemA.4.4} and \eqref{lemA.4.5} yields that
$$\sup_{ t\ge t_{1}}\frac{\e^{-2\lambda_{1} t}}{ t^{2+4r}}\int_{0}^{t}T_{s}\left(\vartheta[\mathrm{Var}_{\delta_{\cdot}}[\langle g,X_{t-s}\rangle]]\right)(x){ \ud s}\lesssim \varphi(x)\left(M^{2}+\e^{2c_{1}t_{1}}\sup_{ t\ge t_{1}}\frac{\e^{-\lambda_{1} t}}{ t^{2+4r}}\right)<+\infty.$$
This together with \eqref{LemA.4.1} and \eqref{lemA.4.6} yields the desired result.
\end{proof}

\end{document}